\documentclass[11pt]{extarticle}
\usepackage{graphicx} 
\usepackage{subcaption} 
\usepackage[labelformat=parens,labelsep=quad,skip=3pt]{caption} 
\usepackage{amsmath} 
\usepackage{amsfonts} 
\usepackage{upgreek} 
\usepackage{textcomp} 
\usepackage{natbib} 
\usepackage[a4paper, total={170mm,257mm}, left=20mm, right=20mm, top=20mm, bottom=20mm]{geometry}
\begin{document}
	\begin{titlepage}
		\begin{center}
			\baselineskip=8pt
			{\LARGE Spatially-optimized finite-difference schemes for numerical dispersion suppression in seismic applications}
			
			\vspace*{2.5cm}
			{\large Edward Caunt \\
				\vspace*{0.5cm}
				Email: ec815@ic.ac.uk} \\
			\vspace*{0.5cm}
			{\small Supervised by Dr Gerard Gorman and Dr Rhodri Nelson}\\
			\vspace*{0.5cm}
			{\small 03/01/19}\\
			\vspace*{5cm}
			
			\centerline{\includegraphics[width=3.0cm]{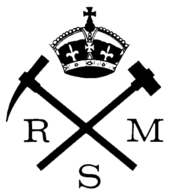}}
			\vspace{2.5cm}
			Department of Earth Science and Engineering \\
			Imperial College London \\
			\vspace*{2.5cm}
			A report submitted in partial fulfilment of the requirements for the degree of Geophysics MSci at Imperial College London and Associateship of the Royal School of Mines. It is substantially the result of my own work except where explicitly indicated in the text. The report may be freely copied and distributed provided the source is explicitly acknowledged.
			\vspace*{0.5cm}
			\\This MSci project was conducted between September 2018 and January 2019
			\\The main body of this dissertation has 5498 words
		\end{center}
		
		\clearpage
	\end{titlepage}
	\pagenumbering{roman}
	\begin{abstract}
		\noindent Propagation characteristics of a wave are defined by the dispersion relationship, from which the governing partial differential equation (PDE) can be recovered. PDEs are commonly solved numerically using the finite-difference (FD) method, discretizing the field to a regular grid and evaluating derivatives using values at a position and adjacent points. FD stencils are conventionally constructed from truncated Taylor series expansions which, whilst typically providing good approximation of the PDE in the space-time domain, often differ considerably from the original partial differential in the wavenumber-frequency domain where the dispersion relationship is defined. Consequentially, stable, high-order FD schemes may not necessarily result in realistic wave behavior, commonly exhibiting numerical dispersion: lagging high-frequency components as a product of discretization. A method for optimizing FD stencil weightings via constrained minimization to better approximate the partial derivative in the wavenumber domain is proposed, allowing for accurate propagation with coarser grids than would be otherwise possible. This was applied to second derivatives on a standard grid and first derivatives on a staggered grid. To evaluate the efficacy of the method, a pair of numerical simulations were devised to compare spatially-optimized stencils with conventional formulations of equivalent extent. A spatially-optimized formulation of the 1D acoustic wave equation with Dirichlet boundary conditions is presented, evaluating performance at a range of grid spacings, examining the interval between the theoretical maximum grid spacings for the conventional and optimized schemes in finer detail. Deviation from the analytical solution was calculated for both schemes. The optimized scheme was found to offer superior performance for undersampled wavefields with minimum wavelengths of $4.4\Delta x$ or smaller and heavily oversampled wavefields with minimum wavelengths of $9.5\Delta x$ or greater. Staggered-grid first derivative stencils were then applied to the P-SV elastic wave formulation, simulating seismic wave propagation for a two-layer, water-over-rock model. For staggered conventional and optimized schemes of 6\textsuperscript{th} and 4\textsuperscript{th} order respectively, general performance was comparable although the optimized scheme provided slightly improved numerical dispersion suppression, particularly apparent for S-wave phases. Performance benefits of the proposed method for staggered grids appear to increase with stencil extent, and are comparable to more complex existing methods, implying significant benefits for high-order FD schemes typically used for seismic processing.
	\end{abstract}
	\newpage
	\tableofcontents
	\newpage
	\pagenumbering{arabic}

	\section{Introduction}
	\subsection{An introduction to finite-difference methods}
	Finite-difference (FD) methods for solving partial differential equations (PDEs) form the cornerstone of much of seismic processing, in applications including seismic migration (e.g. \citealp{Mufti1996}, \citealp{Yan2013}, \citealp{Yan2013a}) and wave modeling (e.g. \citealp{Virieux1986}, \citealp{Tessmer2000}). FD methods benefit from simple implementation, reduced memory requirements, and small computational cost compared to other numerical methods (\citealp{Liu2009}); usage is consequentially widespread. FD methods are not limited to simple domains and can accommodate complex boundary conditions (\citealp{Lele1992}). For time marching schemes, temporal derivatives are typically evaluated using a 2\textsuperscript{nd} order FD scheme for reasons of effectiveness and stability (\citealp{Liu2009}). However, this limits modeling accuracy. Smaller temporal and spatial increments offer a straightforward way to counter this, but incur additional operations, increasing run time (\citealp{Liu2009}). Methods developed to address this without considerably increasing computational expense include high-spatial-order, staggered-grid, and implicit methods.
	\\
	\\The degree of Taylor series truncation used in the derivation of an FD scheme is typically used to express its accuracy, with an error term proportional to the dimensional increment raised to the appropriate power. With high-spatial-order schemes, the error in spatial derivatives is reduced, allowing for coarser grid spacings without introduction of large errors, thereby offsetting inherent costs of larger operators (\citealp{Liu2009}). Practical benefits include reduced memory requirements as wavefields do not require such heavy oversampling (\citealp{Dablain1986}). It is trivial to derive high-spatial-order FD schemes of any order (\citealp{Liu2009}) and their usage in seismic modeling is widespread (e.g. \citealp{Dablain1986}, \citealp{Yan2014}, \citealp{Liu2009b}). Additionally, stencils of sufficiently high order can be truncated to remove small coefficients without notably compromising accuracy (\citealp{Liu2009b}), reducing computational cost. High-order stencils are commonly usd in conjunction with staggered grids. Introduced for seismic modeling by \citealp{Madariaga1976}, staggered-grid methods evaluate derivatives on spatially offset grids, improving precision and stability of the scheme, and have seen widespread use (e.g. \citealp{Virieux1986}, \citealp{Wang2014}, \citealp{Etemadsaeed2016}). However, instabilities may arise at sharp parameter contrasts: a relatively common occurrence in seismic models (\citealp{Liu2009}). Improved stability is achievable with implicit schemes, evaluating spatial using function and derivative values at neighboring points along with function value at the current position (e.g. \citealp{Liu2009a}, \citealp{Lele1992}).
	\\
	\\A stable, high-order FD scheme will adequately approximate the governing PDE for many applications. However, for wave problems this does not guarantee an acceptable solution. Propagation characteristics of a wave are encoded in frequency-wavenumber space as the dispersion relation (\citealp{Tam1993}, \citealp{Whitham1974}), from which the governing PDE can be reconstructed. The dispersion relation takes the form of some function of frequency and wavenumber, typically obtained through spatial and temporal Fourier transforms of the governing PDE. This relationship governs group and phase velocities, along with dispersion characteristics of waves propagating within a medium. It is thus apparent a wave may propagate in an unrealistic manner, even if the FD scheme resembles the governing PDE in the time-space domain (\citealp{Tam1993}). This phenomena is a product of the discretization, referred to as numerical dispersion.

	\subsection{Numerical dispersion and its suppression}
	Numerical dispersion results in unrealistic or spurious dependence of phase velocity on frequency and propagation direction (occasionally referred to as numerical anisotropy) (\citealp{Juntunen2000}), arising from wavefield discretization, causing phase velocity to become a function of grid spacing (\citealp{Kim1990}, \citealp{Yang2003}). This causes high-frequency wave components to lag behind the group, whilst directional dependence causes squaring of circular wavefronts, the net effect being significant distortion of the wave profile. This results in greater dispersion than otherwise expected in dispersive systems, and is particularly egregious in non-dispersive systems. Methods proposed to suppress numerical dispersion include introduction of artificial anisotropy to increase propagation speed of lagging wave components (\citealp{Juntunen2000}), flux-corrected transport algorithms (\citealp{Fei1995}), and nearly analytic discrete methods (\citealp{Yang2003}, \citealp{Yang2010}).
	\\
	\\An alternative approach is to construct an FD scheme such that it has near-identical dispersion characteristics to the PDE which it approximates (e.g. \citealp{Tam1993}). Note that in accordance with Nyquist-Shannon sampling theorem, four grid spacings are required to define a sine wave without ambiguity. Thus the above criteria need only be the case for wave components with wavenumber $\alpha$ such that $\alpha\Delta x \leq \frac{\pi}{2}$ where $\Delta x$ is grid spacing in the $x$ direction. Thus an FD stencil whose Fourier transform has minimal misfit with that of the original partial derivative within the window $-\frac{\pi}{2}\leq\alpha\Delta x\leq\frac{\pi}{2}$ is desired. Such schemes are referred to as dispersion-relation-preserving (DRP) schemes (\citealp{Tam1993}). 
	\\
	\\In this work, a 4\textsuperscript{th} order DRP 1D acoustic wave formulation with zero Dirichlet boundary conditions is presented, extending methods of \citealp{Tam1993}. Performance is compared to the conventional scheme of equal extent, evaluated against the analytic solution at a range of grid spacings. In doing so, the efficacy of spatial optimization is assured and the conditions under which it is effective can be determined. Staggered optimized first-derivative stencils are then derived and applied to a 2D elastic wave test case using the P-SV formulation of \citealp{Virieux1986}, aiming to demonstrate the versatility and universal benefits of spatial optimization.

	\section{Methodology}
	\subsection{Derivation of stencil coefficients}
	Suppose we wish to evaluate the FD approximation of a spatial derivative of the discretized function $u$ at position $l$ on a uniform grid. Stencil extent of $M$ and $N$ nodes from $l$ in the positive and negative directions respectively must be specified. Thus the FD approximation for the second derivative is expressed as
	\begin{center}	
		\begin{equation}
		\frac{\partial^{2} u}{\partial x^{2}} \approx \frac{1}{(\Delta x)^{2}}\sum_{j=-N}^{M}a_{j}u_{l+j}
		\label{eq:2}
		\end{equation}
	\end{center}
	where $a_{j}$ is the stencil coefficient at point $u_{l+j}$.
 	\\
 	\\Values of $a_{j}$ are normally obtained by taking Taylor series expansions of terms on the right hand side of the equation, weighted to cancel undesired derivatives given a symmetric or antisymmetric stencil. Finite difference stencils derived in this manner will henceforth be referred to as conventional schemes. Such schemes will typically only acceptably approximate the partial derivative in the time-space domain, deviating significantly in frequency-wavenumber space. Alternatively, stencil coefficients can be specified by requiring that the Fourier transform of the stencil closely approximates that of the partial derivative. Resemblance between dispersion relationships of the FD scheme and governing PDE is thereby assured.
 	\\
 	\\Equation (\ref{eq:2}) is a special case of the following, trivially recoverable by specifying $x=l\Delta x$.
 	\begin{center}	
 		\begin{equation}
 		\frac{\partial^{2} u}{\partial x^{2}} \approx \frac{1}{(\Delta x)^{2}}\sum_{j=-N}^{M}a_{j}u(x+j\Delta x)
 		\label{eq:3}
 		\end{equation}
 	\end{center}
 	Furthermore, the Fourier transform for a function and its inverse are defined as
 	\begin{center}	
 		\begin{equation}
 		\tilde{u}(\alpha)=\frac{1}{2\pi}\int_{-\infty}^{\infty}u(x)e^{-i\alpha x}\mathrm{d}x
 		\label{eq:4}
 		\end{equation}
 	\end{center}
	\begin{center}	
		\begin{equation}
		u(x)=\int_{-\infty}^{\infty}\tilde{u}(\alpha)e^{i\alpha x}\mathrm{d}\alpha
		\label{eq:5}
		\end{equation}
	\end{center}
 	where $\alpha$ represents the wavenumber. Fourier transforming both sides of (\ref{eq:3}) yields
 	\begin{center}	
 		\begin{equation}
 		-\alpha^{2}\tilde{u}(\alpha) \approx \left(\frac{1}{(\Delta x)^{2}}\sum_{j=-N}^{M}a_{j}e^{i\alpha j\Delta x}\right)\tilde{u}(\alpha)
 		\label{eq:6}
 		\end{equation}
 	\end{center}
 	By inspection of (\ref{eq:6}), it is apparent that the wavenumber of the Fourier transform of the FD stencil can be expressed as
 	\begin{center}	
 		\begin{equation}
 		\bar{\alpha}^{2} = \frac{1}{(\Delta x)^{2}}\sum_{j=-N}^{M}a_{j}e^{i\alpha j\Delta x}
 		\label{eq:7}
 		\end{equation}
 	\end{center}
 	It is also apparent that $(\bar{\alpha}\Delta x)^{2}$ is a function of $\alpha\Delta x$ with  period $2\pi$. To achieve minimal misfit between the Fourier transforms of the partial derivative and FD scheme, stencil weighting are chosen such that the integrated error function
 	\begin{center}
 		\begin{equation}
	 		\begin{split}
	 		E = \int_{-\frac{\pi}{2}}^{\frac{\pi}{2}}\left|(\alpha\Delta x)^{2} - (\bar{\alpha}\Delta x)^{2}\right|^{2}\mathrm{d}(\alpha\Delta x) \\ \\ = \int_{-\frac{\pi}{2}}^{\frac{\pi}{2}}\left|\kappa^{2}+\sum_{j=-N}^{M}a_{j}e^{ij\kappa}\right|^{2}\mathrm{d}\kappa
	 		\label{eq:8}
	 		\end{split}
 		\end{equation}	
 	\end{center}
 	is minimized. The approximation need only be valid for waves with wavelengths greater than $4\Delta x$, corresponding with $|\alpha\Delta x| < \frac{\pi}{2}$. $E$ is minimized when conditions
 	\begin{center}	
 		\begin{equation}
 		\frac{\partial E}{\partial a_{j}} = 0; \;-N \leq j \leq M
 		\label{eq:9}
 		\end{equation}
 	\end{center}
	are met. Note that whilst \citealp{Tam1993} outlined optimization of first derivatives, error functions for stencil optimization differ extensively with derivative-order and grid, requiring unique analysis for each derivative type.
	\\
	\\For asymmetric stencils ($N\neq M$), $\bar{\alpha}$ is complex. Applying such stencils throughout the computational region will typically result in instabilities, with solutions prone to blowup (\citealp{Tam1993}). Consequentially, central-difference stencils are used exclusively within the main body of the computational domain in this work. Boundaries necessitate limited use of asymmetric stencils, although accumulated numerical instability in these regions is minimal and unlikely to cause issues.
	\\
	\\This minimization is combined with the truncated Taylor series method to create an optimized FD scheme. Specifying stencil extent as $N=M=3$ and 4\textsuperscript{th} order accuracy, a symmetric stencil can be formulated from pairs of Taylor series expansions at appropriate points. Values of $a_{j}$ are constrained, becoming functions of a single coefficient. Using $a_{1}$ as our free parameter, we obtain
	\begin{center}
		\begin{equation}
			\begin{split}                      
			a_{0}=-\frac{4a_{1}}{3}-\frac{13}{18}
			\\ \\a_{2}=a_{-2}=\frac{9}{20}-\frac{2a_{1}}{5}
			\\ \\a_{3}=a_{-3}=\frac{3a_{1}}{45}-\frac{4}{45} 
			\end{split}
		\end{equation}	
	\end{center}
	Constrained minimization of $E$ via Gauss-Newton or comparable methods, yields weights for all stencil points.
	\\Thus
	\begin{center}
		\begin{equation}
			\begin{split}                      
			a_{0}=-2.81299833
			\\ \\a_{1}=a_{-1}=1.56808208
			\\ \\a_{2}=a_{-2}=-0.17723283
			\\ \\a_{3}=a_{-3}=0.01564992
			\end{split}
		\end{equation}	
	\end{center}
	are the coefficients obtained for the 4\textsuperscript{th} order second derivative stencil.

	\subsection{Dispersion characteristics of finite-difference schemes}
\begin{figure*}[h!]
		\begin{center}
			\includegraphics[angle=0, width=0.75\textwidth, keepaspectratio=true]{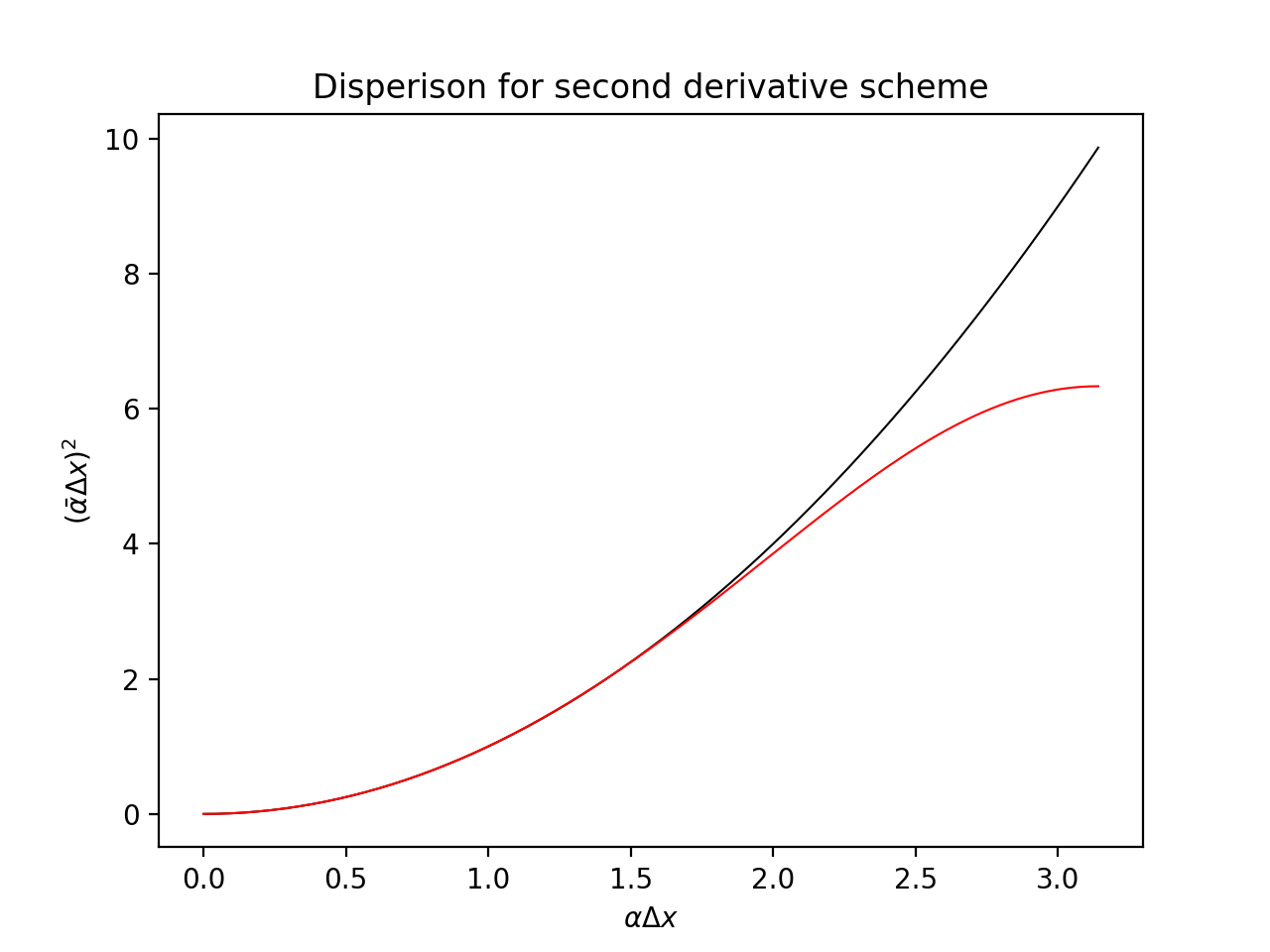} 
			\caption{$(\bar{\alpha}\Delta x)^{2}$ versus $\alpha\Delta x$ for the 4\textsuperscript{th} order optimized central-difference scheme. The solid black line designates the ideal: $(\bar{\alpha}\Delta x)^{2}=(\alpha\Delta x)^{2}$. The curve and ideal are effectively overlaid for $\alpha\Delta x\leq 1.55$, beyond which the two begin to diverge. Therefore the optimized stencil will accurately represent the partial derivative for  wavelengths greater than $4.1\Delta x$. Beyond this point, misfit between the ideal and optimized curves increases with larger values of $\alpha\Delta x$. Wave components with wavelengths shorter than $4.1\Delta x$ will not have their propagation characteristics accurately preserved and values beyond the inflection point at $\alpha\Delta x = \pi$ will produce dispersion characteristics that differ considerably from those of the PDE.}
			\label{2nd_dispersion_plot_single}
		\end{center}
	\end{figure*}

\begin{figure*}[h!]
		\begin{center}
			\includegraphics[angle=0, width=0.75\textwidth, keepaspectratio=true]{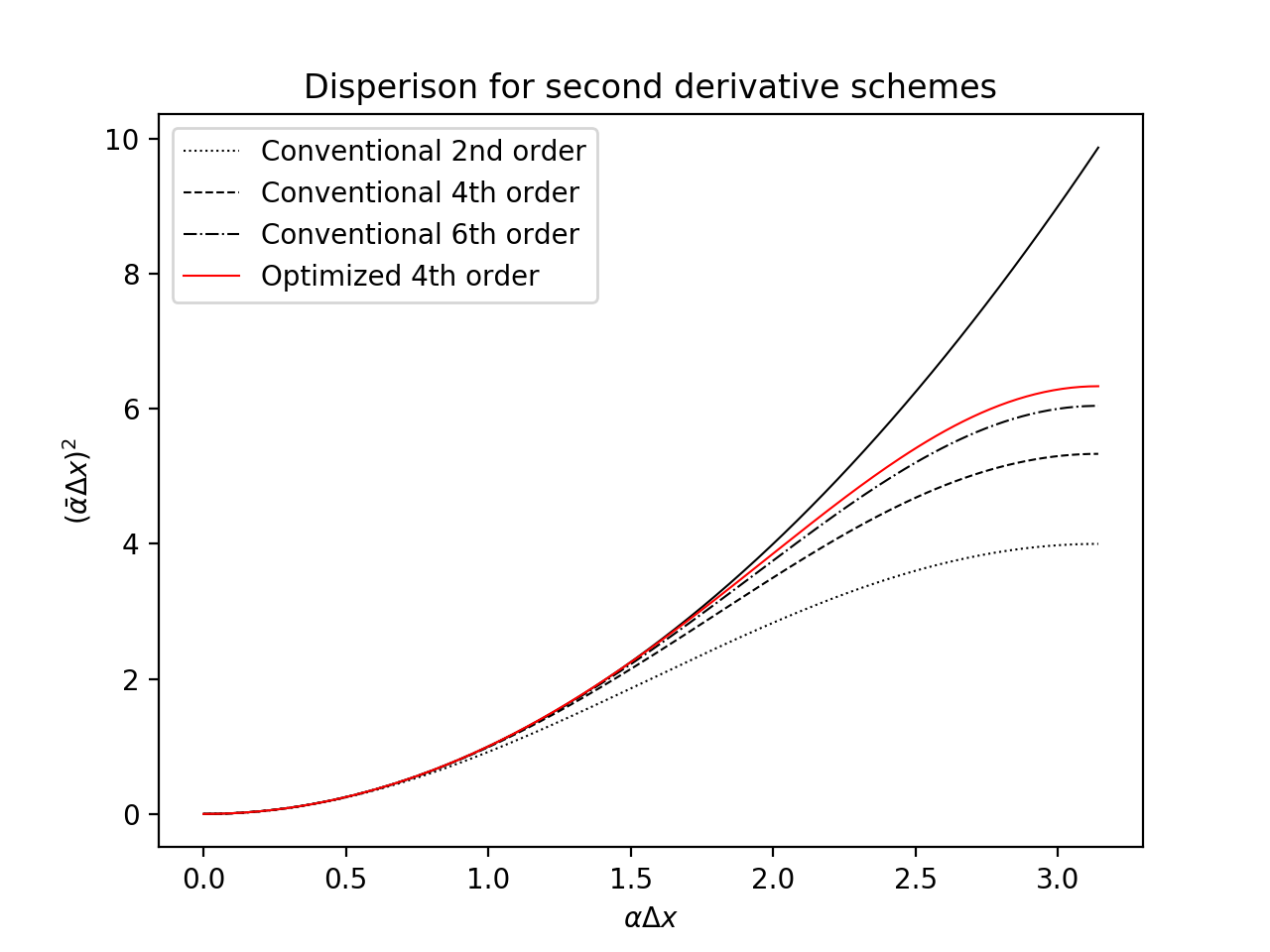} 
			\caption{$(\bar{\alpha}\Delta x)^{2}$ versus $\alpha\Delta x$ for a selection of central-difference schemes. Note the considerable reduction in misfit for the 4\textsuperscript{th} order optimized scheme versus other discretizations. Considering the dispersion curve of the conventional 6\textsuperscript{th} order scheme, it is apparent that deviation from the ideal begins at  $\alpha\Delta x = 1.35$, implying adequate approximation only if the shortest wavelength present is greater than $4.7\Delta x$. Thus it can be concluded that a 4\textsuperscript{th} order optimized scheme has greater resolving power than its conventional 6\textsuperscript{th} order counterpart with equal stencil extent ($N=M=3$). Also shown are dispersion curves for 2\textsuperscript{nd} and 4\textsuperscript{th} order conventional schemes ($N=M=1$ and $N=M=2$ respectively). Both schemes deviate from the ideal at much smaller $\alpha\Delta x$: unable to resolve waves with wavelengths smaller than $9.7\Delta x$ and $6.3\Delta x$ respectively, emphasizing the necessity of high-order and optimized schemes for handling high-frequency components.}
			\label{2nd_dispersion_plot}
		\end{center}
	\end{figure*}

	The relationship between $(\bar{\alpha}\Delta x)^{2}$ and $\alpha\Delta x$ over the first half period for the above scheme is shown in figure \ref{2nd_dispersion_plot_single}, and alongside other FD schemes in figure \ref{2nd_dispersion_plot}.
	\clearpage

 	\subsection{Finite difference implementation}
 	Implementation of FD schemes was carried out using Devito, a domain-specific Python module and compiler for FD applications (\citealp{Luporini2018}). Devito allows for model specification with a handful of high-level symbolic Python objects to build an FD operator, used to generate highly optimized C++ code at runtime via a series of intermediate representations, allowing for complex multi-stage optimizations (\citealp{Luporini2018}). The generated code is automatically executed, with the solution returned.
 	\\
 	\\Specification of FD stencils is reasonably straightforward, using index notation to access nodes based on their position relative to the central stencil point. Coefficients are specified for each stencil point within the optimized stencil, which can then be used to build the FD operator. By defining ``time" and ``x" as their respective dimensions, the 4\textsuperscript{th} order accurate spatially-optimized implementation for the 1D acoustic wave equation can be expressed as follows
 	\begin{verbatim}
 	# Optimized stencil coefficients
 	a_0 = -2.81299833
 	a_1 = 1.56808208
 	a_2 = -0.17723283
 	a_3 = 0.01564992
 	
 	eq_opt =  (a_3*u_opt[time, x - 3]
	             + a_2*u_opt[time, x - 2]
	             + a_1*u_opt[time, x - 1]
	             + a_0*u_opt[time, x]
	             + a_1*u_opt[time, x + 1]
	             + a_2*u_opt[time, x + 2]
	             + a_3*u_opt[time, x + 3])/h_x**2- u_opt.dt2
 	\end{verbatim}
 	Devito simplifies stencil specification by automatically handling attempts to access out of bounds grid points, negating the need for additional stencils at the boundaries of the computational domain, and can handle a suite of common boundary conditions (\citealp{Luporini2018}). Applying the FD operator built with this stencil for a specified number of timesteps allows the updated wavefield to be obtained.
 	\\
 	\\Symbolic computation to implement the test cases expedites workflow compared to model building with low-level languages, enabling rapid prototyping in hours as opposed to weeks or months. Performance optimization during compilation results in generated C++ of comparable quality to hand-optimized implementations (\citealp{Luporini2018}), whilst allowing attention to be focused on the problem at hand. The improved readability of high-level code allows for straightforward transfer from governing PDEs to the computational model. There is presently a dearth of literature using symbolic computation as a tool for implementation DRP schemes, with the method remaining unnoticed as far.
 	\\
 	\\Test cases and example code are located in the supplementary materials.

	\section{Test case: 1D acoustic wave propagation}
	\subsection{Model formulation}
	To evaluate the effectiveness of the 4\textsuperscript{th} order spatially-optimized 2\textsuperscript{nd} derivative, a 1D acoustic wave problem was formulated and solved using optimized and conventional stencils. The 1D acoustic wave equation has the form
	\begin{center}	
		\begin{equation}
		\frac{\partial^{2}u}{\partial t^{2}}=c_{0}^{2}\frac{\partial^{2}u}{\partial x^{2}}
		\label{eq:10}
		\end{equation}
	\end{center}
	where $c_{0}$ is the wavespeed. Solutions consist of superpositioned wavetrains of the form
	\begin{center}	
		\begin{equation}
		u(x, t)=Ae^{i\alpha x - i\omega t}
		\label{eq:11}
		\end{equation}
	\end{center}
	where $A$ represents amplitude and $\omega$ is frequency. Inserting this into (\ref{eq:10}) yields
	\begin{center}	
		\begin{equation}
		\omega^{2}=c_{0}^{2}\alpha^{2}
		\label{eq:12}
		\end{equation}
	\end{center}
	This is the dispersion relation for this equation (\citealp{Whitham1974}). As the wave will propagate with phase velocity $c=\frac{\omega}{\alpha}=\pm c_{0}$, it is apparent that the 1D acoustic wave equation is non-dispersive. Any numerical dispersion will thus be particularly obvious.  
	\\
	\\Straightforward implementation and and a simple analytic solution also contributed to the selection of this initial test case. The 1D domain minimized gridpoint count, reducing memory requirements and operations per timestep. Resultant short runtimes enabled rapid debugging and parameter tuning. Zero Dirichlet boundary conditions were applied, recreating the classic ``wave on a string" setup, thereby ensuring intuitive behavior, and further simplifying implementation and the analytical solution. Initial displacement was defined as a square wave with amplitude 0.1 and wavelength equal to half the domain length, approximated via Fourier sine series. This ensured obvious numerical dispersion would affect high-frequency components. The series was limited to 100 terms to prevent requirement of an excessively large number of gridpoints to discretize high-frequency wave components without ambiguity, limiting runtime. Initial $\frac{\partial u}{\partial t}=0$ for all $x$, again for ease of implementation. A value of $c_{0}=1ms^{-1}$ was selected to simplify the PDE, and domain extent was specified as 10m, producing a standing wave with a period of 5s. A Courant number of 0.2 was empirically determined, achieving acceptable trade-off between accuracy and expediency.

	\subsection{Analytical solution}
	As a benchmark for comparison, the system was solved analytically. Consider equation (\ref{eq:10}): with $c_{0}=1ms^{-1}$, this becomes
	\begin{center}	
		\begin{equation}
		\frac{\partial^{2}u}{\partial t^{2}}=\frac{\partial^{2}u}{\partial x^{2}}
		\label{eq:13}
		\end{equation}
	\end{center}
	Assuming the solution $u(x, t)$ can be expressed as
	\begin{center}	
		\begin{equation}
		u(x, t)=T(t)X(x)
		\label{eq:14}
		\end{equation}
	\end{center}
	then via separation of variables, it is straightforward to find that the general solution takes the form
	\begin{center}	
		\begin{equation}
		u(x, t)= \left( A\sin(ax)+B\cos(ax)\right)\left( C\sin(at)+D\cos(at)\right)
		\label{eq:15}
		\end{equation}
	\end{center}
	where $a, A, B, C, D$ are constants. Applying the boundary conditions
	\begin{center}
		\begin{equation}
		u(0,t) = 0; \; u(L,t) = 0
		\label{eq:16}
		\end{equation}	
	\end{center}
	where $L$ is domain length yields $B=0$ and $a=\frac{n\pi}{L}$ for $n\in\mathbb{Z}^{+}$. Dividing through by $A$:
	\begin{center}
		\begin{equation}
		u(x, t)=\sin\left(\frac{n\pi x}{L}\right)\left(\tilde{C}\sin\left(\frac{n\pi t}{L}\right)+\tilde{D}\cos\left(\frac{n\pi t}{L}\right)\right)
		\label{eq:17}
		\end{equation}	
	\end{center}
	Applying the first initial condition
	\begin{center}
		\begin{equation}
		\begin{split}
		u(x,0) = \sum_{n=1}^{100}b_{n}\sin\left(\frac{2n\pi x}{L}\right) \\ \\
		b_{n} = \frac{4}{L}\int_{0}^{\frac{L}{2}}f(x)\sin\left(\frac{2n\pi x}{L}\right)\mathrm{d}x \\ \\
		f(x)=
		\begin{cases}
		1,& \text{if } 0\leq x\leq\frac{L}{4} \text{ or }  \frac{L}{2}\leq x\leq\frac{3L}{4}\\ \\
		-1,& \text{otherwise}
		\end{cases}
		\label{eq:18}
		\end{split}
		\end{equation}	
	\end{center}
	$\tilde{D}$ is found to take the form
	\begin{center}
		\begin{equation}
		\tilde{D}_{n} = 
		\begin{cases}
		\frac{4}{L}\int_{0}^{\frac{L}{2}}f(x)\sin\left(\frac{n\pi x}{L}\right)\mathrm{d}x,& \text{for all even $n$}\\ \\
		0,& \text{otherwise}
		\end{cases}
		\label{eq:19}
		\end{equation}	
	\end{center}
	In applying the second initial condition
	\begin{center}
		\begin{equation}
		\frac{\partial u}{\partial t}(x, 0) = 0
		\label{eq:20}
		\end{equation}	
	\end{center}
	the particular solution for the PDE is found to be
	\begin{center}
		\begin{equation}
		u(x,t)=\sum_{n=1}^{200}\tilde{D}_{n}\sin\left(\frac{n\pi x}{L}\right)\cos\left(\frac{n\pi t}{L}\right)
		\label{eq:21}
		\end{equation}	
	\end{center}

	\subsection{Performance at near-maximum grid spacings}
\begin{figure*}[h!]
		\begin{center}
			\includegraphics[angle=0, width=0.75\textwidth, keepaspectratio=true]{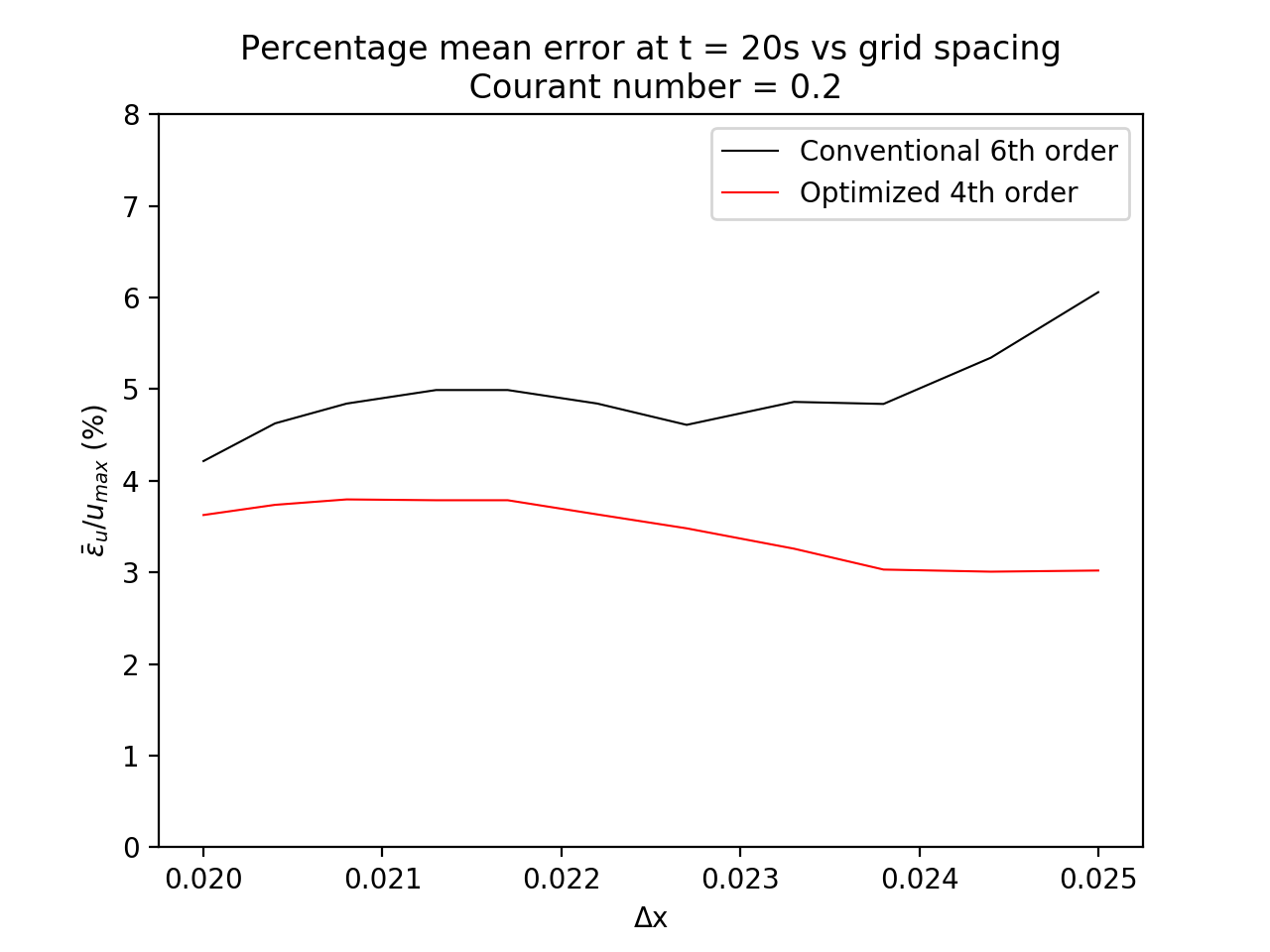} 
			\caption{Mean absolute error ($\bar{\epsilon_{u}}$) evaluated at $t=20s$, normalized against maximum amplitude versus $\Delta x$. Performance for the optimized stencil is considerably better than for the conventional at $\Delta x=0.0250\mathrm{m}$, with $\frac{\bar{\epsilon_{u}}}{u_{max}}=3.0\%$ and $\frac{\bar{\epsilon_{u}}}{u_{max}}=6.1\%$ respectively. This broadly agrees with the relationship in figure \ref{2nd_dispersion_plot}: the optimized scheme is expected to offer superior performance in this region. Error for the two schemes converges as expected towards $\Delta x= 0.0213\mathrm{m}$. However, the somewhat rough relationship between $\frac{\bar{\epsilon_{u}}}{u_{max}}$ and $\Delta x$ throughout the window studied prevents meaningful conclusions based solely upon this data. For adequately small $\Delta x$, the conventional scheme should offer equivalent or better performance. Surprisingly, for all $\Delta x$ within the region of interest, the optimized scheme generates smaller $\bar{\epsilon_{u}}$: unexpected given the conventional scheme should accurately propagate at the smallest $\Delta x$ tested. Furthermore, the smallest error for the conventional scheme is greater than the largest error for the optimized within this window. Contrary to expectation, error in the optimized scheme decreases between $\Delta x=0.0200\mathrm{m}$ and $\Delta x=0.0250\mathrm{m}$. Given the $\Delta x$ interval used, it is unlikely that this reflects small-scale roughness in the relationship between $\Delta x$ and $\bar{\epsilon_{u}}$ (although this relationship is seen to be somewhat rough for the conventional scheme). However, it is inconsistent with expected behavior of both truncation and numerical dispersion errors, inferring the presence of some additional factor, possibly concerning the formulation of the test case.}
			\label{percent_mean_err_dx_zoomed}
		\end{center}
	\end{figure*}

	From figure \ref{2nd_dispersion_plot}, the optimized scheme should accurately approximate propagation characteristics of the standing wave for $\Delta x\leq0.0244\mathrm{m}$, with $\Delta x\leq0.0213\mathrm{m}$ required for the conventional scheme. Within the region $0.0213\mathrm{m}\leq \Delta x\leq 0.0244\mathrm{m}$ the optimized stencil should yield reduced error; the error of the two schemes converging towards $\Delta x=0.0213\mathrm{m}$ with similar or superior performance for the conventional scheme thereafter. To verify, both schemes were applied with grid spacings between $0.0200\mathrm{m}$ and $0.0250\mathrm{m}$, incremented corresponding to addition of 10 gridpoints. The standing wave was propagated over 4 periods (total length 20s) to ensure numerical dispersion error would be apparent at simulation end. Error accumulated by longer simulations resulted in excessively noisy data, limiting propagation time. Absolute misfit between the numerically propagated wavefields and analytic solution ($\epsilon_{u}$) was evaluated at each grid point, with the mean calculated for each scheme (see figure \ref{percent_mean_err_dx_zoomed}). Mean error was selected over maximum, as small phase differences were capable of producing large maximum errors, failing to reflect overall simulation quality.
	\begin{figure*}[h!]
		\begin{center}
			\begin{subfigure}{0.49\textwidth}
				\centering\includegraphics[width=\textwidth]{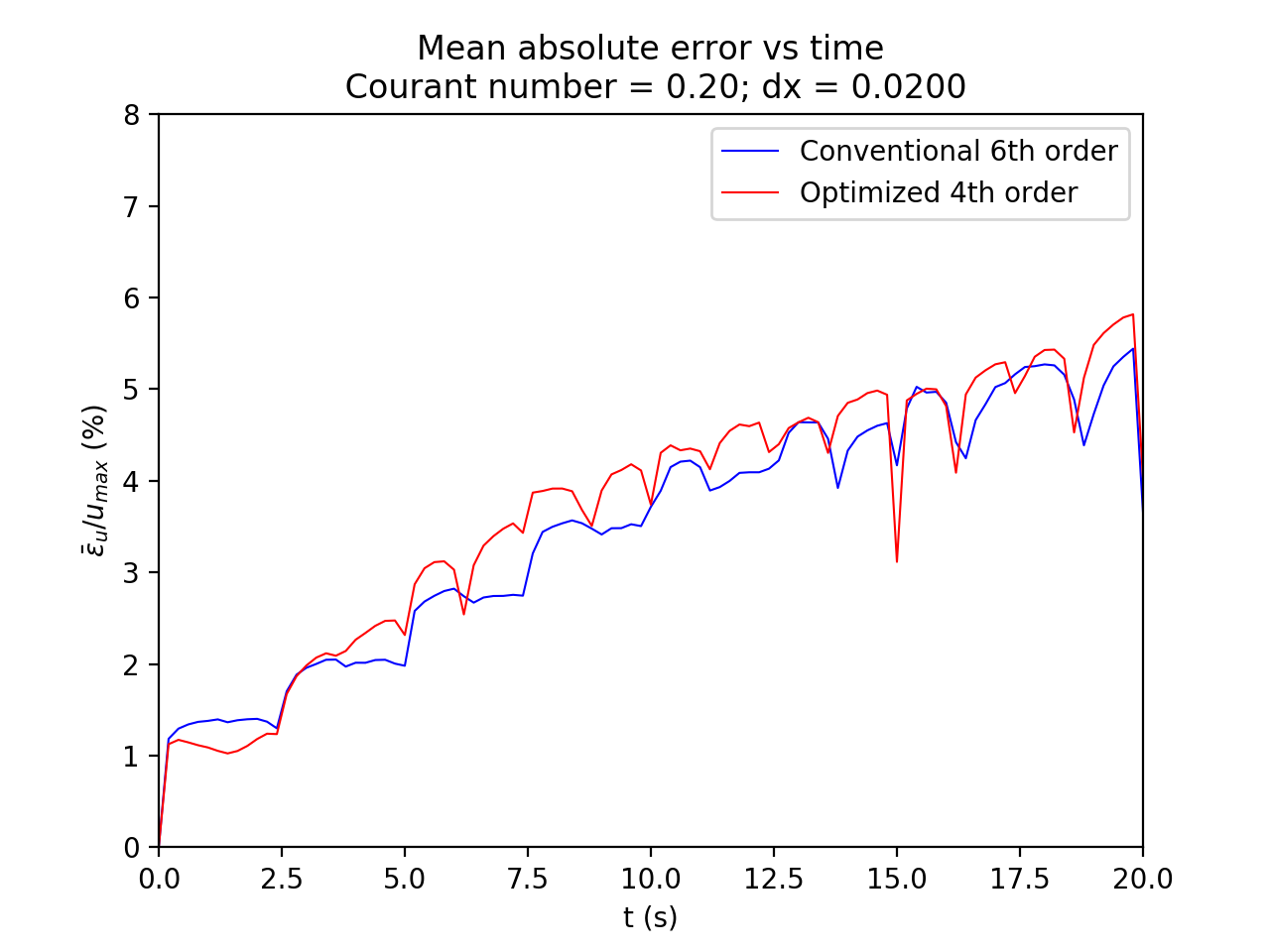}
				\caption{}
			\end{subfigure}
			\begin{subfigure}{0.49\textwidth}
				\centering\includegraphics[width=\textwidth]{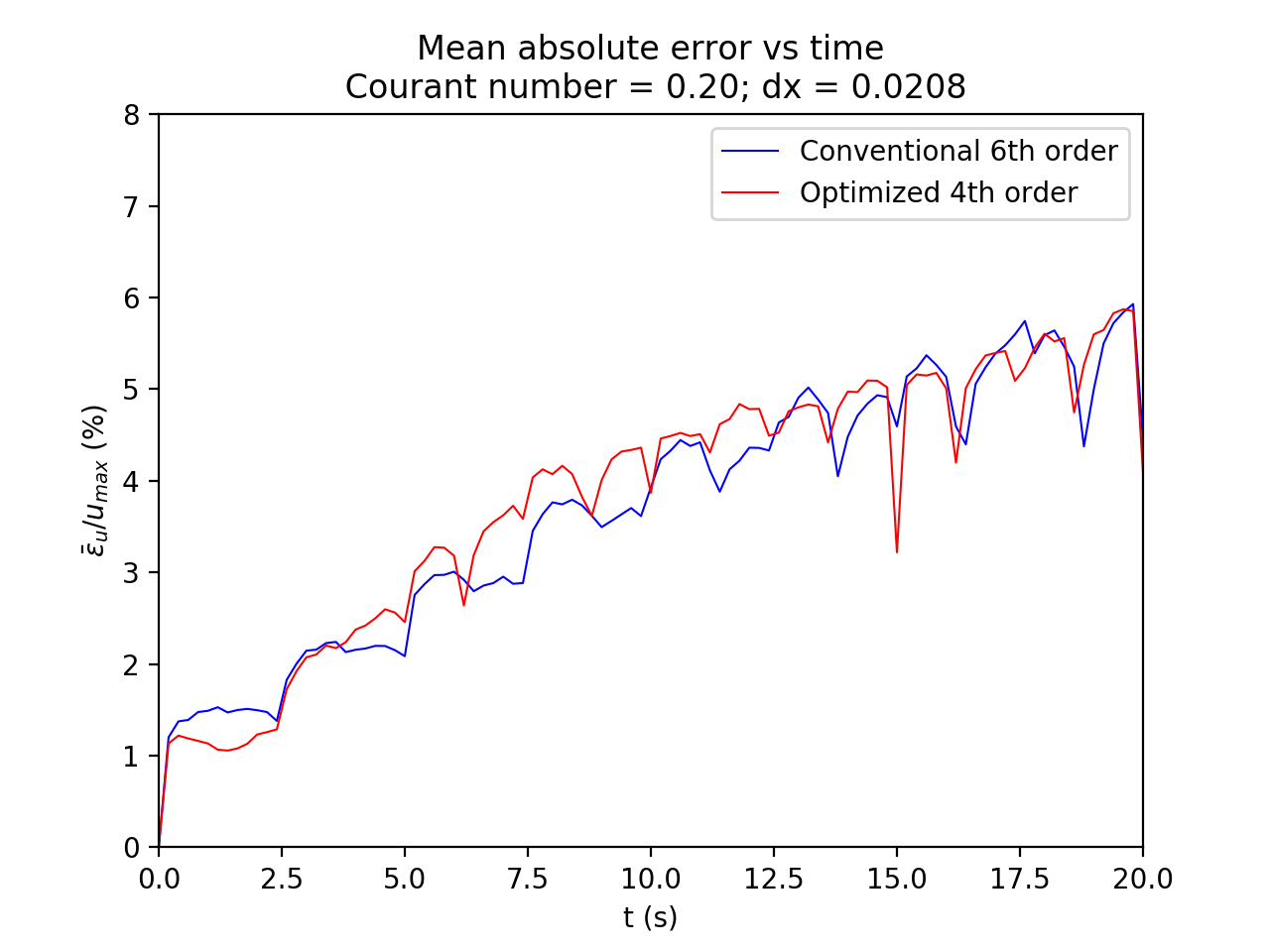}
				\caption{}
			\end{subfigure}
		
			\begin{subfigure}{0.49\textwidth}
				\centering\includegraphics[width=\textwidth]{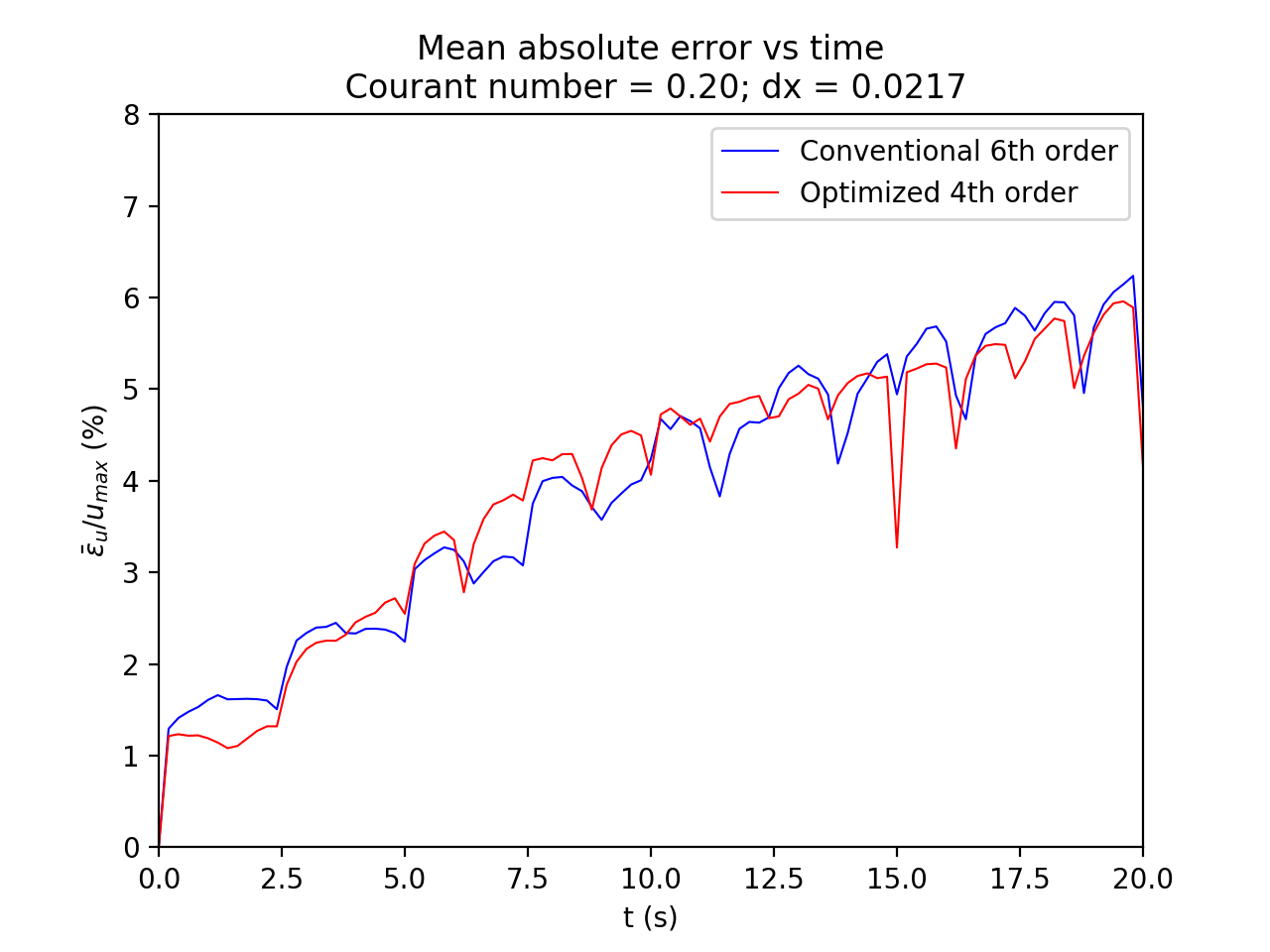}
				\caption{}
			\end{subfigure}
			\begin{subfigure}{0.49\textwidth}
				\centering\includegraphics[width=\textwidth]{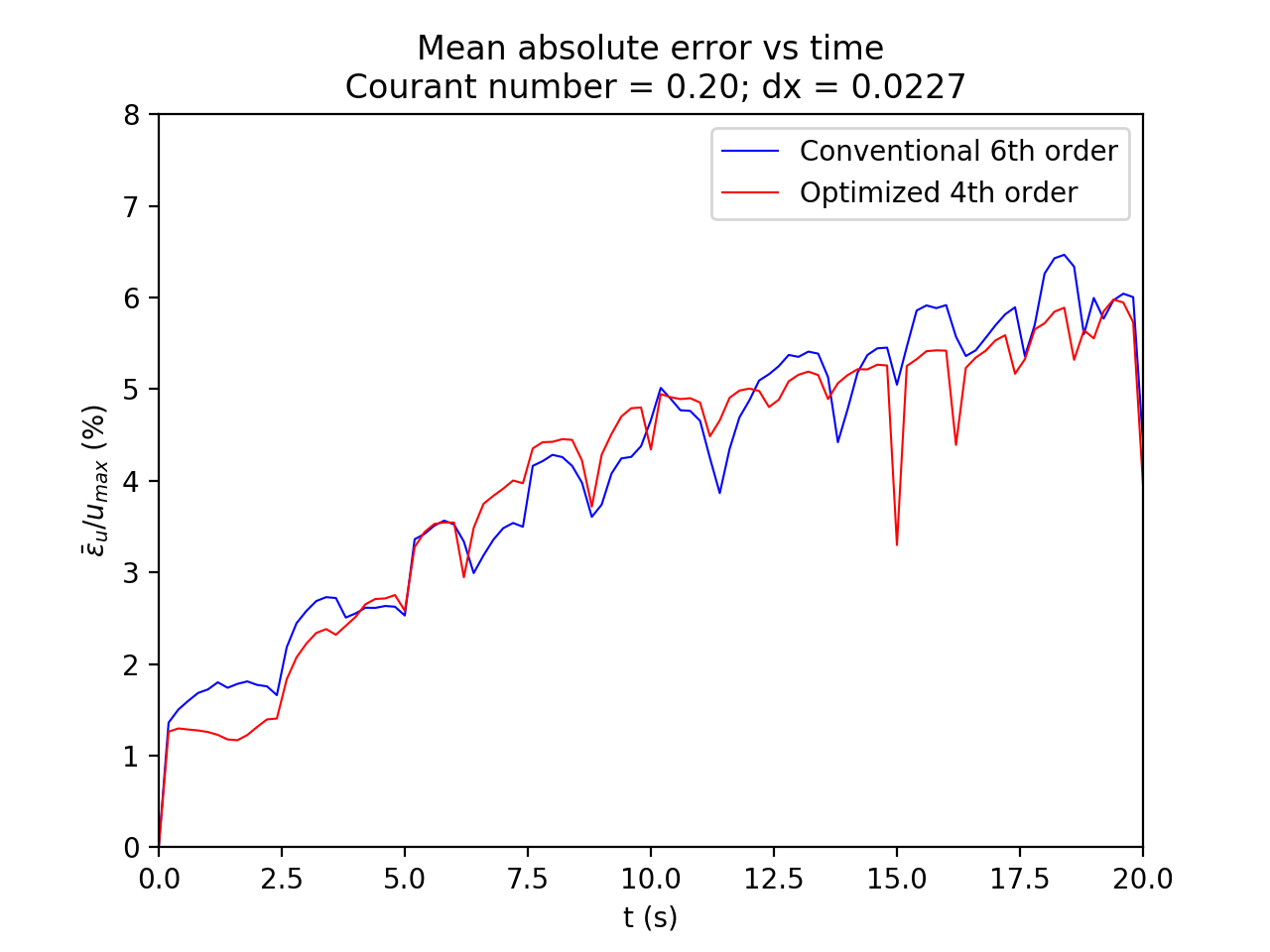}
				\caption{}
			\end{subfigure}
		\caption{$\frac{\bar{\epsilon_{u}}}{u_{max}}$ versus time for various $\Delta x$. (Figure continued overleaf.)}
		\end{center}
	\end{figure*}
	\begin{figure*}[h!]
		\begin{center}
		\ContinuedFloat
			\begin{subfigure}{0.49\textwidth}
				\centering\includegraphics[width=\textwidth]{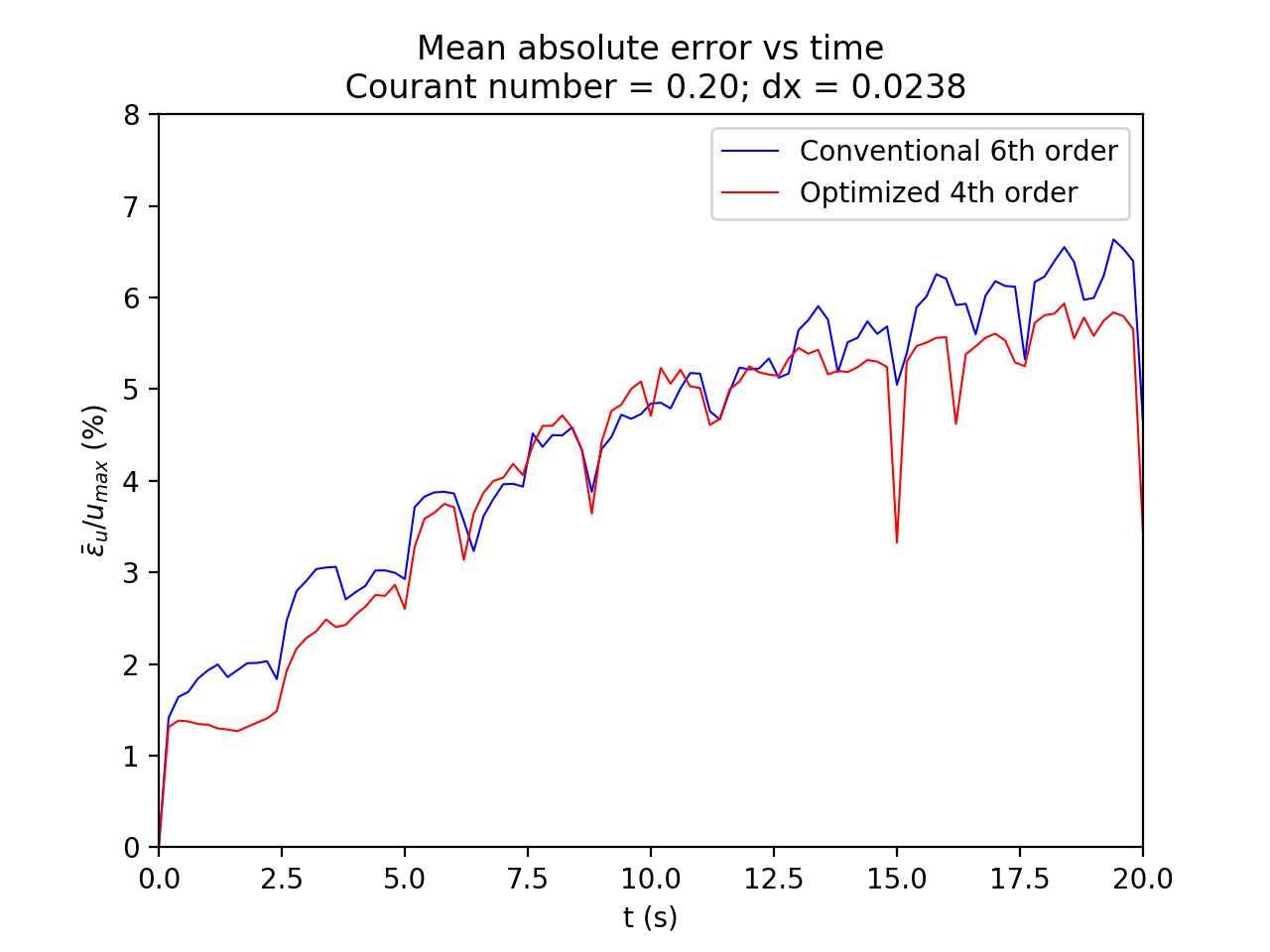}
				\caption{}
			\end{subfigure}
			\begin{subfigure}{0.49\textwidth}
				\centering\includegraphics[width=\textwidth]{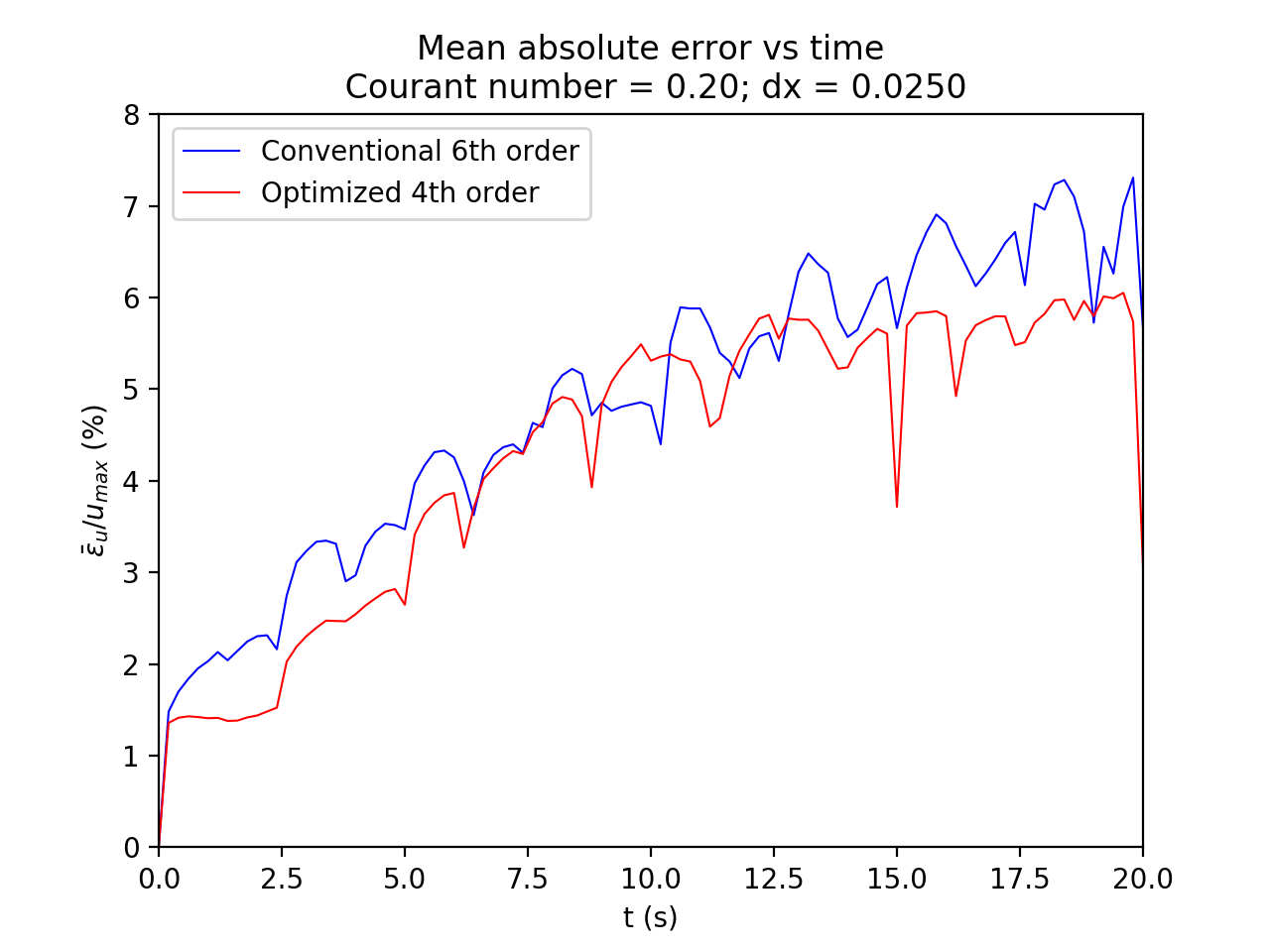}
				\caption{}
			\end{subfigure}
			\caption{$\frac{\bar{\epsilon_{u}}}{u_{max}}$ versus time for various $\Delta x$. Note the periodic dips in error coinciding with half the period of the standing wave, explaining the increased upper values of error when compared to figure \ref{percent_mean_err_dx_zoomed}.  The optimized scheme is clearly superior when sampling is coarse. Over the course of a simulation using $\Delta x=0.0250\mathrm{m}$, the optimized scheme better fits the analytic solution at 64.0\% of sampled timesteps, and 51.8\% for $\Delta x=0.0238\mathrm{m}$, in agreement with figure \ref{percent_mean_err_dx_zoomed}. Convergence of the schemes as $\Delta x$ is reduced is strongly apparent: both error-time curves trace similar paths in figures \ref{percent_mean_err_time}c and \ref{percent_mean_err_time}d. For $\Delta x=0.0217\mathrm{m}$ and $\Delta x=0.0227\mathrm{m}$, the optimized scheme is more accurate at 23.5\% and 41.9\% of sampled timesteps respectively, implying the performance crossover of the optimized 4\textsuperscript{th} order and conventional 6\textsuperscript{th} order schemes occurs at larger $\Delta x$ than implied by figure \ref{2nd_dispersion_plot}. For $\Delta x=0.0200\mathrm{m}$ and $\Delta x=0.0208\mathrm{m}$ (figures \ref{percent_mean_err_time}a and \ref{percent_mean_err_time}b respectively), the optimized scheme produces the larger error, in agreement with predictions of behavior in this region, although implying that errors at 20s (see figure \ref{percent_mean_err_dx_zoomed}) are misleading and do not accurately represent overall propagation quality. This is affirmed by the optimized scheme having the smaller error at 8.6\% and 11.4\% of sampled timesteps for $\Delta x=0.0200\mathrm{m}$ and $\Delta x=0.0208\mathrm{m}$  respectively.}
			\label{percent_mean_err_time}
		\end{center}
	\end{figure*}
	\\
	\\To further investigate, $\frac{\bar{\epsilon_{u}}}{u_{max}}$ over time was plotted for each $\Delta x$. To avoid impractical runtimes, error was evaluated not at timesteps, but 0.2s intervals. Error plots over the region of interest are illustrated in figure \ref{percent_mean_err_time}.
	\\
	\\Examining wavefield evolution at $\Delta x=0.0227\mathrm{m}$ (the apparent crossover point in figure \ref{percent_mean_err_time}), nature of error in the two solutions appears to differ, despite being of comparable magnitude. For the conventional solution, high-frequency components exhibit phase shifts and distortion, apparent in figure \ref{wavefield_crossover_1}b, particularly on ``flat" sections of the wave. In contrast, error in the spatially-optimized solution appears to originate at steep gradients, generating slight misfit and spurious high frequencies which propagate outwards. This implies differing main error sources at this $\Delta x$, potentially explaining the discrepancy between predicted convergence at $\Delta x\leq0.0213\mathrm{m}$ and actual at approximately $\Delta x\leq0.0227\mathrm{m}$. 
	\begin{figure*}[h!]
		\begin{center}
			\begin{subfigure}{0.75\textwidth}
				\centering\includegraphics[width=\textwidth]{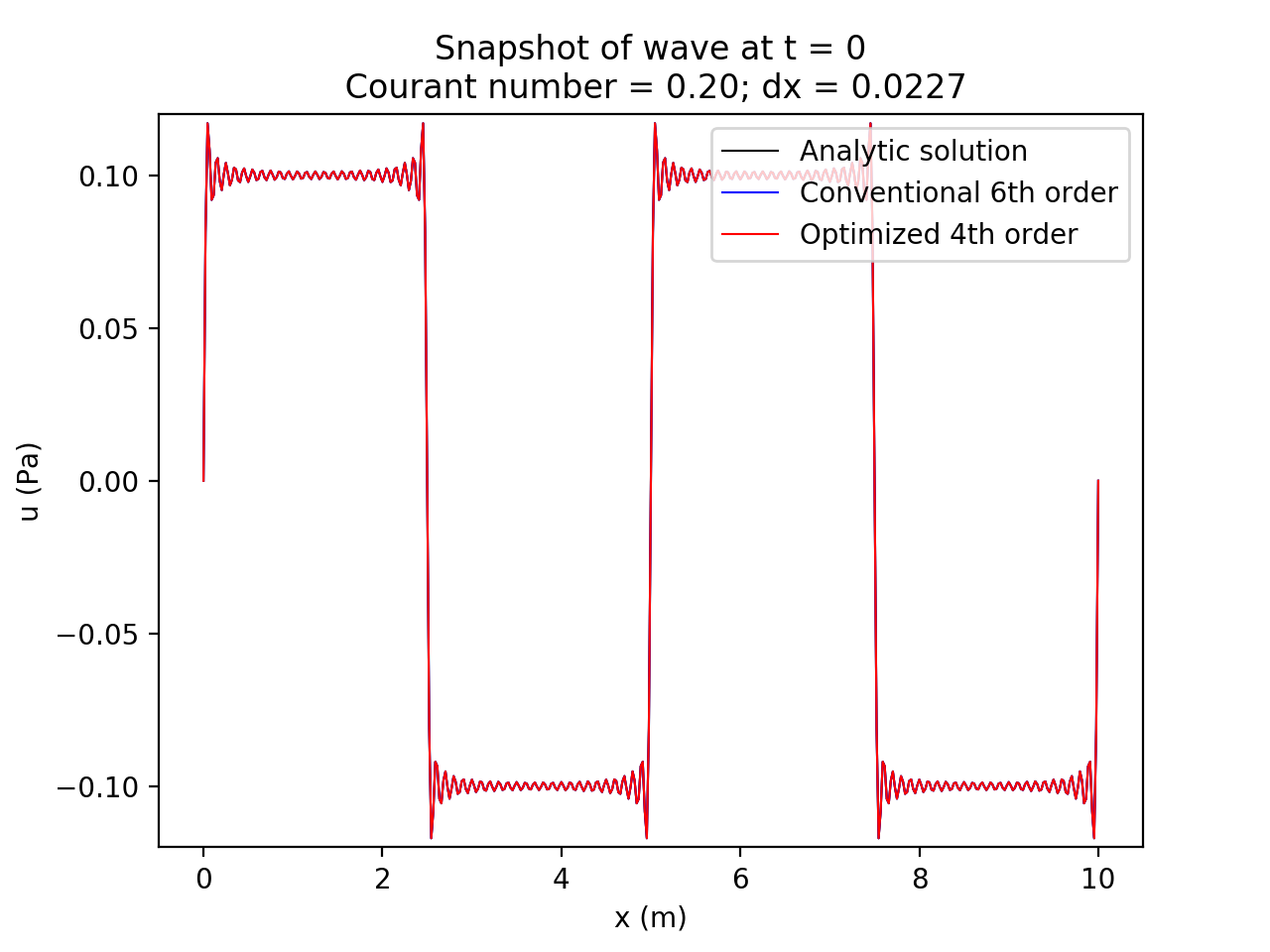}
				\caption{}
			\end{subfigure}
			\begin{subfigure}{0.75\textwidth}
				\centering\includegraphics[width=\textwidth]{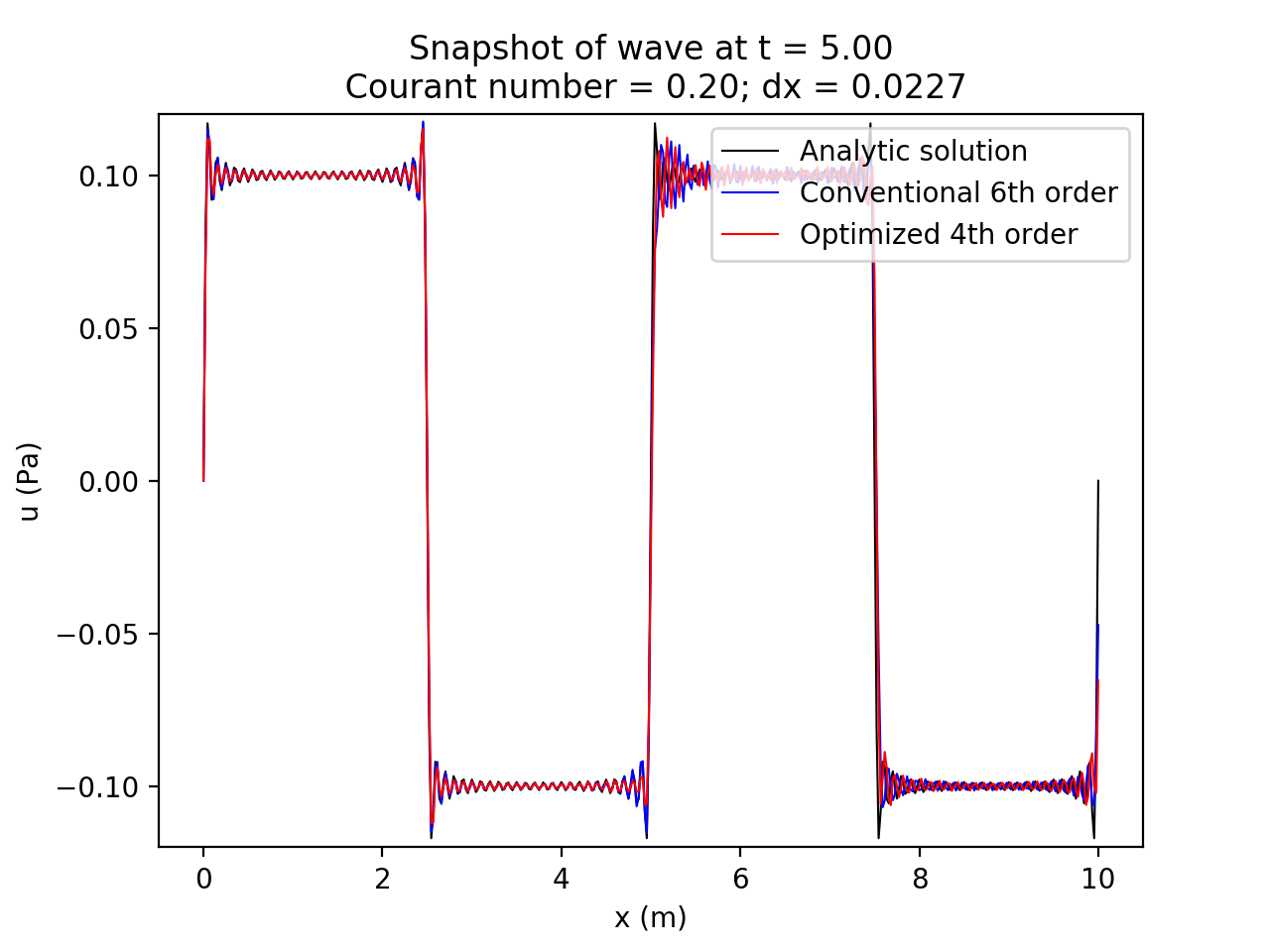}
				\caption{}
			\end{subfigure}
			\caption{(Figure continued overleaf.)}
		\end{center}
	\end{figure*}
	\begin{figure*}[h!]
		\begin{center}
			\ContinuedFloat
			\begin{subfigure}{0.75\textwidth}
				\centering\includegraphics[width=\textwidth]{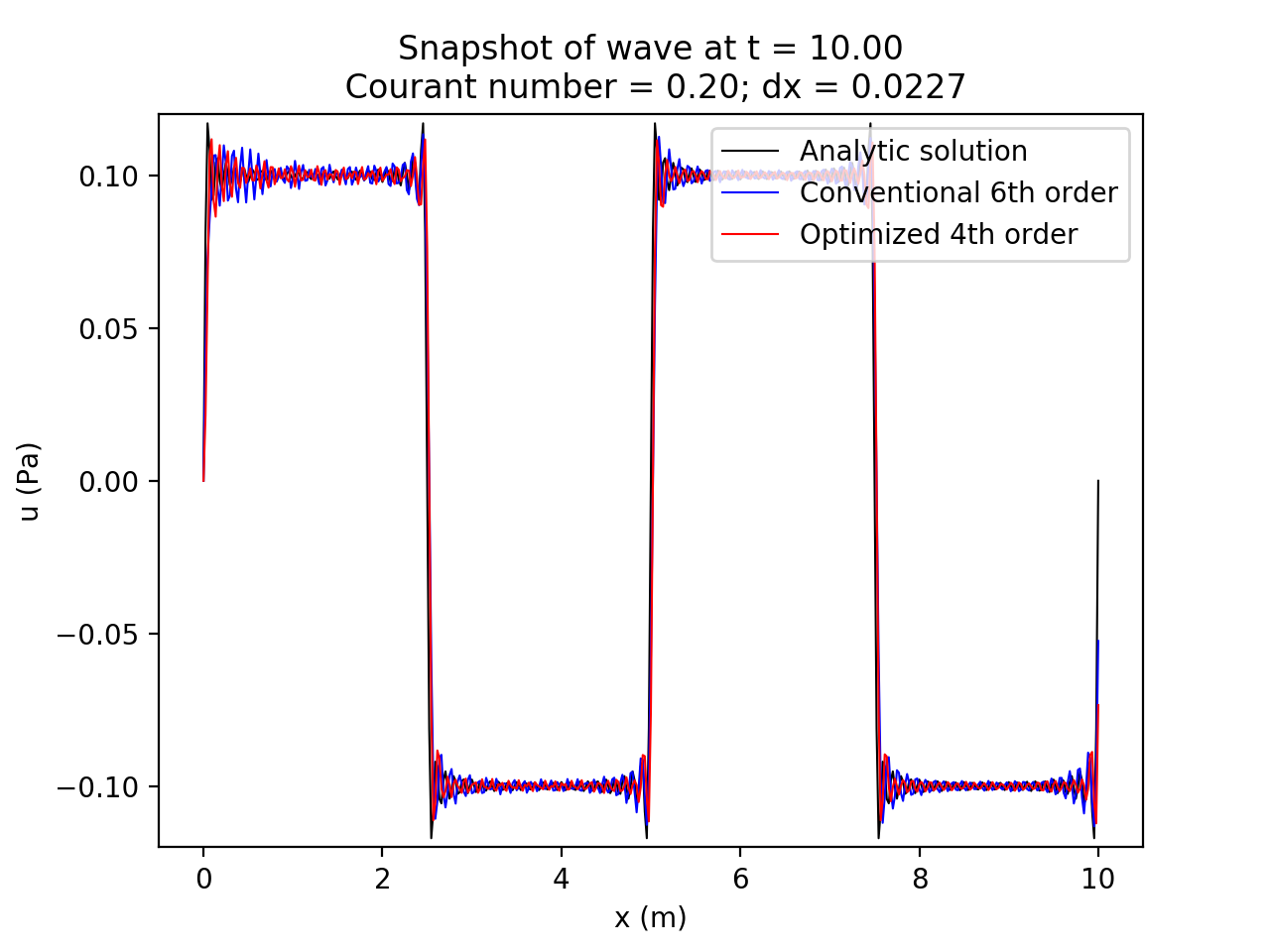}
				\caption{}
			\end{subfigure}
			\begin{subfigure}{0.75\textwidth}
				\centering\includegraphics[width=\textwidth]{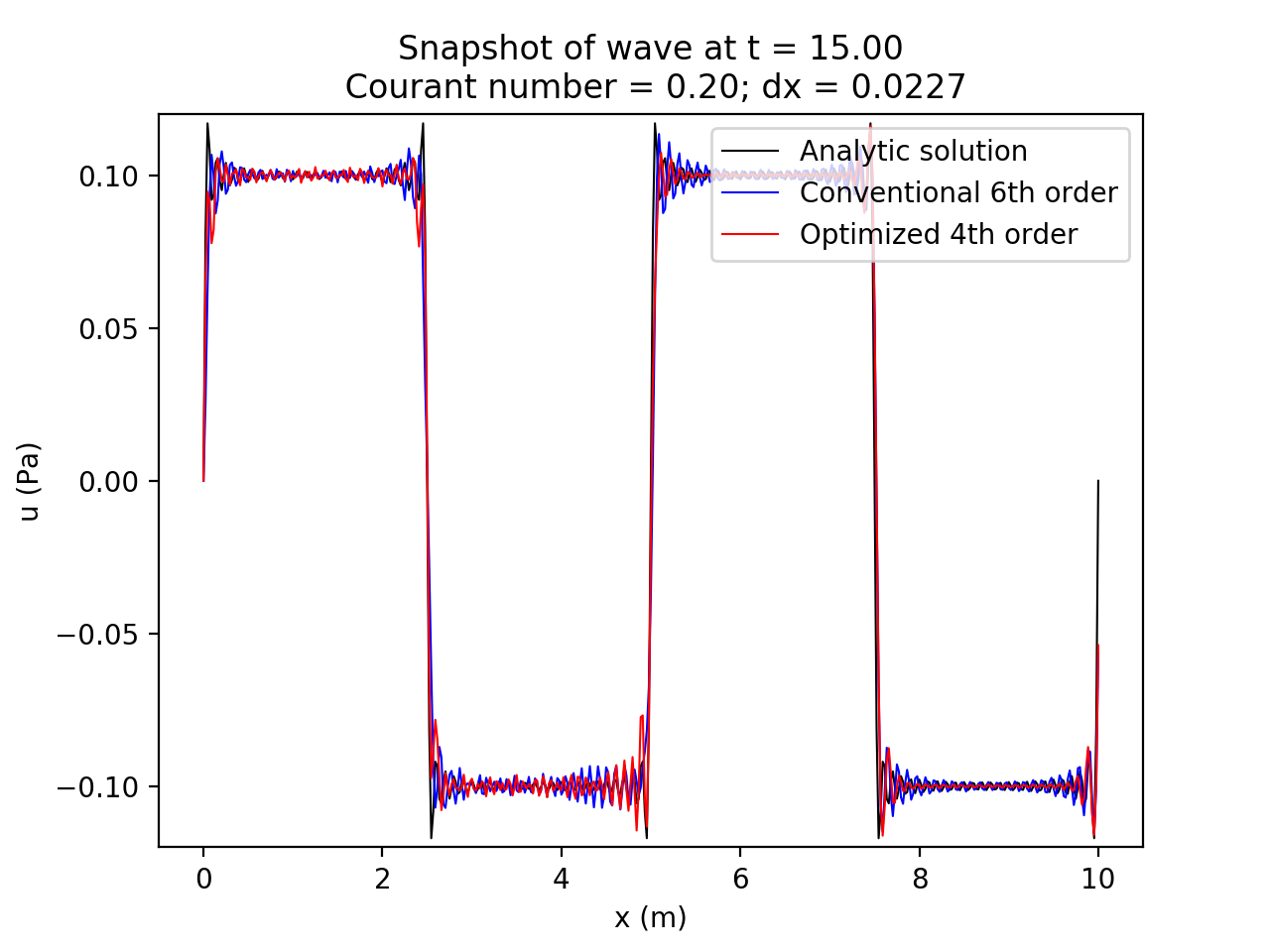}
				\caption{}
			\end{subfigure}
			\caption{(Figure continued overleaf.)}
		\end{center}
	\end{figure*}
	\begin{figure*}[h!]
		\begin{center}
			\ContinuedFloat
			\begin{subfigure}{0.75\textwidth}
				\centering\includegraphics[width=\textwidth]{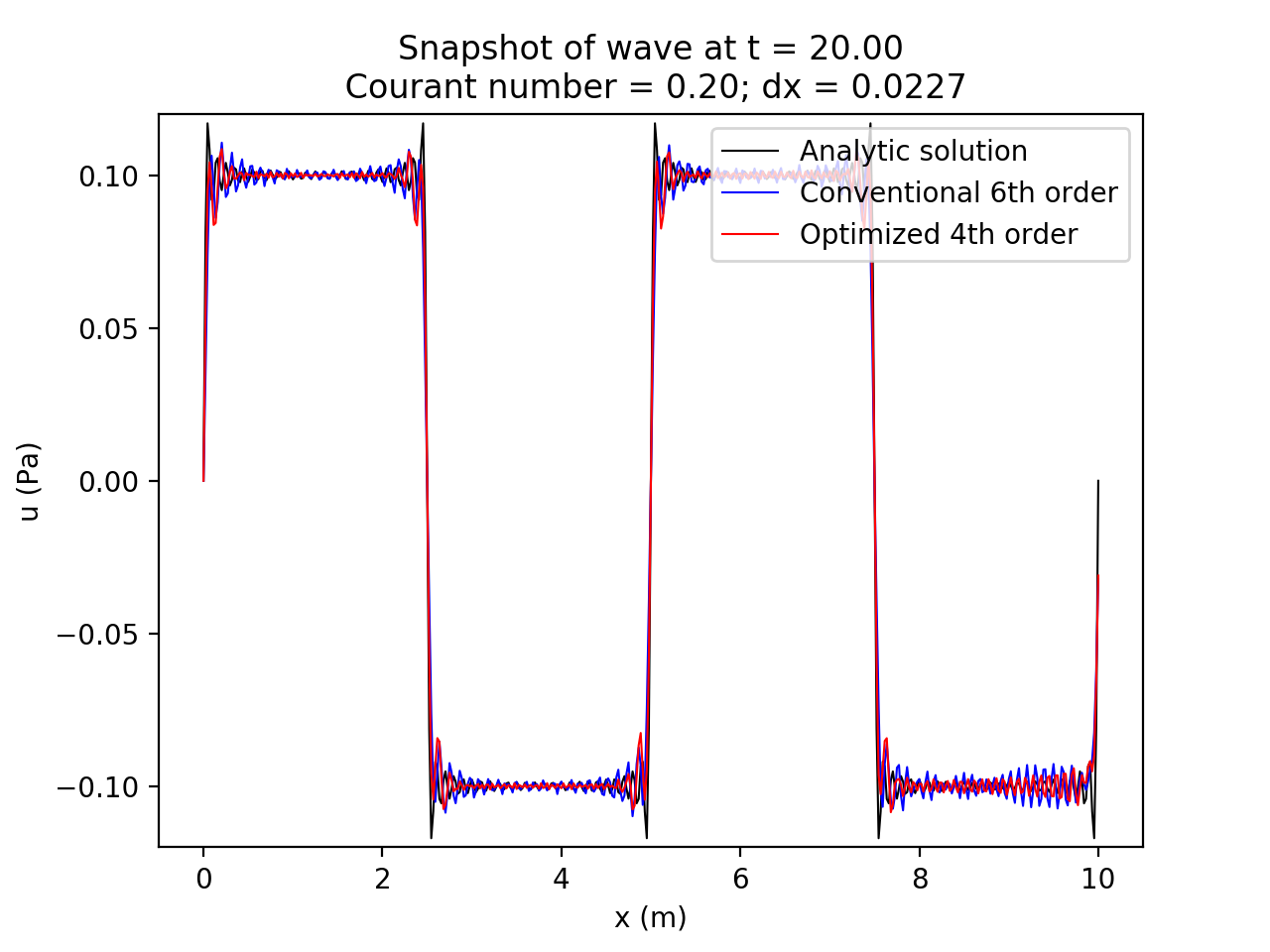}
				\caption{}
			\end{subfigure}
			\caption{Evolution of the analytical, spatially-optimized FD, and conventional FD solutions at 5s intervals for $\Delta x=0.0227\mathrm{m}$.}
			\label{wavefield_crossover_1}
		\end{center}
	\end{figure*}
	\begin{figure*}[h!]
		\begin{center}
			\centering\includegraphics[width=0.75\textwidth]{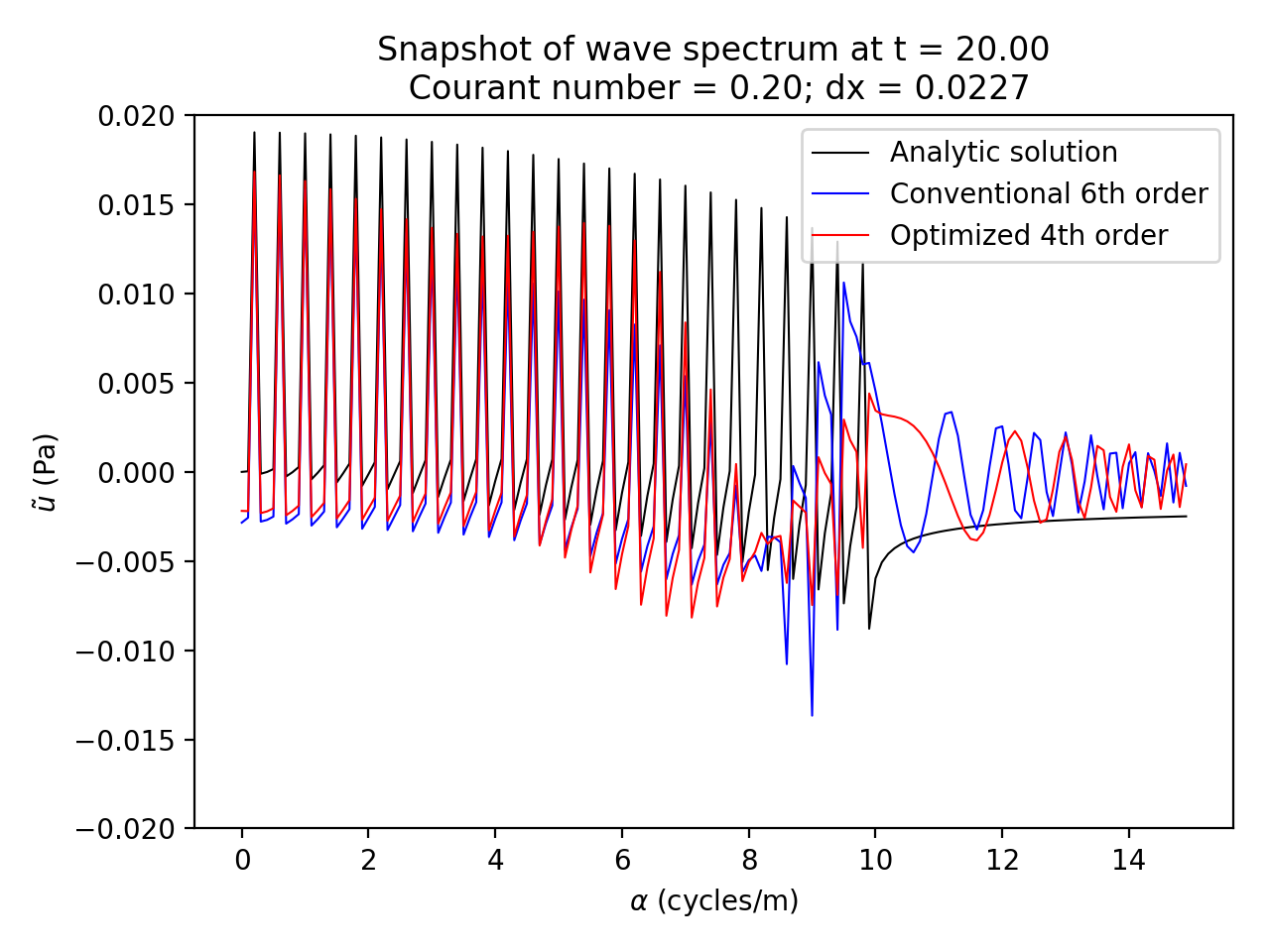}
			\caption{Spatial frequency amplitude spectrum of the analytical, spatially-optimized FD, and conventional FD solutions at 20s for $\Delta x=0.0227\mathrm{m}$. The normalized fast Fourier transform of the wavefield reveals significant distortion and phase shifts in the conventional FD solution, with considerable offset and some high-frequency components near polarity reversed. Whilst the optimized solution deviates from the exact spectrum in this region, the magnitude of misfit is reduced. Induced high frequencies are of approximately equivalent amplitude for both schemes, although they are present in different bands, becoming closer spaced at higher frequencies for both schemes. Overall spectral misfit is slightly greater for the conventional scheme, despite similar misfit in space.}
			\label{wavefield_crossover_spectrum}
		\end{center}
	\end{figure*}
	\clearpage

	\subsection{General performance trends}
\begin{figure*}[h!]
		\begin{center}
			\includegraphics[angle=0, width=0.75\textwidth, keepaspectratio=true]{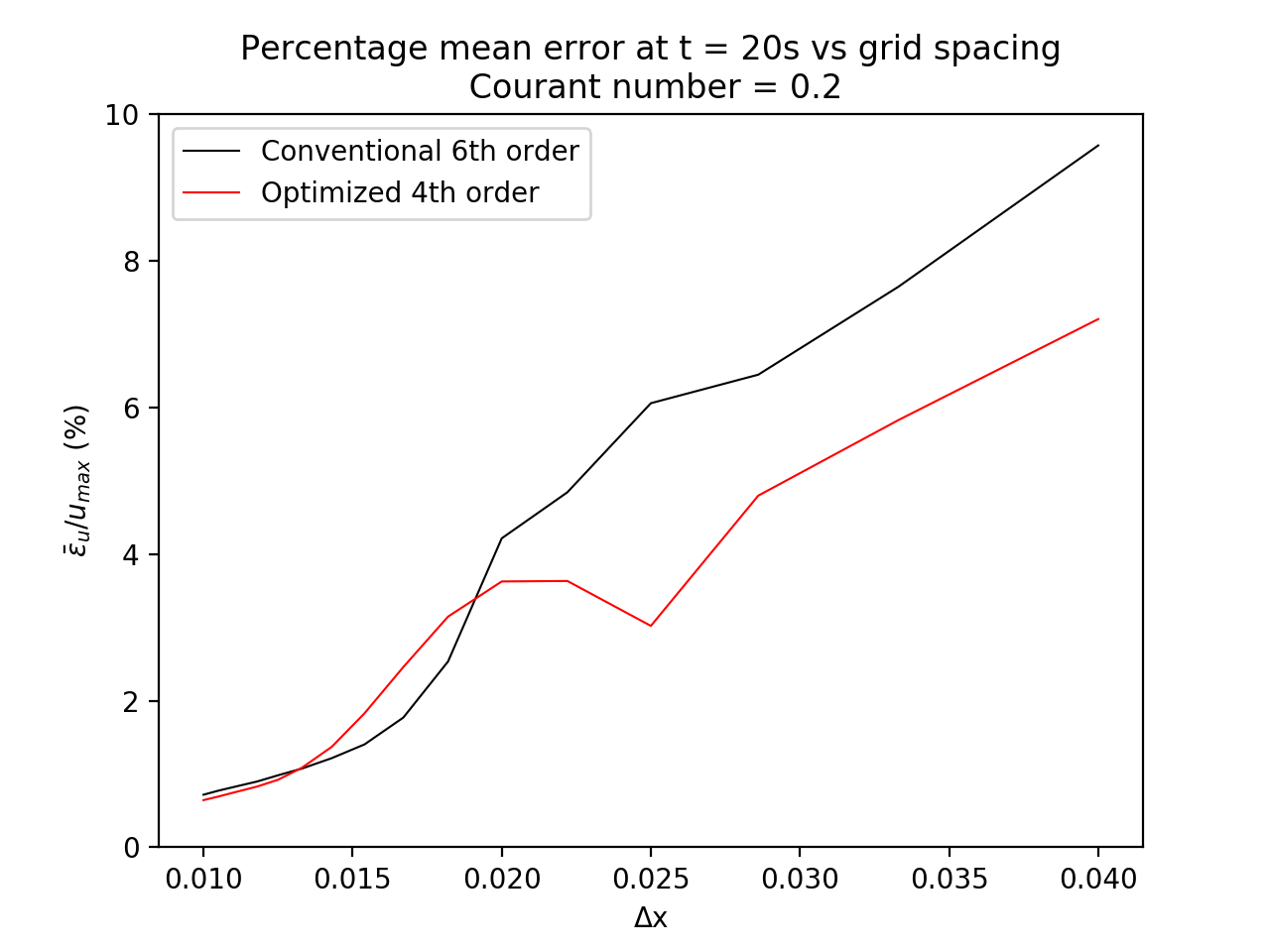} 
			\caption{Mean absolute error ($\bar{\epsilon_{u}}$) evaluated at $t=20s$ normalized against maximum amplitude over a wide range of $\Delta x$. This confirms that the optimized scheme offers considerable simulation quality improvements for undersampled wavefields, with $\frac{\bar{\epsilon_{u}}}{u_{max}}=7.2\%$ versus $\frac{\bar{\epsilon_{u}}}{u_{max}}=9.6\%$ for the conventional scheme at $\Delta x=0.0400\mathrm{m}$. However, overall performance difference between the two schemes is fairly rough when determined by this metric. Lower error in the optimized scheme for coarse grids corresponds with its reduced deviation from ideal dispersion characteristics for short wavelengths shown in figure \ref{2nd_dispersion_plot}. Error magnitude is considerably greater for both schemes when the wavefield is spatially undersampled due to dependence of truncation error on grid spacing. For moderately oversampled wavefields, the conventional FD solution contains reduced error at simulation end, matching expected behavior as minimal numerical dispersion is eclipsed by truncation error. However, for very small $\Delta x$, the two schemes exhibit similar error with marginally less error in the optimized FD solution, implying that either the optimized scheme performs better than expected, or the conventional scheme experiences some undetermined drawback when both are applied to fine grids. The reduction of error in the optimized scheme between $\Delta x=0.0200\mathrm{m}$ and $\Delta x=0.0250\mathrm{m}$ apparent in figure \ref{percent_mean_err_dx_zoomed} is present: the only decrease in error with increased $\Delta x$ within the window studied. The mechanism behind this anomalous behavior is illuminated no further by the wider region of study, although it is of interest that it appears to cease at the maximum $\Delta x$ required to avoid ambiguity in the wavefield.}
			\label{percent_mean_err_dx}
		\end{center}
	\end{figure*}

	For grid spacings where both schemes are subject to numerical dispersion, the 4\textsuperscript{th} order optimized scheme should achieve greater accuracy than the 6\textsuperscript{th} order conventional scheme; deviation from the ideal at the inflection point is smaller (see figure \ref{2nd_dispersion_plot}). For small $\Delta x$, wavefield sampling becomes saturated minimizing or eliminating numerical dispersion. In this circumstance, the 6\textsuperscript{th} order conventional scheme should yield more accurate results as net error for both conventional and optimized schemes will be dominated by truncation error.
	\\
	\\To verify, the test case was repeated over a wider range of grid spacings between $\Delta x=0.0100\mathrm{m}$ and $\Delta x=0.0400\mathrm{m}$, chosen to over and undersample the wavefield at respective extremes. The $\Delta x$ interval was widened to coincide with the addition of 50 grid points to the computational domain, as smaller intervals were deemed prohibitively time consuming.
	\\
	\\Error at simulation end versus $\Delta x$ is shown in figure \ref{percent_mean_err_dx}.
	\begin{figure*}[h!]
		\begin{center}
			\begin{subfigure}{0.75\textwidth}
				\centering\includegraphics[width=\textwidth]{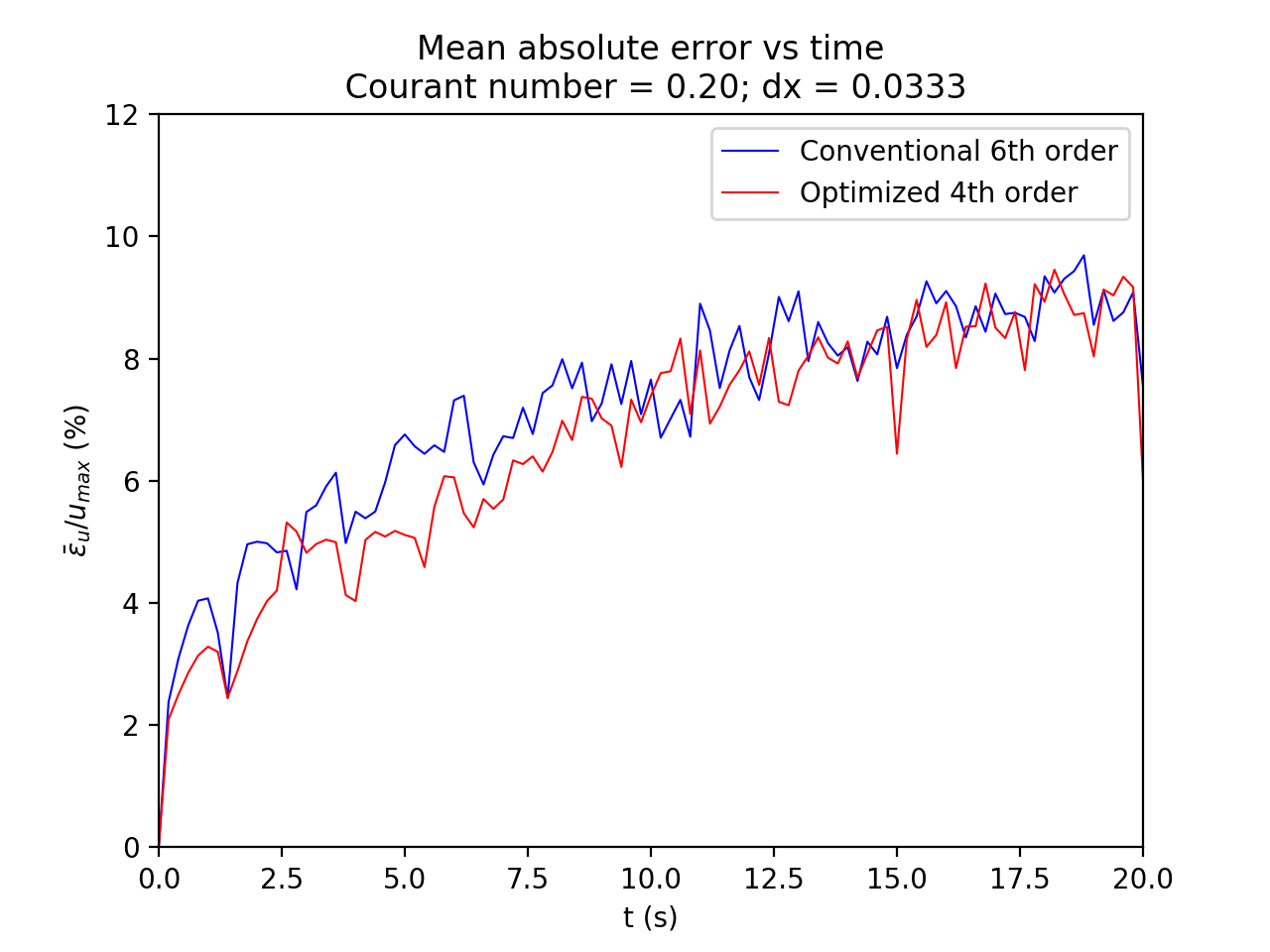}
				\caption{}
			\end{subfigure}
			\begin{subfigure}{0.75\textwidth}
				\centering\includegraphics[width=\textwidth]{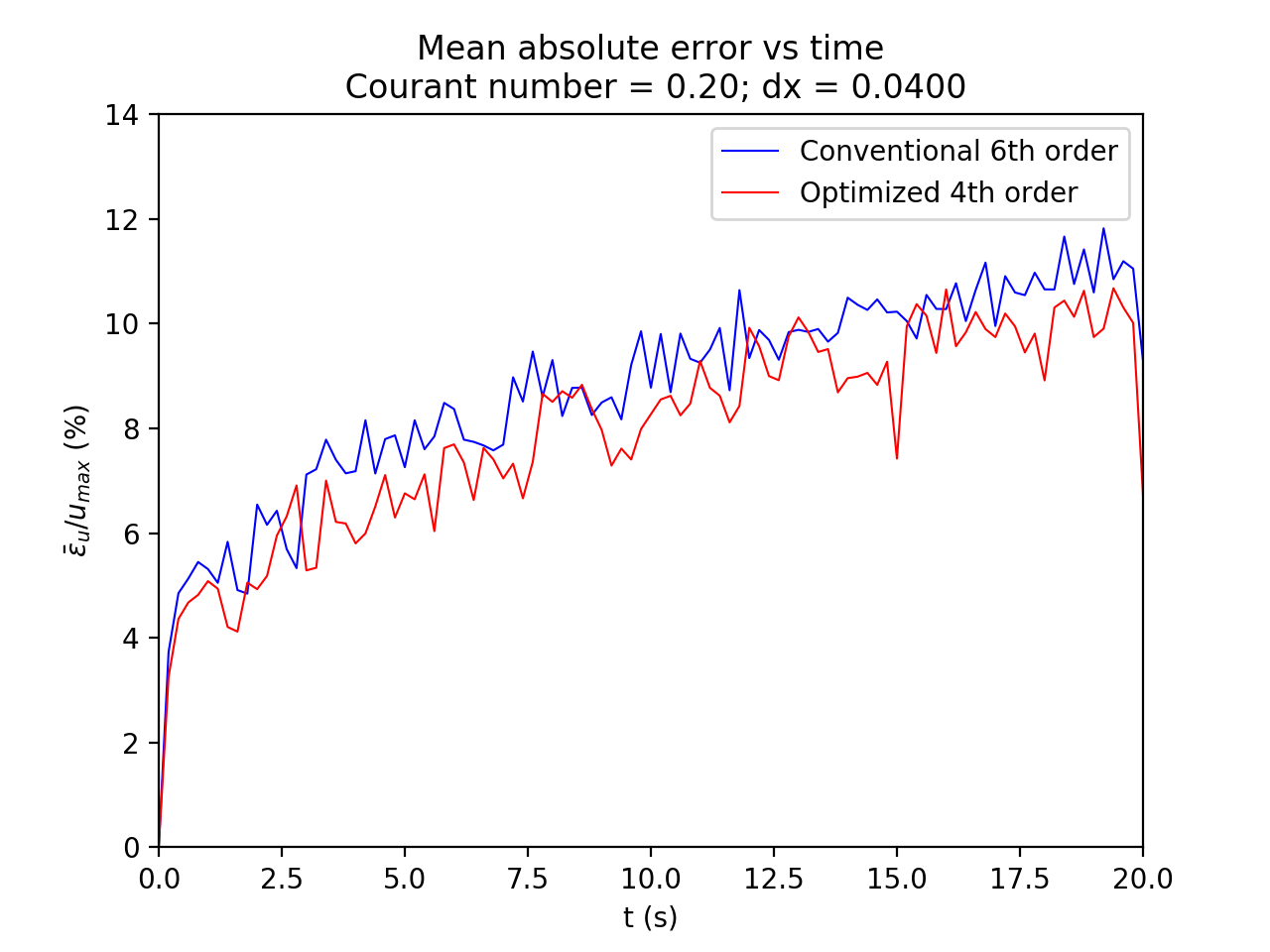}
				\caption{}
			\end{subfigure}
			\caption{Evolution of mean absolute error over time for undersampled wavefields. It is apparent that the spatially-optimized scheme offers consistently better performance for undersampled wavefields, although perhaps not to the degree implied by figure \ref{percent_mean_err_dx}. Errors in this plot are considerably larger than in figure \ref{percent_mean_err_time}, demonstrating that whilst spatial optimization has considerable benefits at large $\Delta x$, the waveform will still be distorted by undersampling.}
			\label{percent_mean_err_under}
		\end{center}
	\end{figure*}
	\begin{figure*}[h!]
		\begin{center}
			\begin{subfigure}{0.75\textwidth}
				\centering\includegraphics[width=\textwidth]{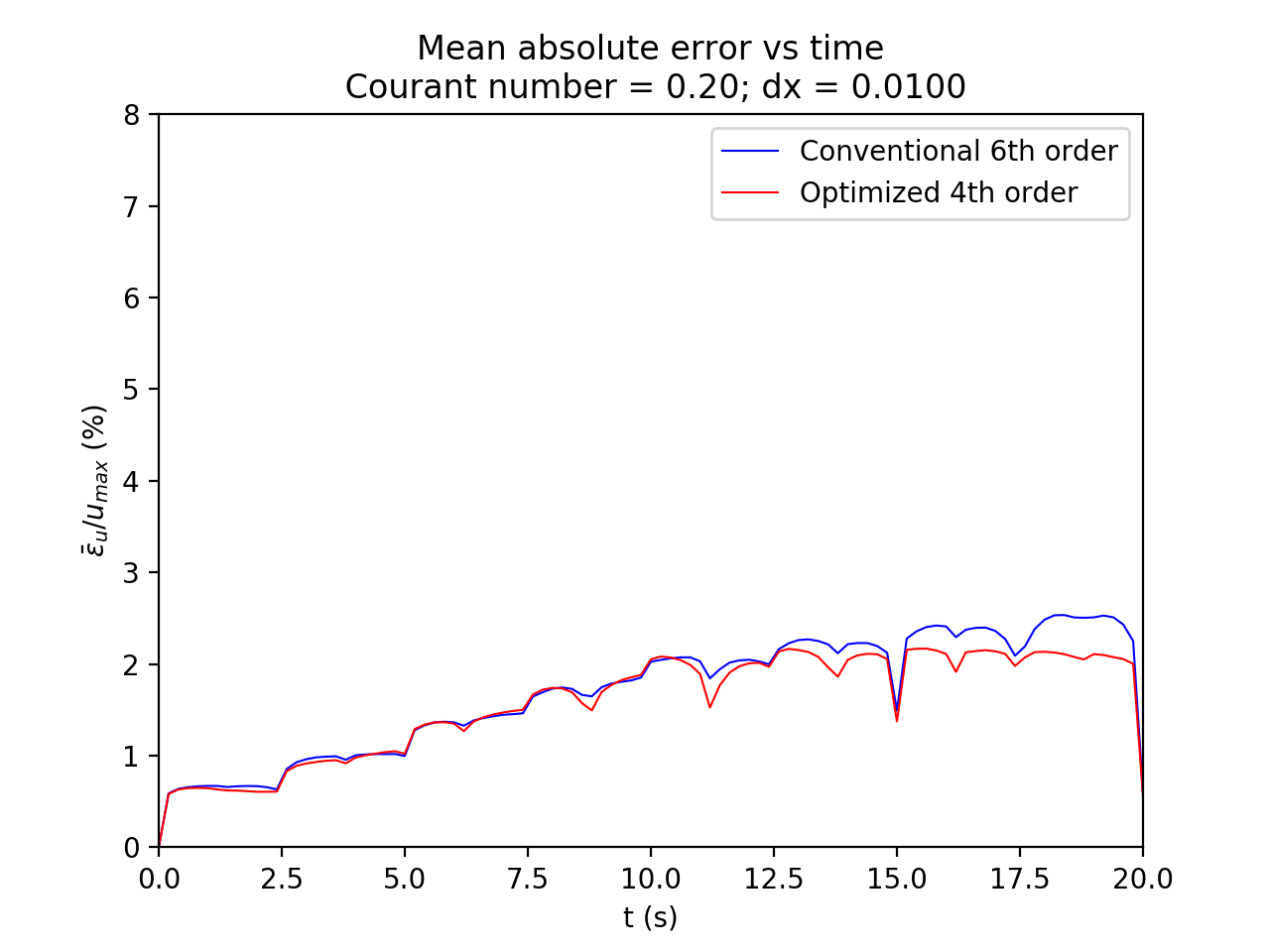}
				\caption{}
			\end{subfigure}
			\begin{subfigure}{0.75\textwidth}
				\centering\includegraphics[width=\textwidth]{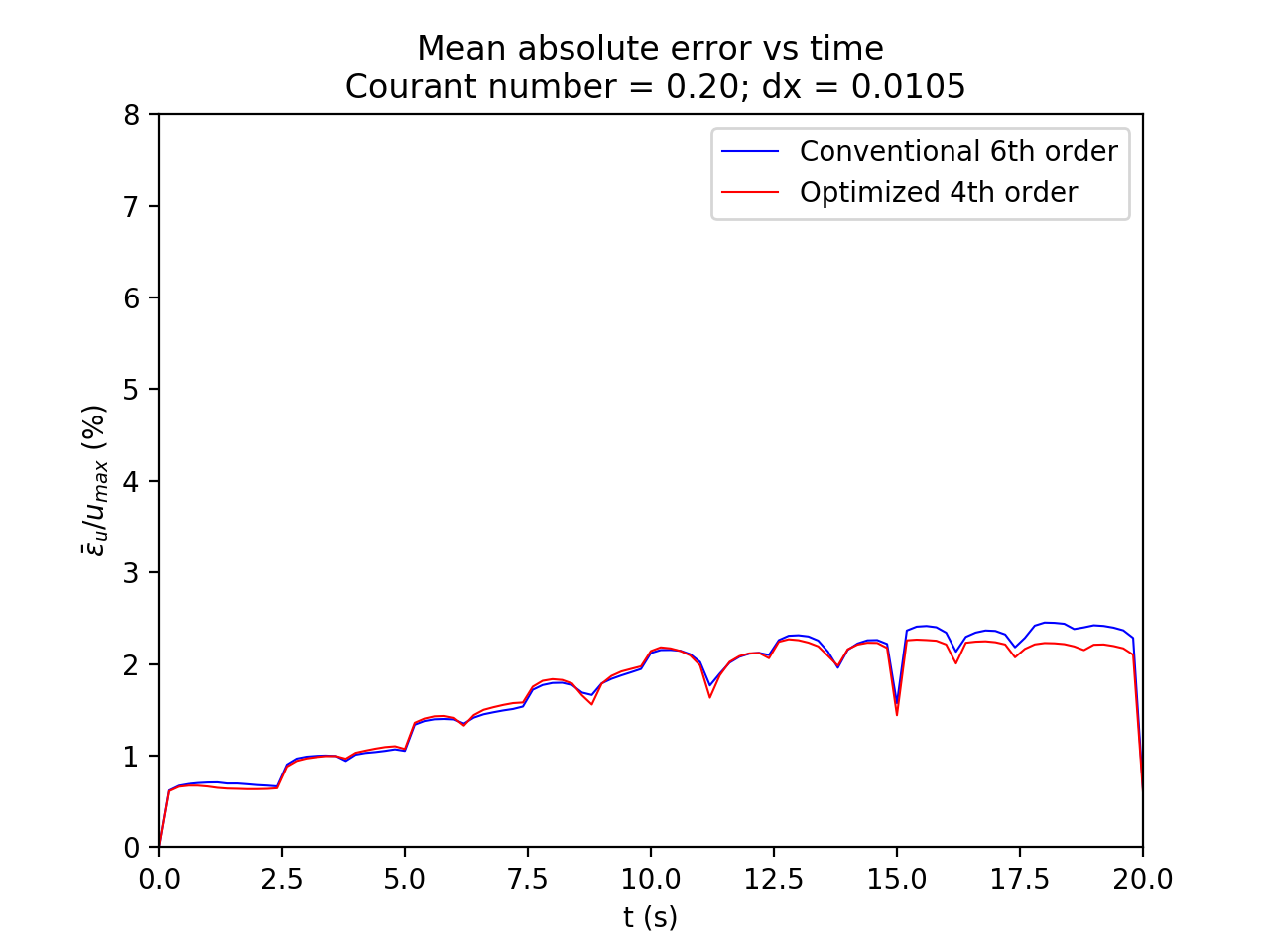}
				\caption{}
			\end{subfigure}
			\caption{Evolution of mean absolute error over time for oversampled wavefields. Oversampled wavefields show minimal error generated by either scheme, with similar error up to approximately 10s, after which error in the optimized scheme does not appear to significantly increase. Error in the conventional scheme continues to accumulate after this point (albeit at a reduced rate), resulting in greater error for the conventional solution at simulation end. This contradicts expectations of superior performance for the conventional scheme at very small $\Delta x$, potentially implying some other error source inherent in conventional FD formulations. Note that error over time for the oversampled wavefields is considerably smoother than for the undersampled wavefields, possibly resulting from dependence of timestep size on $\Delta x$, leading to coarser temporal sampling when spatial sampling is coarse.}
			\label{percent_mean_err_over}
		\end{center}
	\end{figure*}
	\clearpage

	\section{Test case: 2D elastic wave propagation}
	\subsection{Model formulation}
	For seismic modeling, the elastic wave equation is often more useful than the acoustic wave equation, as both P and S waves are supported along with conversion between the two. As the elastic wave equation is a second-order, hyperbolic system, it is commonly substituted with a velocity-stress formulation based on the elastodynamic equations (e.g. \citealp{Virieux1984}, \citealp{Virieux1986}). Additionally, the second-order elastic wave equation suffers from issues arising from an unacceptable degree of dispersion in media with high values of Poisson's ratio (\citealp{Bartolo2015}, \citealp{Stephen1988}). The system can be expressed in terms of particle velocities and stresses to form a series of coupled equations, and is thereby transformed into the first-order, hyperbolic system
	\begin{center}
		\begin{equation}
		\begin{split}
		\frac{\partial v_{x}}{\partial t} = b\left(\frac{\partial\uptau_{xx}}{\partial x}+\frac{\partial\uptau_{xz}}{\partial z}\right) \\ \\
		\frac{\partial v_{z}}{\partial t} = b\left(\frac{\partial\uptau_{xz}}{\partial x}+\frac{\partial\uptau_{zz}}{\partial z}\right) \\ \\
		\frac{\partial\uptau_{xx}}{\partial t}=(\lambda+2\mu)\frac{\partial v_{x}}{\partial x}+\lambda\frac{\partial v_{z}}{\partial z} \\ \\
		\frac{\partial\uptau_{zz}}{\partial t}=(\lambda+2\mu)\frac{\partial v_{z}}{\partial z}+\lambda\frac{\partial v_{x}}{\partial x} \\ \\
		\frac{\partial\uptau_{xz}}{\partial t}=\mu\left(\frac{\partial v_{x}}{\partial z}+\frac{\partial v_{z}}{\partial x}\right)
		\end{split}
		\label{eq:22}
		\end{equation}	
	\end{center} 
	where $(v_{x},v_{y})$ is particle velocity, $\uptau_{xx}$ and $\uptau_{zz}$ are normal stresses, and $\uptau_{xz}$ is shear stress. First and second Lam\'e parameters are referred to as $\lambda$ and $\mu$ respectively, and $b(x,z)$ is buoyancy: density's inverse.
	\\
	\\To improve precision and stability, it is commonplace to use a staggered grid for elastic modeling applications (\citealp{Liu2009}), evaluating derivatives on grids offset by $\frac{\Delta x}{2}$ from the reference grid. Velocity components are discretized on citealpate grids with shear stress discretized at a citealpate node to normal stresses in accordance with \citealp{Virieux1986}. Use of non-staggered grids with the P-SV formulation produces artifacts of similar amplitude to the data, and is thus inadvisable.
	\\
	\\To prove the efficacy of spatially-optimized stencils for seismic modeling applications, optimized first derivative stencils for staggered grids have been derived for this system. To test feasibility of the method in practical applications, a 2 layer model consisting of water overlying shale, evocative of offshore surveys, was used. A 2km by 2km computational domain with $\Delta x=\Delta z=10\mathrm{m}$, divided into two layers each 1km thick was used. Layers were defined in terms of P-wave velocity, S-wave velocity, and density, from which Lam\'e parameters were obtained. The upper water layer had a P-wave velocity of 1400ms\textsuperscript{-1} and a density of 1000kgm\textsuperscript{-3}. As fluids do not support shearing, S-wave velocity was 0ms\textsuperscript{-1}. The lower shale layer had a P-wave velocity of 4000ms\textsuperscript{-1}, S-wave velocity of 2400ms\textsuperscript{-1} and a density of 2600kgm\textsuperscript{-3}. These values were informed by measurements from \citealp{Castagna1985}, \citealp{Hamilton1978}, and \citealp{Greenspan1959} to ensure results analogous to real-world phenomena. The source function was a Ricker wavelet with a peak frequency of 11.2Hz, coinciding with the theoretical maximum frequency which can be propagated accurately by both the conventional and optimized staggered schemes (see figure \ref{1st_dispersion_plot_staggered}). Note that as this frequency choice was based on S-wave velocity, longer-wavelength P-waves should be subject to reduced dispersion (\citealp{Moczo2000}). This source was injected into both normal stress fields at a centered position, 750m from the top of the domain. This position was chosen such that incidence of the wave upon the layer boundary occurred prior to reflections from domain boundaries. External boundaries were kept as free surfaces to simplify implementation. A Courant number of 0.47 was used in accordance with stability conditions specified in \citealp{Virieux1986} and a P-wave velocity of 6000ms\textsuperscript{-1}: well in excess of the model maximum, ensuring a good-quality solution.

	\subsection{Derivation of stencil coefficients}
\begin{figure*}[h!]
		\begin{center}
			\includegraphics[angle=0, width=0.75\textwidth, keepaspectratio=true]{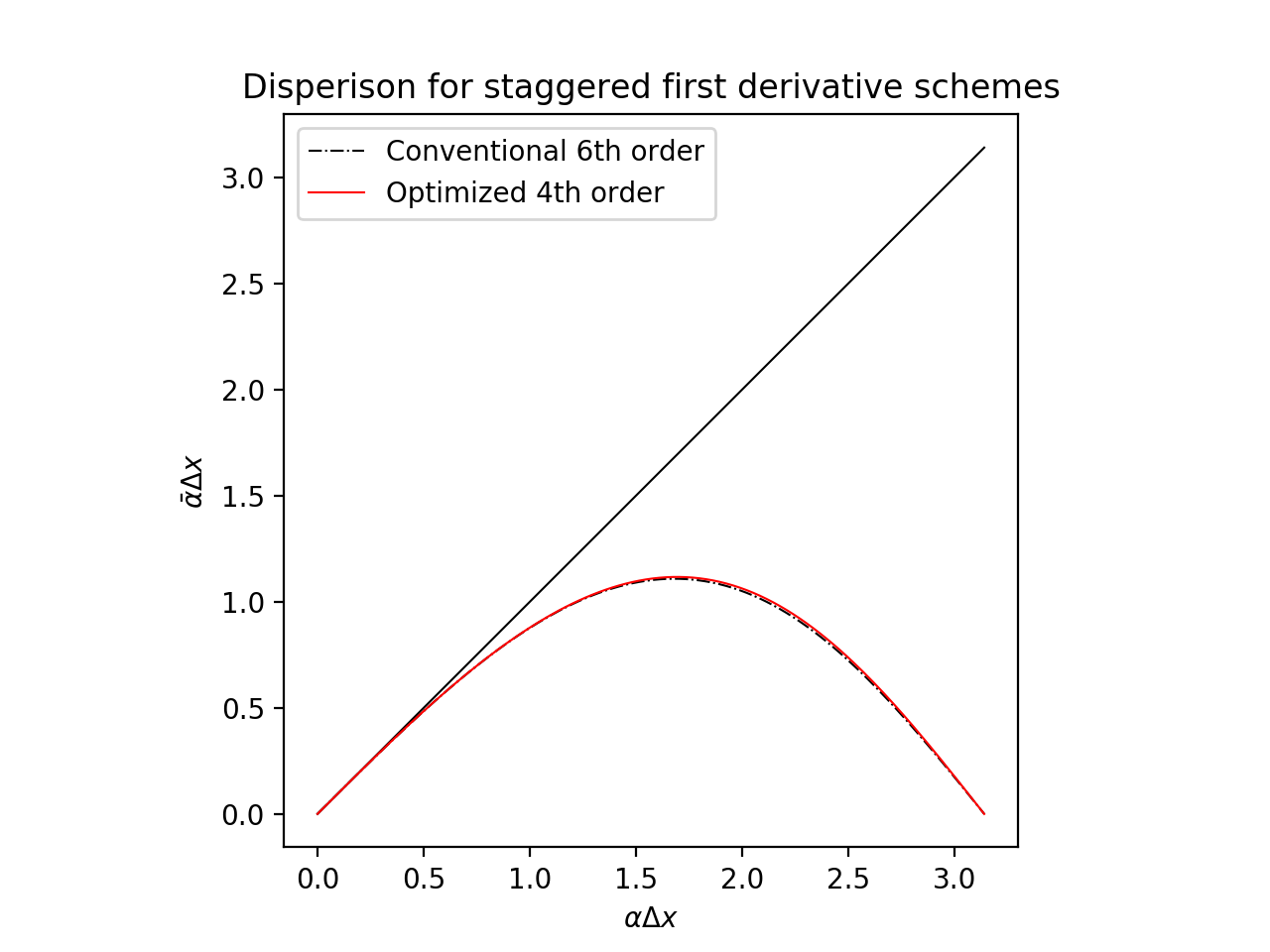} 
			\caption{$\bar{\alpha}\Delta x$ versus $\alpha\Delta x$ for the conventional 6\textsuperscript{th} order staggered and spatially-optimized 4\textsuperscript{th} order staggered central-difference scheme. 	Note the only marginally better performance of the optimized scheme in terms of dispersion-relation-preservation. This is unsurprising, as the two have near-identical coefficients; the staggered conventional scheme is approximately optimized. Interestingly, both schemes deviate from the ideal at relatively small values of $\alpha\Delta x$ compared to second derivative schemes shown in figure \ref{2nd_dispersion_plot}. The conventional and optimized staggered schemes both require highest-frequency component wavelengths above $12.6\Delta x$. For shorter wavelengths, the optimized scheme will slightly better approximate propagation characteristics of a wave.}
			\label{1st_dispersion_plot_staggered}
		\end{center}
	\end{figure*}

	The 4\textsuperscript{th} order spatially-optimized staggered first derivative has two stencil variants, $M=3,\: N=2$ and $M=2,\: N=3$, referenced henceforth as types ``A" and ``B" respectively. The former evaluates a derivative offset by half a grid spacing in the positive direction, whilst the latter yields the derivative half a grid spacing in the negative direction. For these stencils, equation \ref{eq:3} becomes
	\begin{center}	
		\begin{equation}
		\frac{\partial u}{\partial x}\approx\frac{1}{\Delta x}\sum_{j=-N}^{M}a_{j}u\left(x+\left(j-\frac{1}{2}\right)\Delta x\right)
		\label{eq:23}
		\end{equation}
	\end{center}
	and
	\begin{center}	
		\begin{equation}
		\frac{\partial u}{\partial x}\approx\frac{1}{\Delta x}\sum_{j=-N}^{M}a_{j}u\left(x+\left(j+\frac{1}{2}\right)\Delta x\right)
		\label{eq:24}
		\end{equation}
	\end{center}
	respectively. Values of $E$ are thus expressed as
	\begin{center}
		\begin{equation}
		E_{A}=\int_{-\frac{\pi}{2}}^{\frac{\pi}{2}}\left|i\kappa+\sum_{j=-N}^{M}a_{j}e^{i\left(j-\frac{1}{2}\right)\kappa}\right|^{2}\mathrm{d}\kappa
		\label{eq:25}
		\end{equation}	
	\end{center}
	and
	\begin{center}
		\begin{equation}
		E_{B}=\int_{-\frac{\pi}{2}}^{\frac{\pi}{2}}\left|i\kappa+\sum_{j=-N}^{M}a_{j}e^{i\left(j+\frac{1}{2}\right)\kappa}\right|^{2}\mathrm{d}\kappa
		\label{eq:26}
		\end{equation}	
	\end{center}
	Coefficients of the type ``A" stencil are constrained by
	\begin{center}
		\begin{equation}
		\begin{split}
		a_{0}=-a_{1} \\ \\
		a_{2}=-a_{-1}=\frac{25}{48}-\frac{a_{1}}{2} \\ \\
		a_{3}=-a_{-2}=\frac{a_{1}}{10}-\frac{9}{80}
		\end{split}
		\label{eq:27}
		\end{equation}	
	\end{center}
	whilst constraints for the type ``B" stencil are
	\begin{center}
		\begin{equation}
		\begin{split}
		a_{-1}=-a_{0} \\ \\
		a_{1}=-a_{-2}=\frac{25}{48}-\frac{a_{0}}{2} \\ \\
		a_{2}=-a_{-3}=\frac{a_{0}}{10}-\frac{9}{80}
		\end{split}
		\label{eq:28}
		\end{equation}	
	\end{center}
	When minimization is carried out for the two stencils, the coefficients
	\begin{center}
		\begin{equation}
		\begin{split}
		a_{0}=-a_{1}=-1.1890097 \\ \\
		a_{2}=-a_{-1}=-0.07367152 \\ \\
		a_{3}=-a_{-2}=0.00640097
		\end{split}
		\label{eq:29}
		\end{equation}	
	\end{center}
	and
	\begin{center}
		\begin{equation}
		\begin{split}
		a_{-1}=-a_{0}=-1.1890097 \\ \\
		a_{1}=-a_{-2}=-0.07367152 \\ \\
		a_{2}=-a_{-3}=0.00640097
		\end{split}
		\label{eq:30}
		\end{equation}	
	\end{center}
	are obtained for stencils ``A" and ``B" respectively. Note that these two stencils have identical coefficients shifted by a single grid point. Consequentially, both stencils are assumed to have identical dispersion characteristics. 
	\\
	\\Dispersion curves for conventional 6\textsuperscript{th} order staggered and spatially-optimized 4\textsuperscript{th} order staggered schemes are shown in figure \ref{1st_dispersion_plot_staggered}.
	\clearpage

	\subsection{Performance at near-maximum grid spacings}
	\begin{figure*}[h!]
		\begin{center}
			\begin{subfigure}{0.45\textwidth}
				\centering\includegraphics[width=\textwidth]{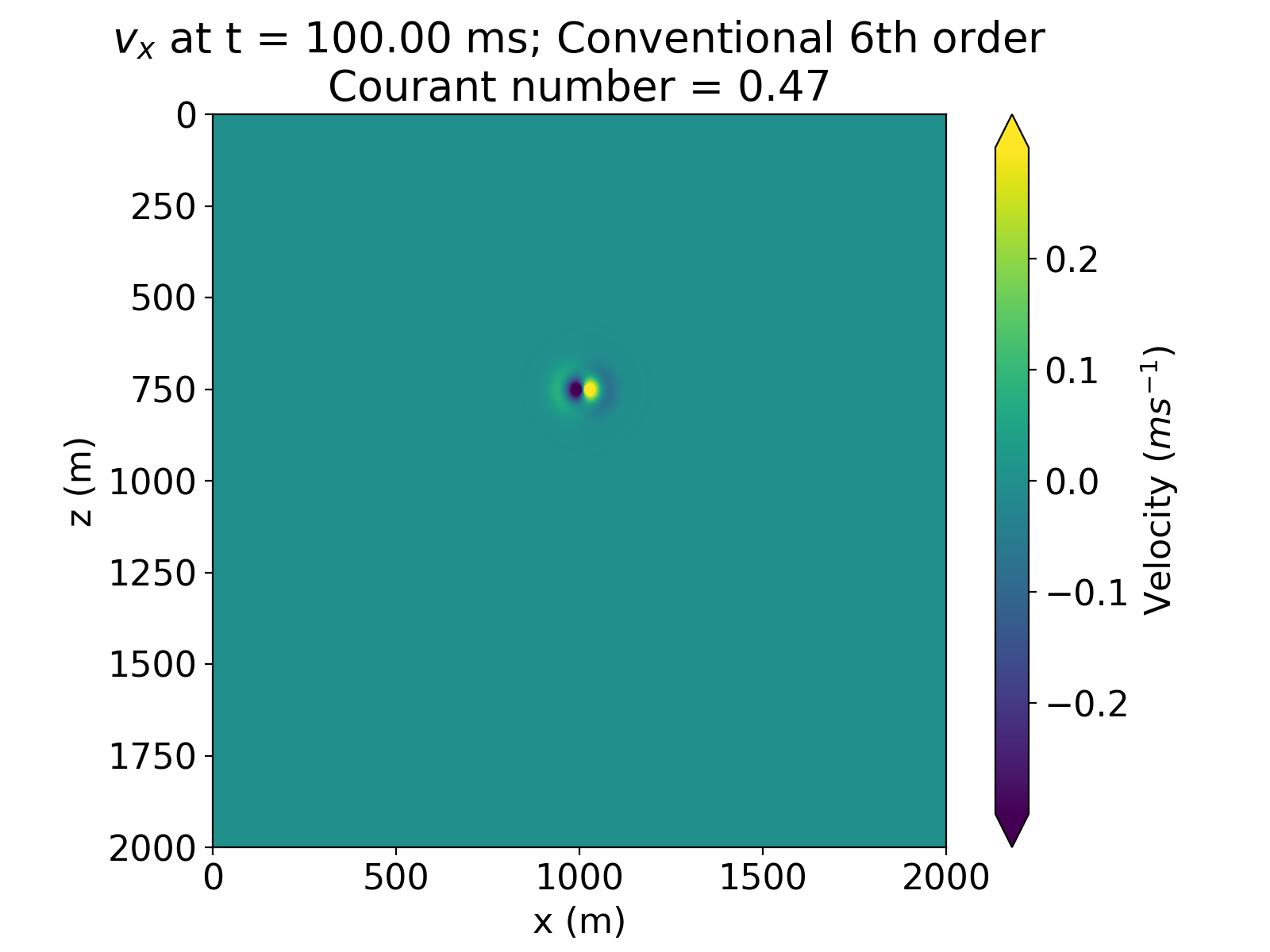}
				\caption{}
			\end{subfigure}
			\begin{subfigure}{0.45\textwidth}
				\centering\includegraphics[width=\textwidth]{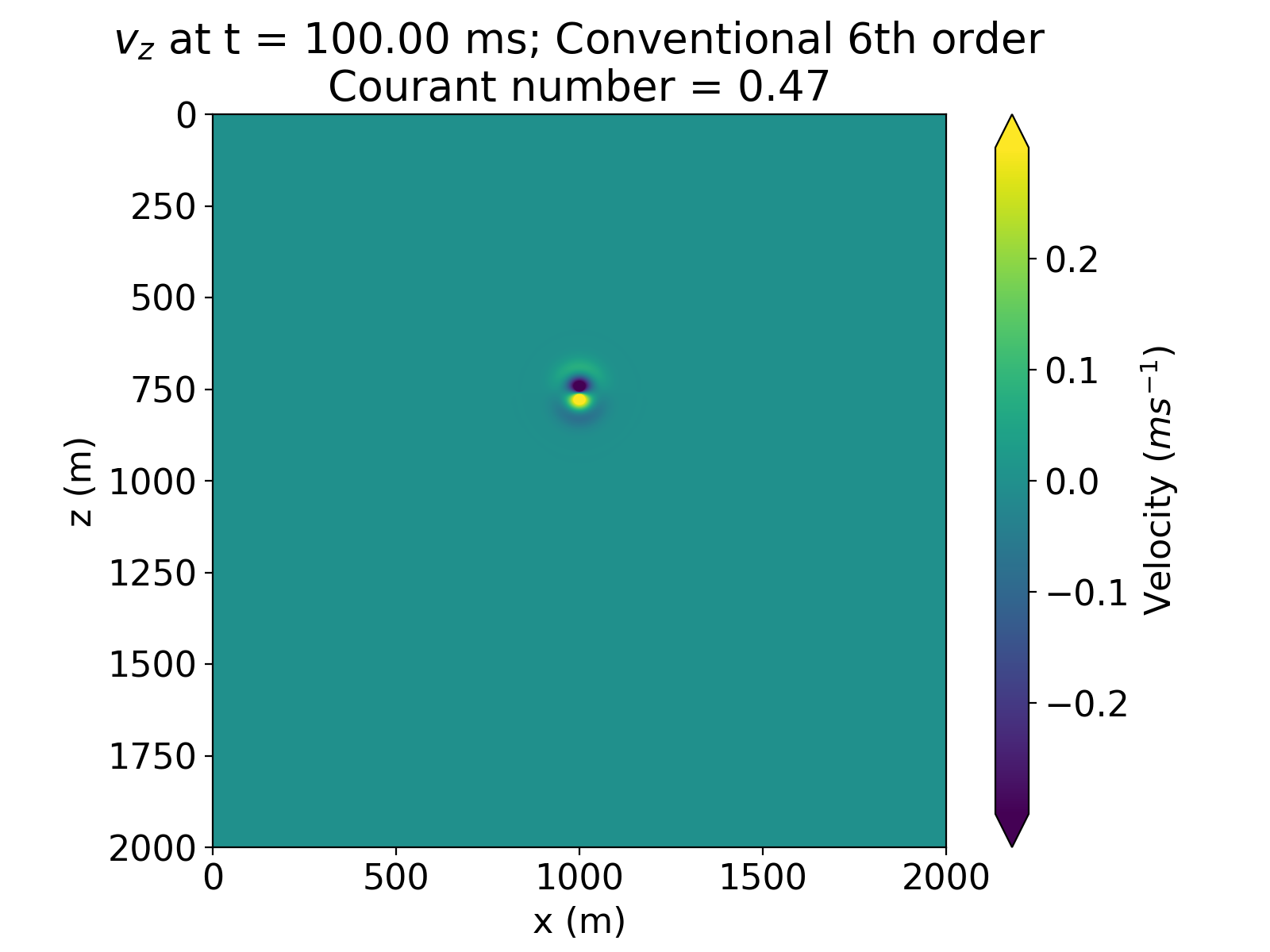}
				\caption{}
			\end{subfigure}\\
			\begin{subfigure}{0.45\textwidth}
				\centering\includegraphics[width=\textwidth]{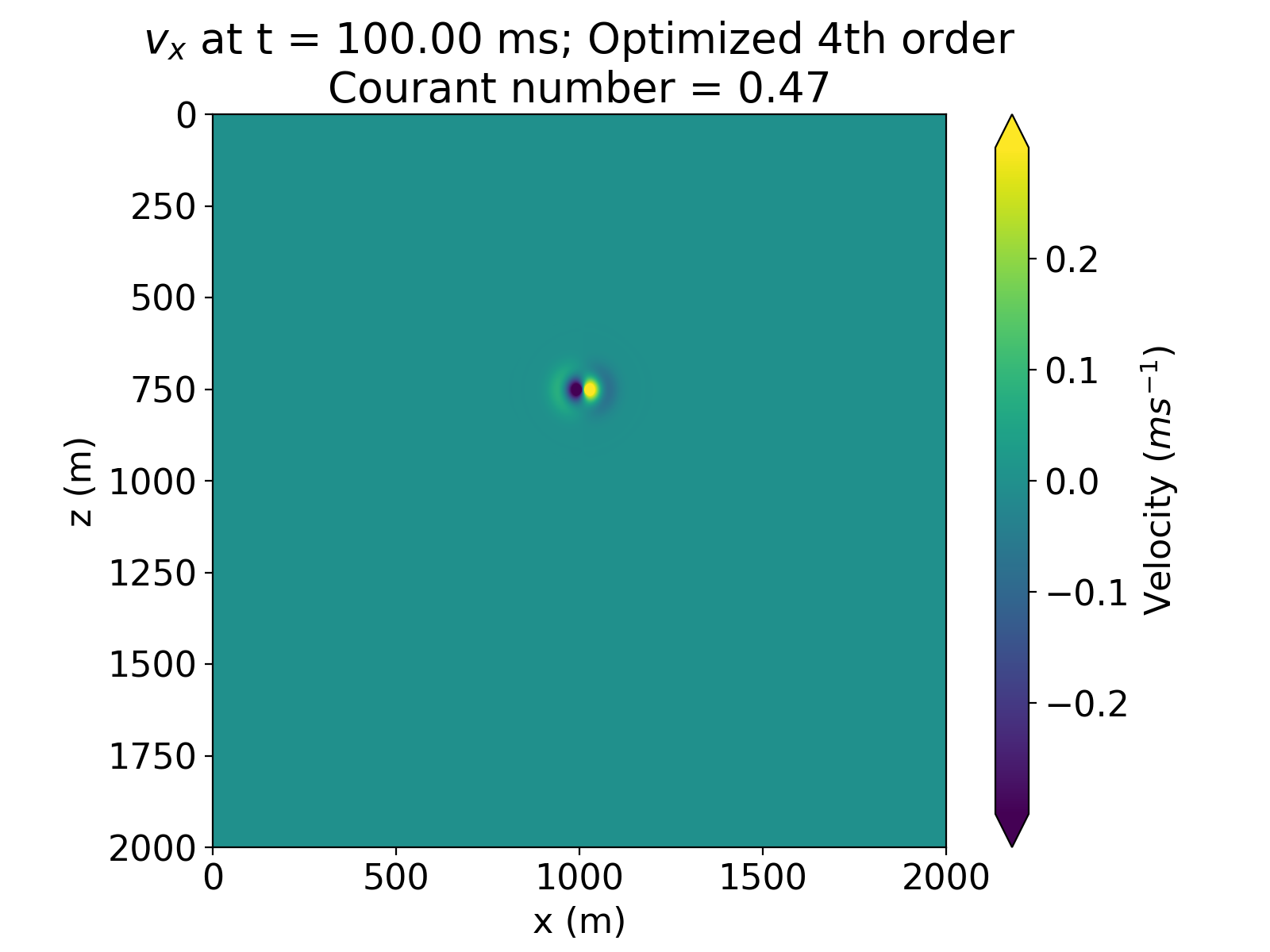}
				\caption{}
			\end{subfigure}
			\begin{subfigure}{0.45\textwidth}
				\centering\includegraphics[width=\textwidth]{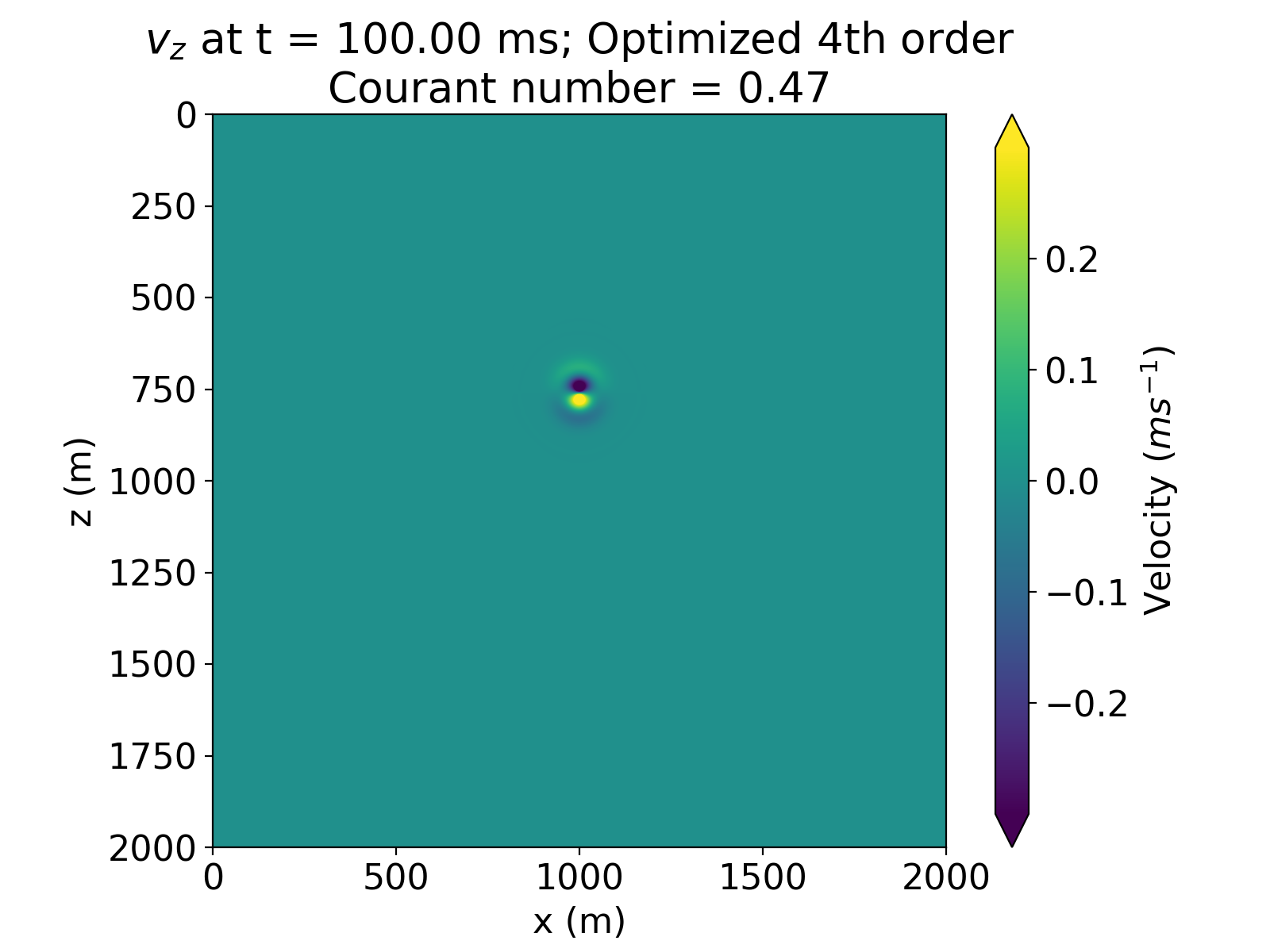}
				\caption{}
			\end{subfigure}
			\begin{subfigure}{0.45\textwidth}
				\centering\includegraphics[width=\textwidth]{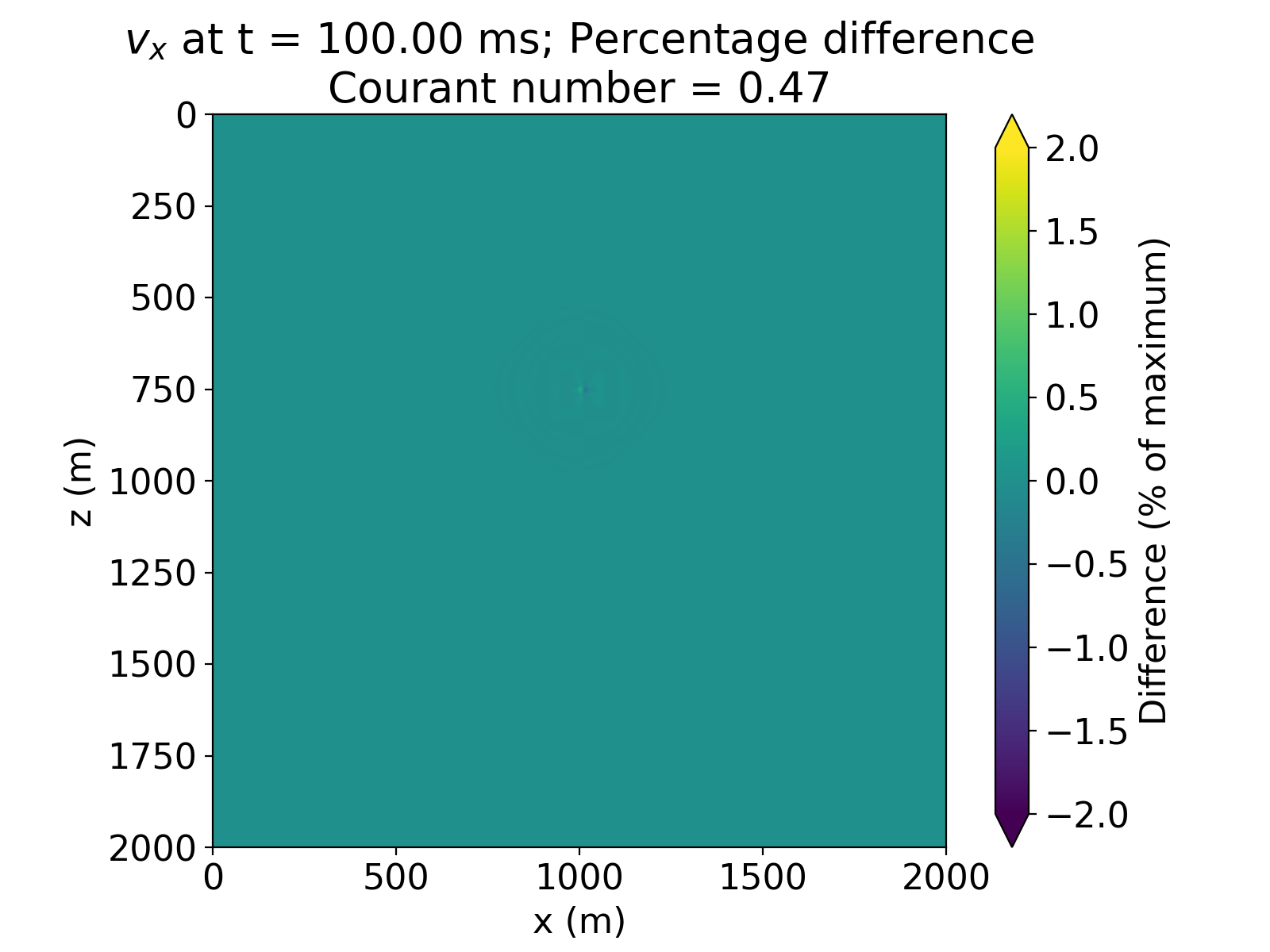}
				\caption{}
			\end{subfigure}
			\begin{subfigure}{0.45\textwidth}
				\centering\includegraphics[width=\textwidth]{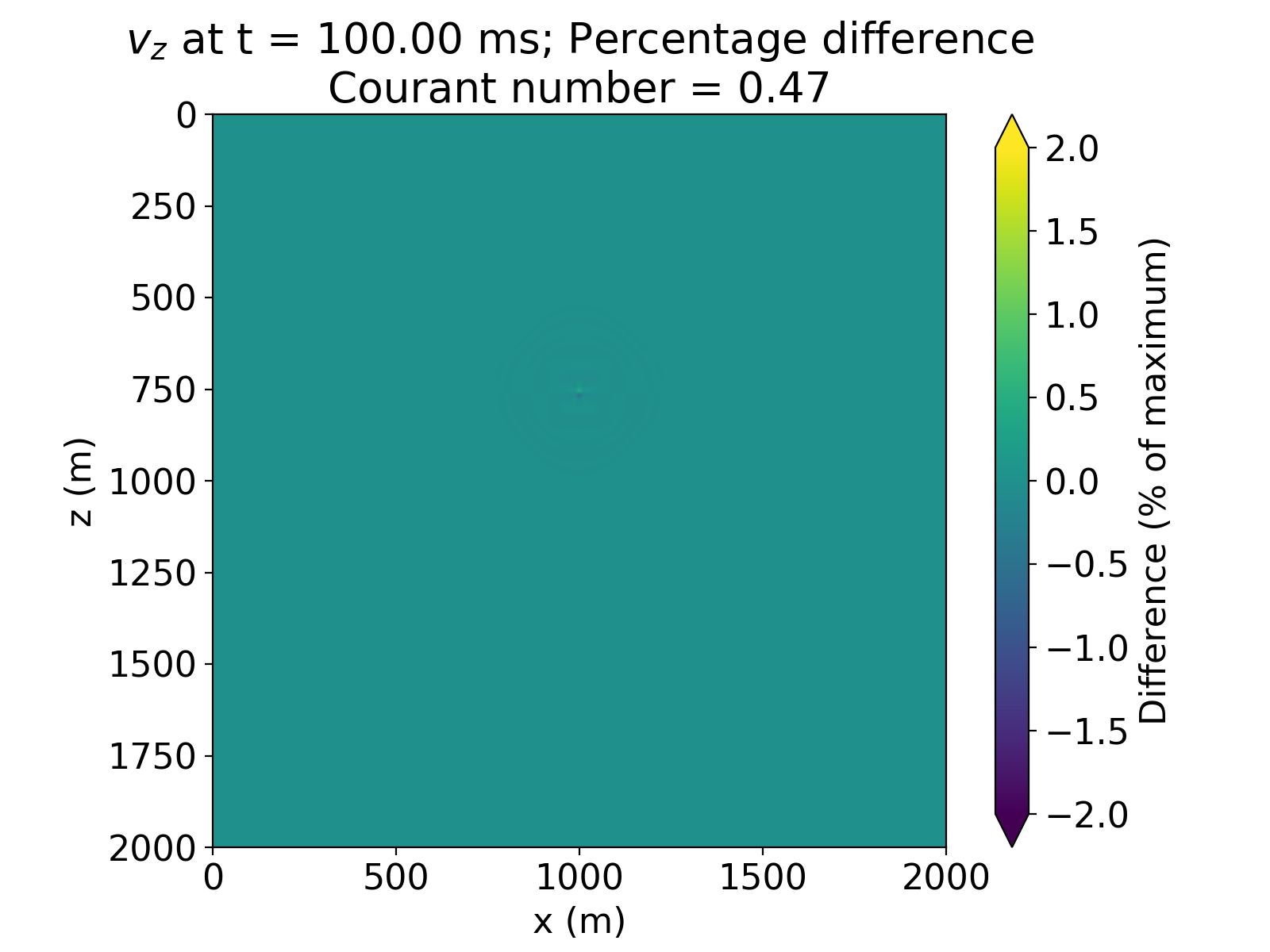}
				\caption{}
			\end{subfigure}
			\caption{Wavefields for velocity components $v_{x}$ and $v_{z}$ calculated using conventional and optimized schemes at 100ms. Difference is normalized against the largest amplitude present in the wavefield.}
			\label{2D_elastic_100}
		\end{center}
	\end{figure*}
	\begin{figure*}[h!]
		\begin{center}
			\begin{subfigure}{0.45\textwidth}
				\centering\includegraphics[width=\textwidth]{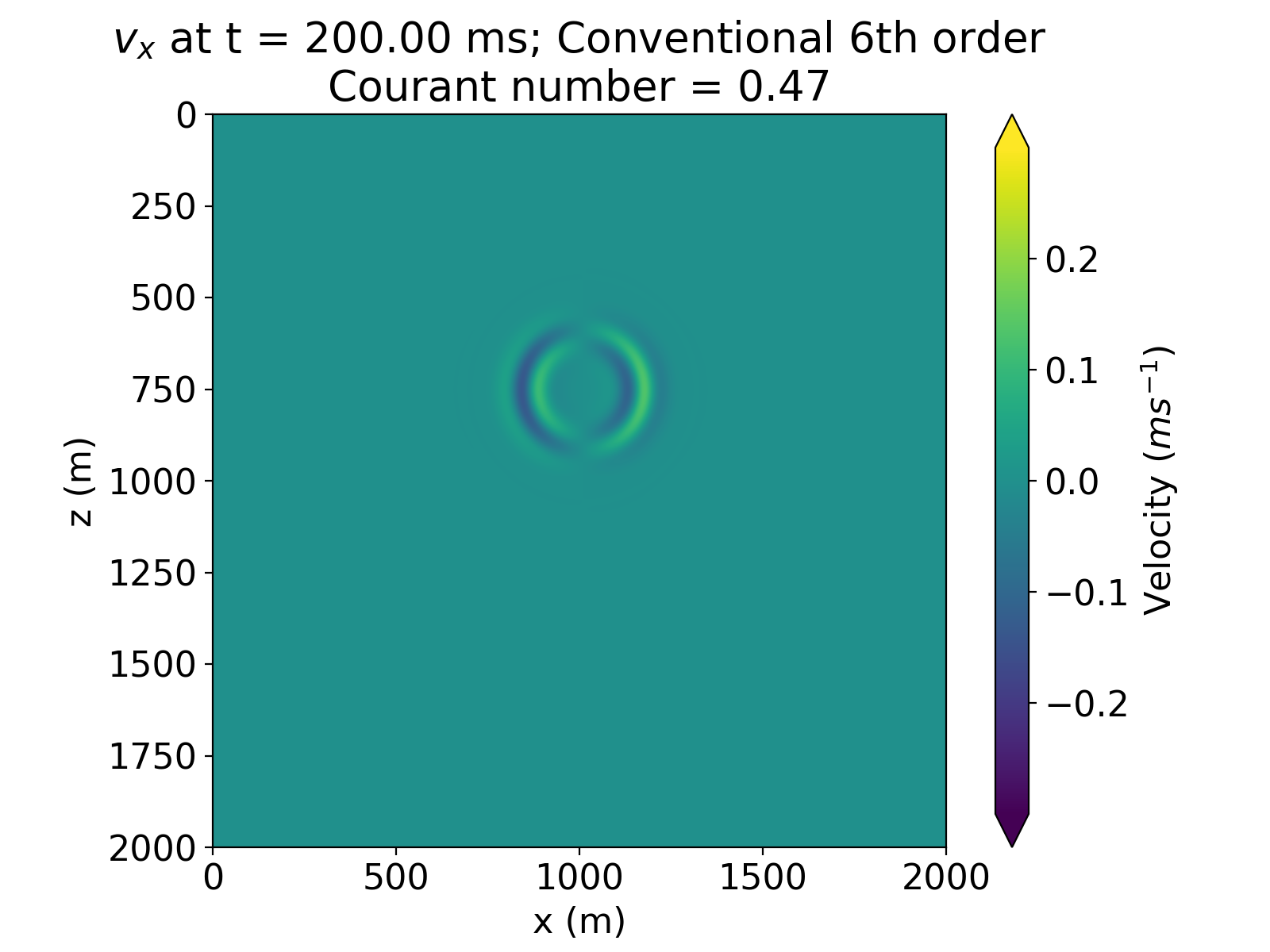}
				\caption{}
			\end{subfigure}
			\begin{subfigure}{0.45\textwidth}
				\centering\includegraphics[width=\textwidth]{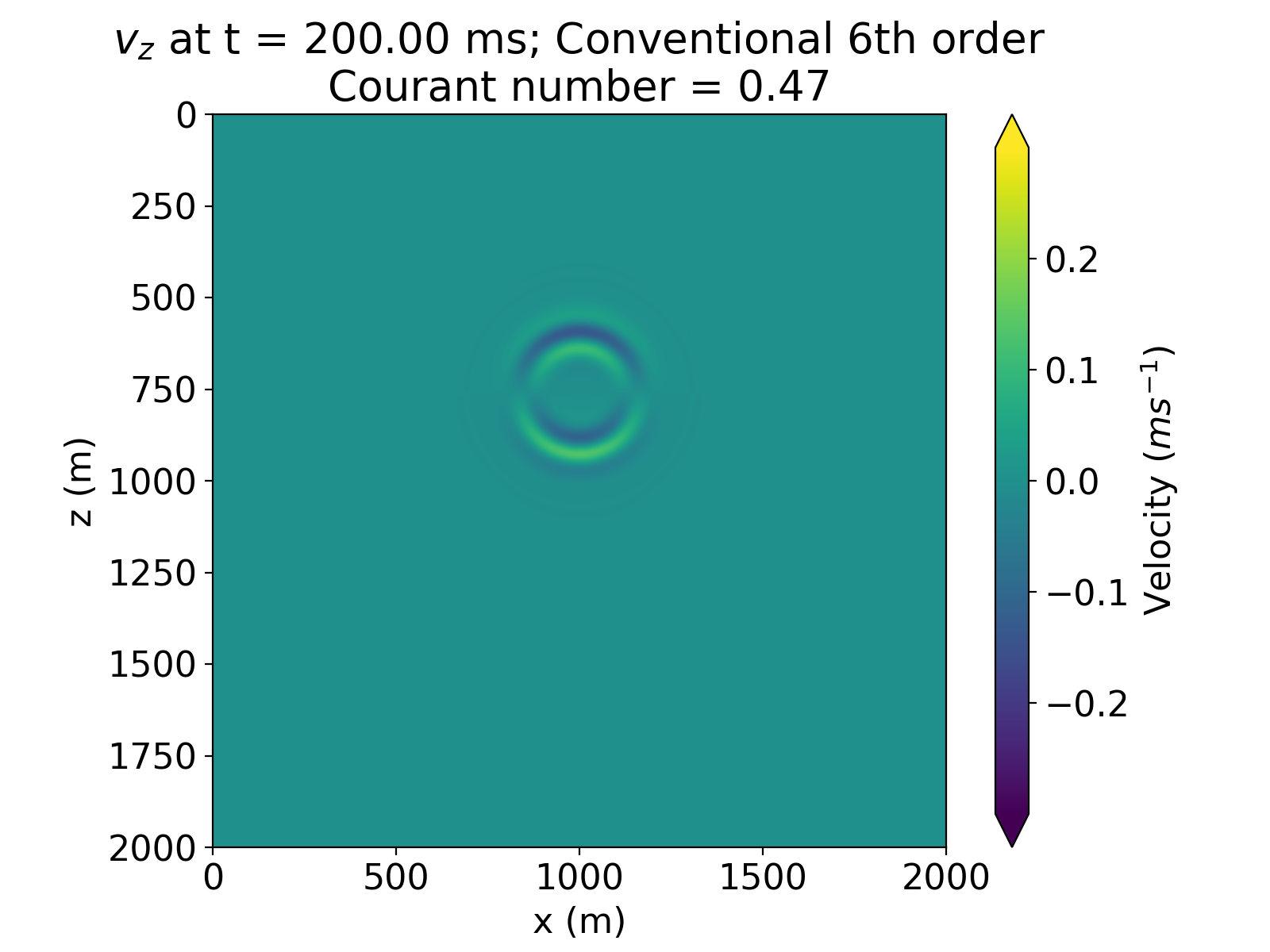}
				\caption{}
			\end{subfigure}\\
			\begin{subfigure}{0.45\textwidth}
				\centering\includegraphics[width=\textwidth]{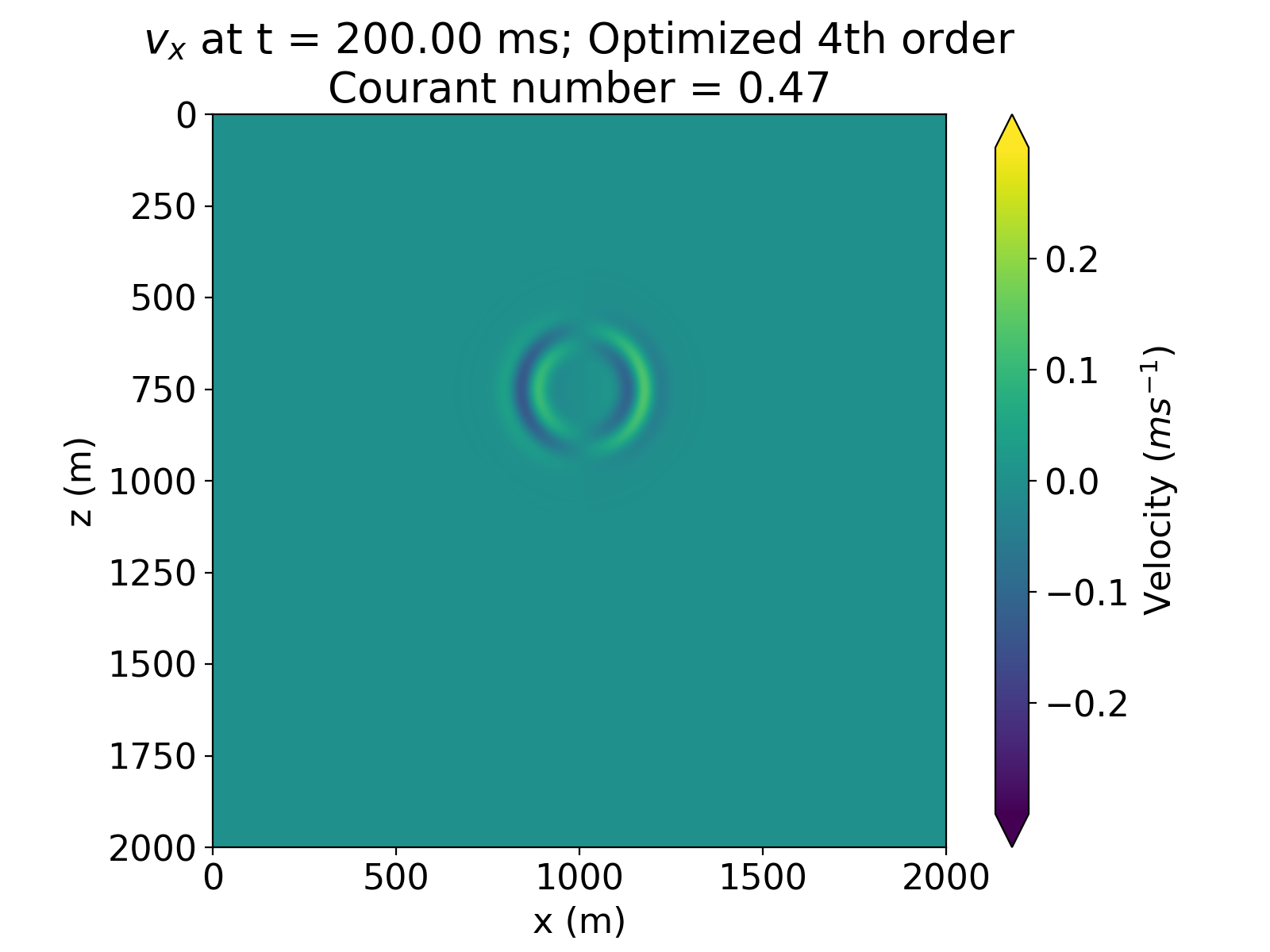}
				\caption{}
			\end{subfigure}
			\begin{subfigure}{0.45\textwidth}
				\centering\includegraphics[width=\textwidth]{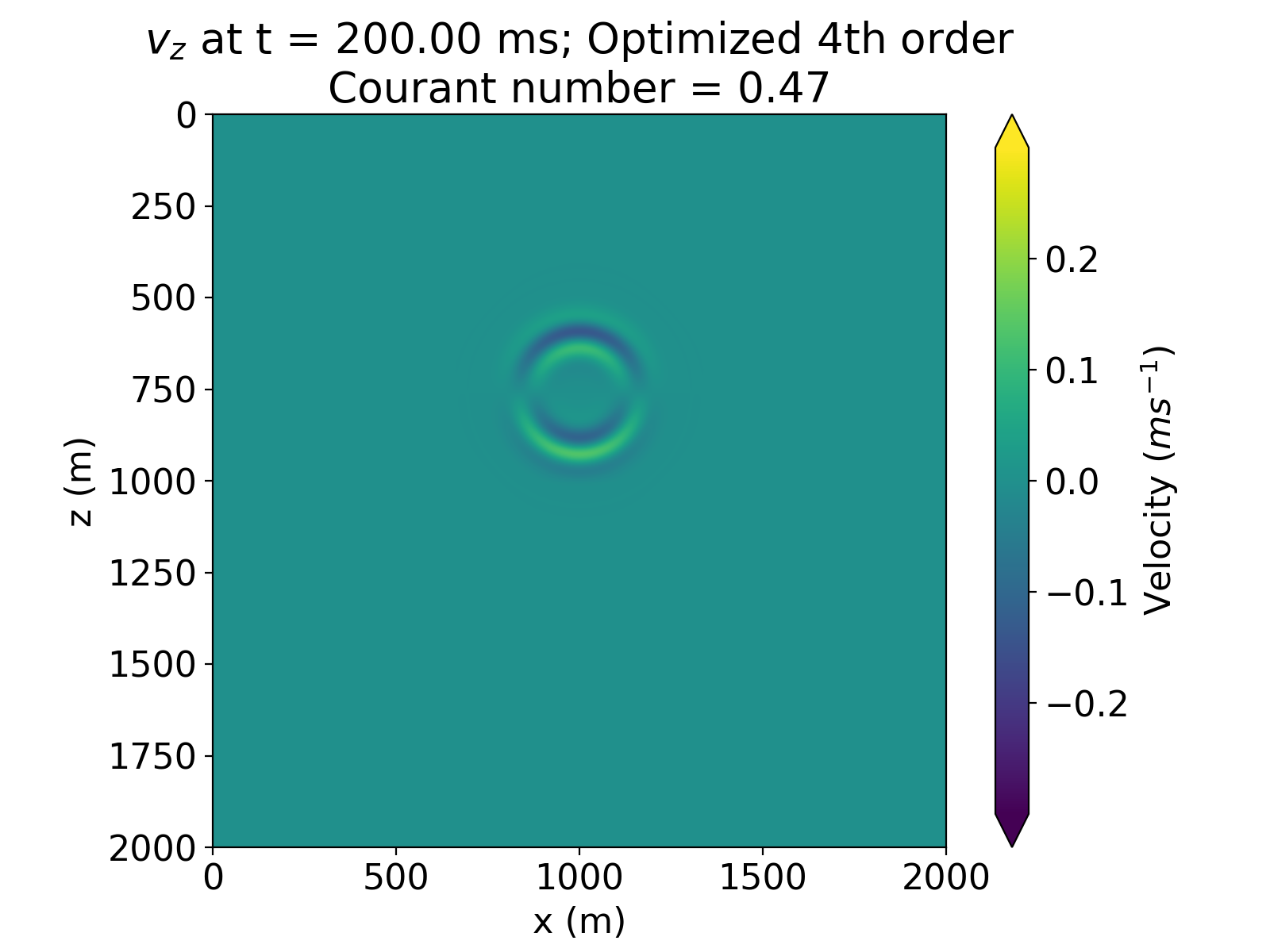}
				\caption{}
			\end{subfigure}
			\begin{subfigure}{0.45\textwidth}
				\centering\includegraphics[width=\textwidth]{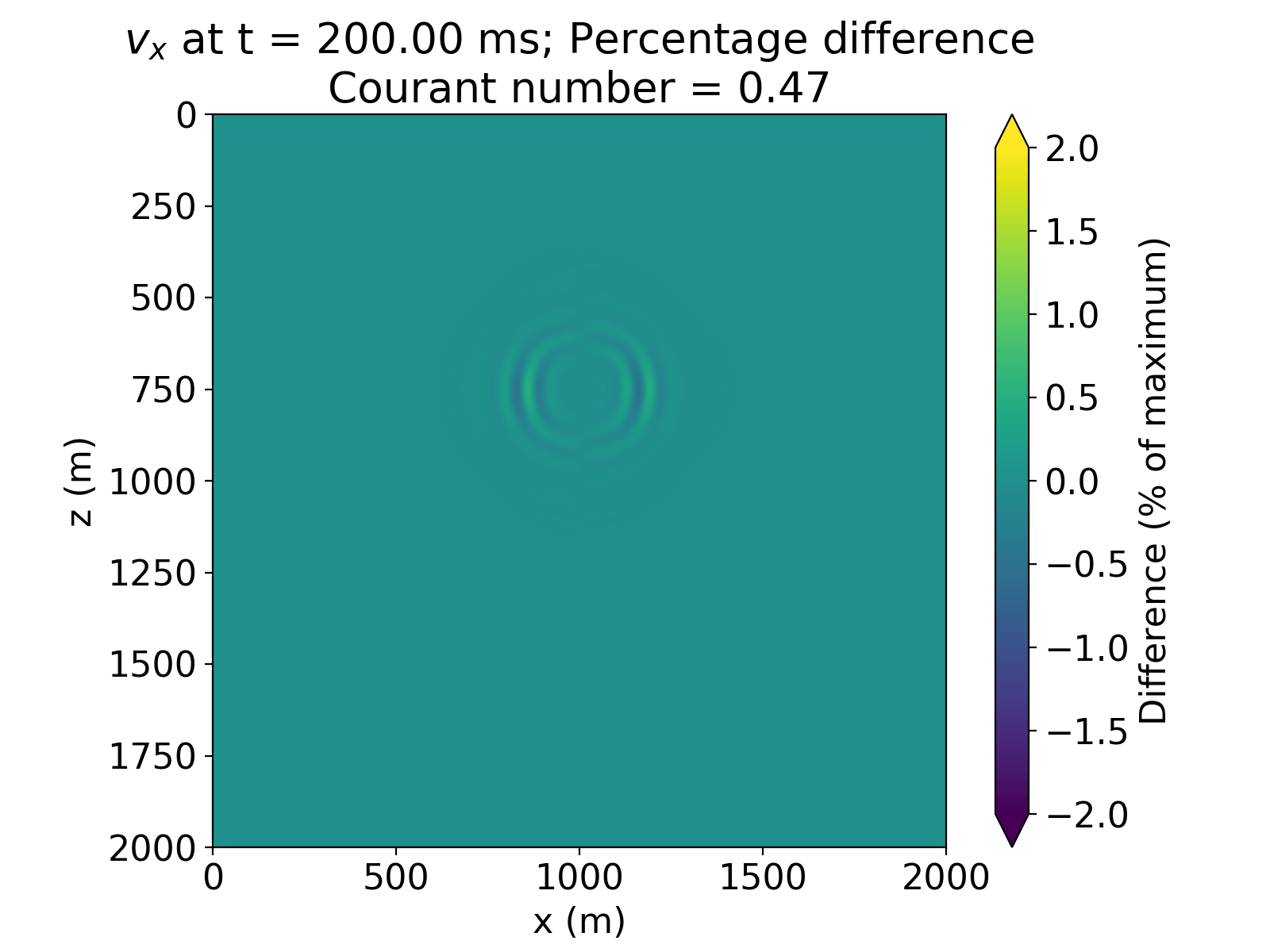}
				\caption{}
			\end{subfigure}
			\begin{subfigure}{0.45\textwidth}
				\centering\includegraphics[width=\textwidth]{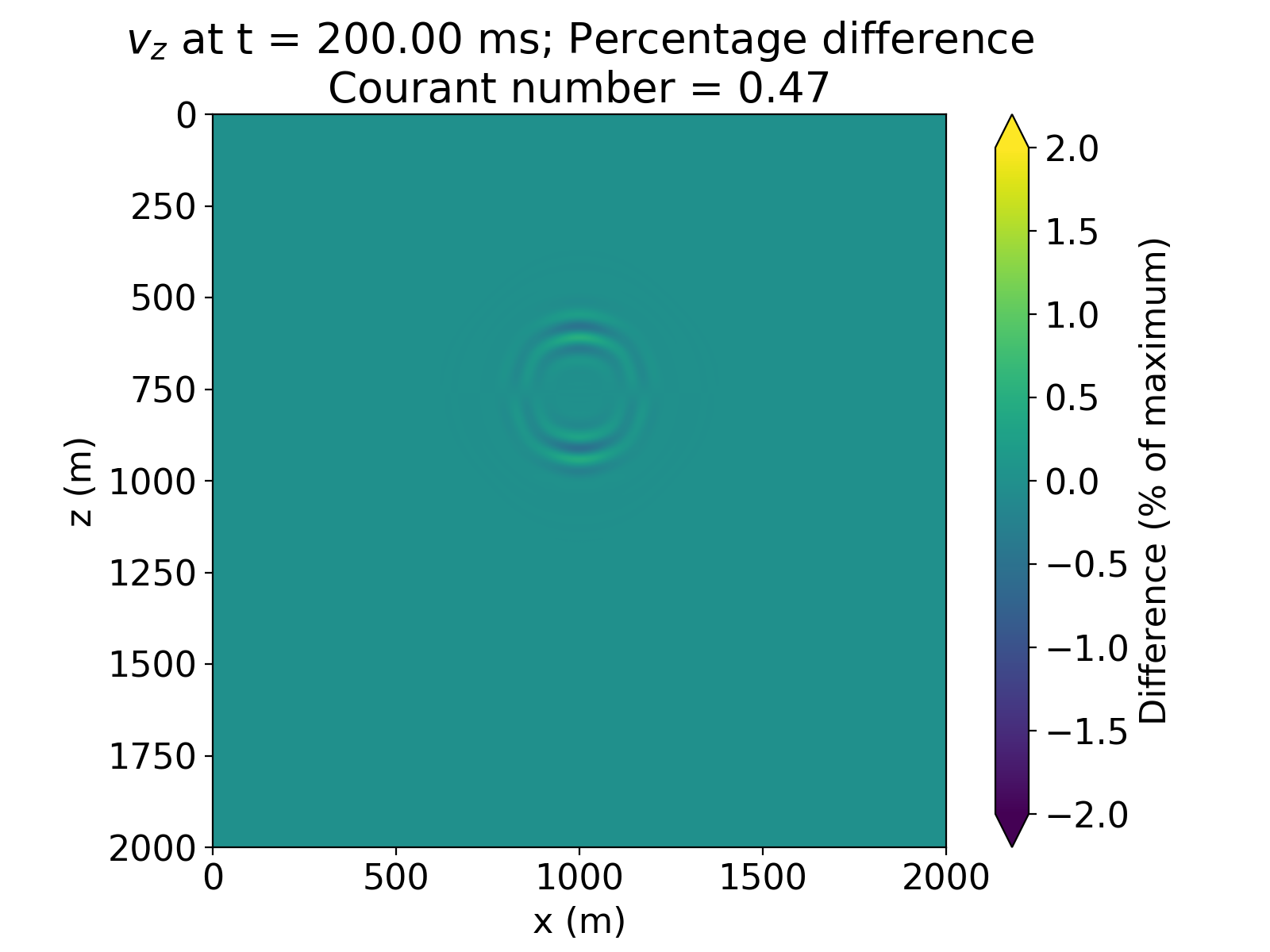}
				\caption{}
			\end{subfigure}
			\caption{Wavefields for velocity components $v_{x}$ and $v_{z}$ calculated using conventional and optimized schemes at 200ms. Difference is normalized against the largest amplitude present in the wavefield.}
			\label{2D_elastic_200}
		\end{center}
	\end{figure*}
	\begin{figure*}[h!]
		\begin{center}
			\begin{subfigure}{0.45\textwidth}
				\centering\includegraphics[width=\textwidth]{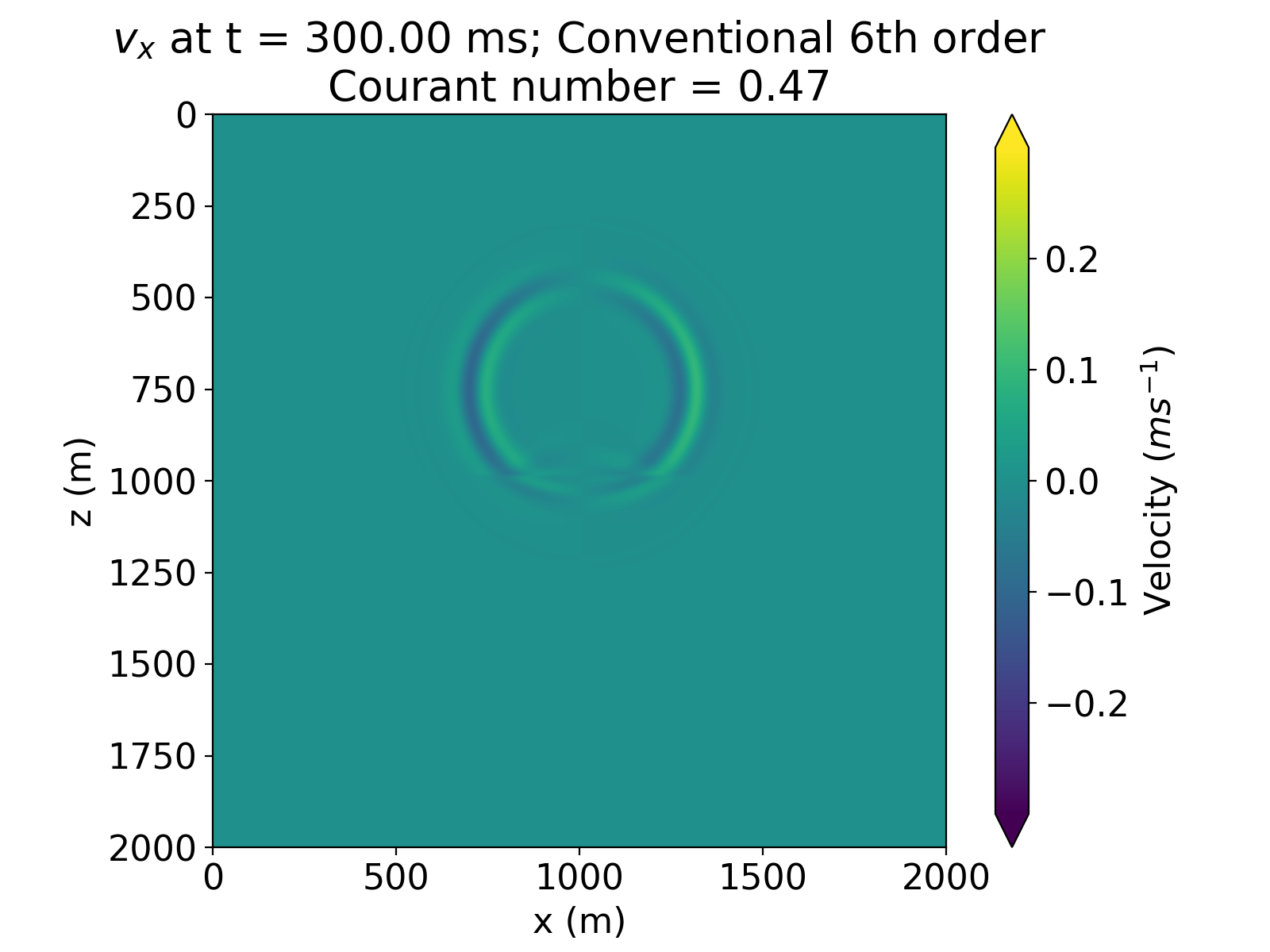}
				\caption{}
			\end{subfigure}
			\begin{subfigure}{0.45\textwidth}
				\centering\includegraphics[width=\textwidth]{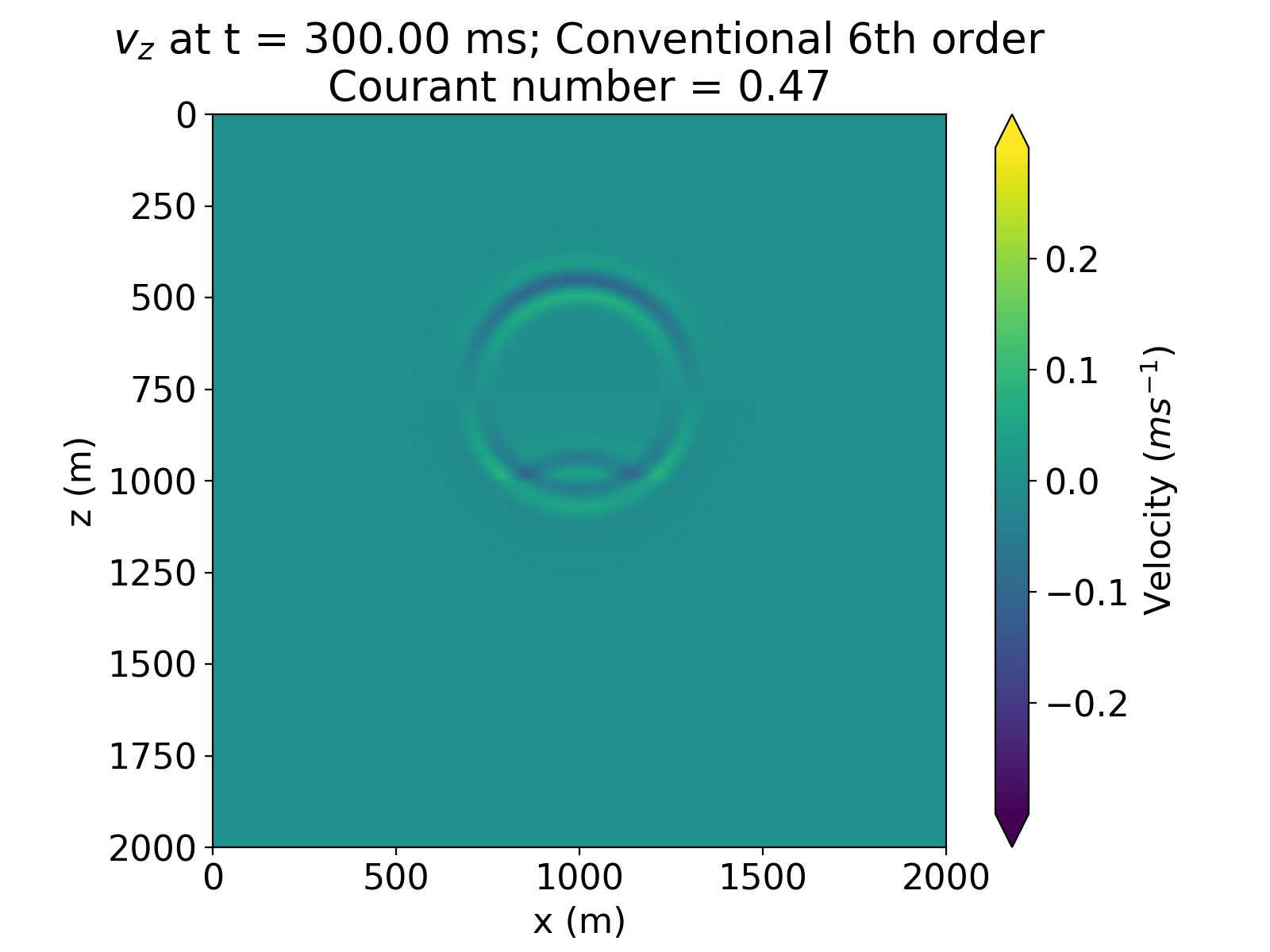}
				\caption{}
			\end{subfigure}\\
			\begin{subfigure}{0.45\textwidth}
				\centering\includegraphics[width=\textwidth]{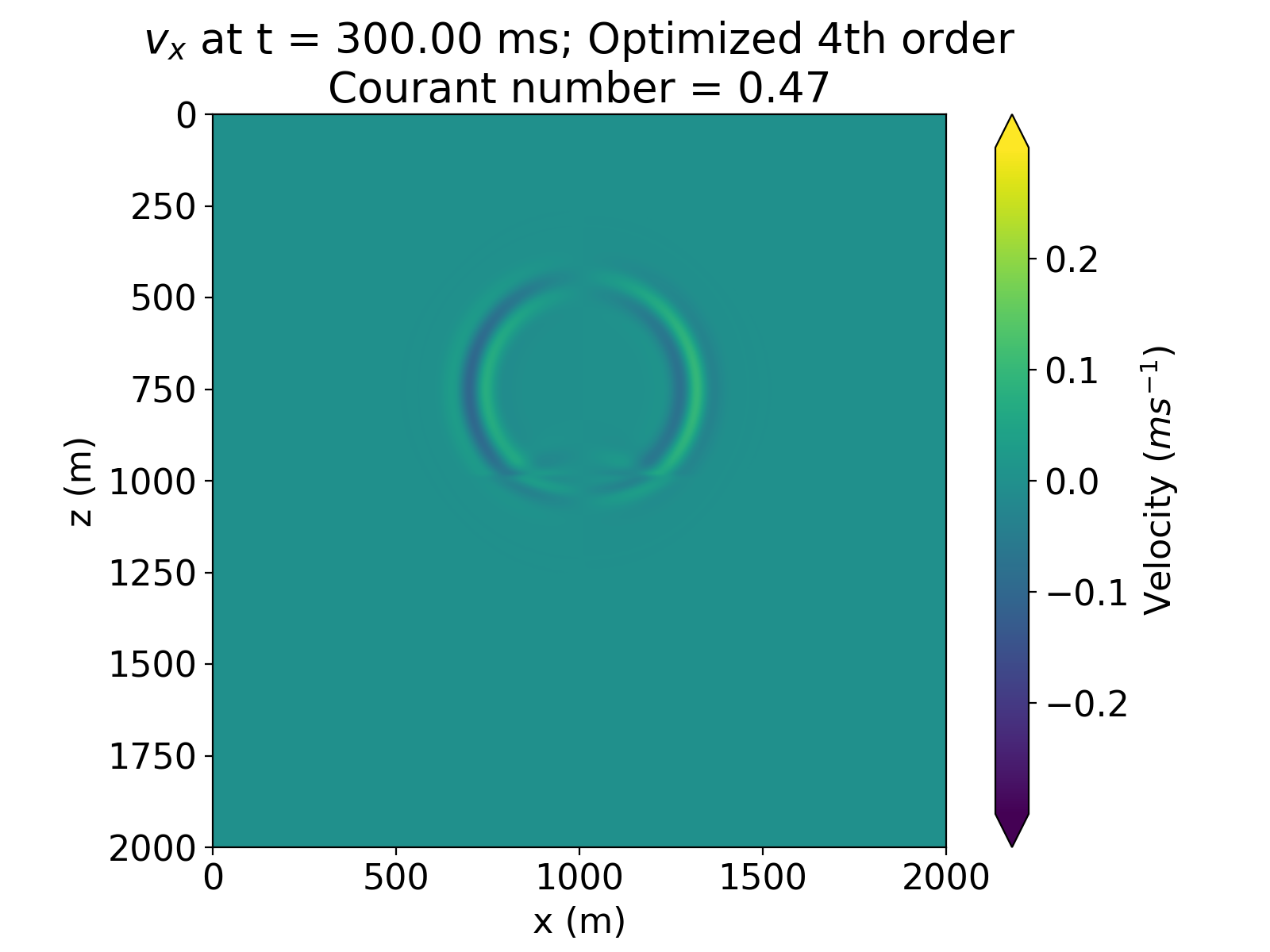}
				\caption{}
			\end{subfigure}
			\begin{subfigure}{0.45\textwidth}
				\centering\includegraphics[width=\textwidth]{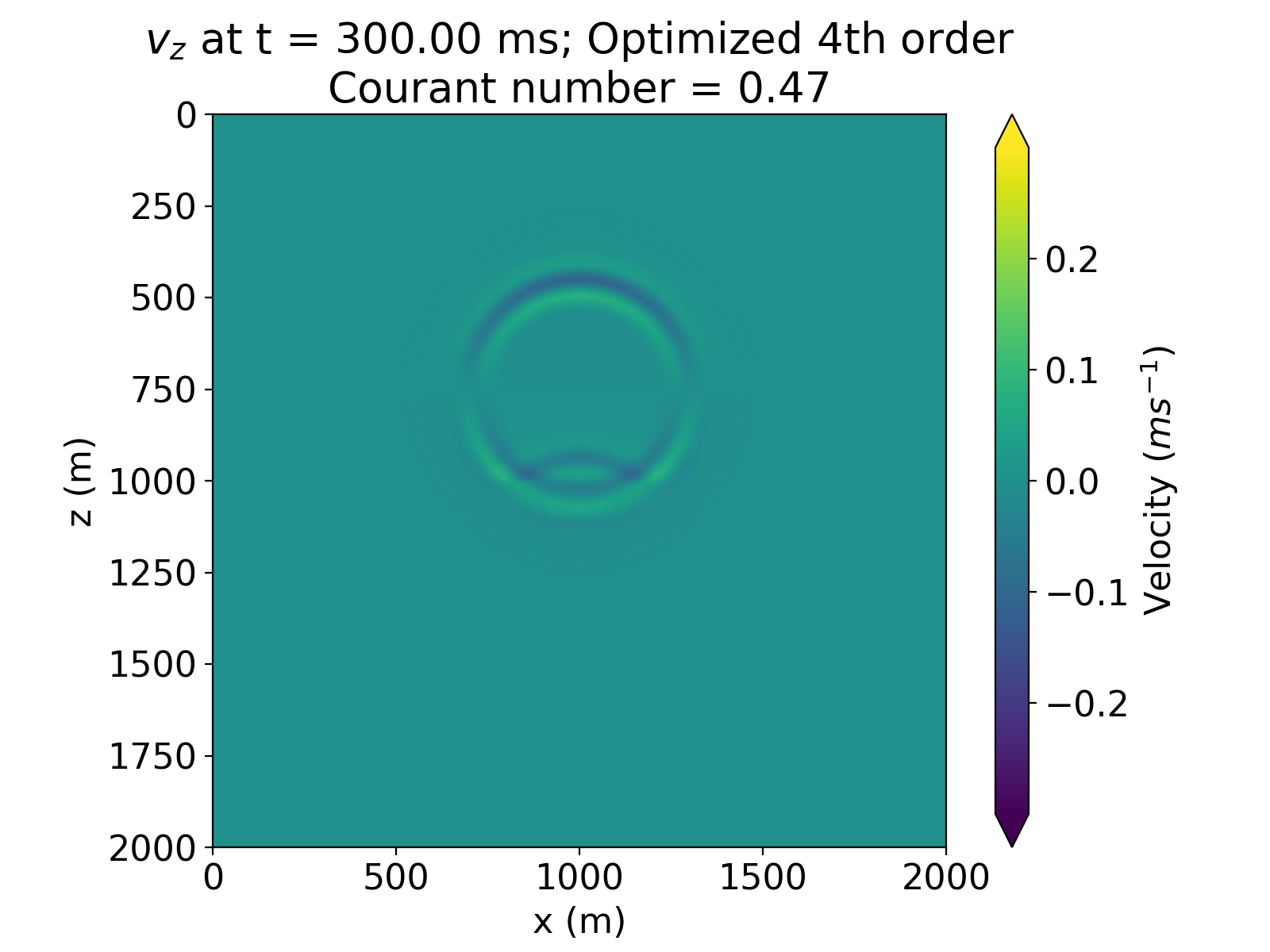}
				\caption{}
			\end{subfigure}
			\begin{subfigure}{0.45\textwidth}
				\centering\includegraphics[width=\textwidth]{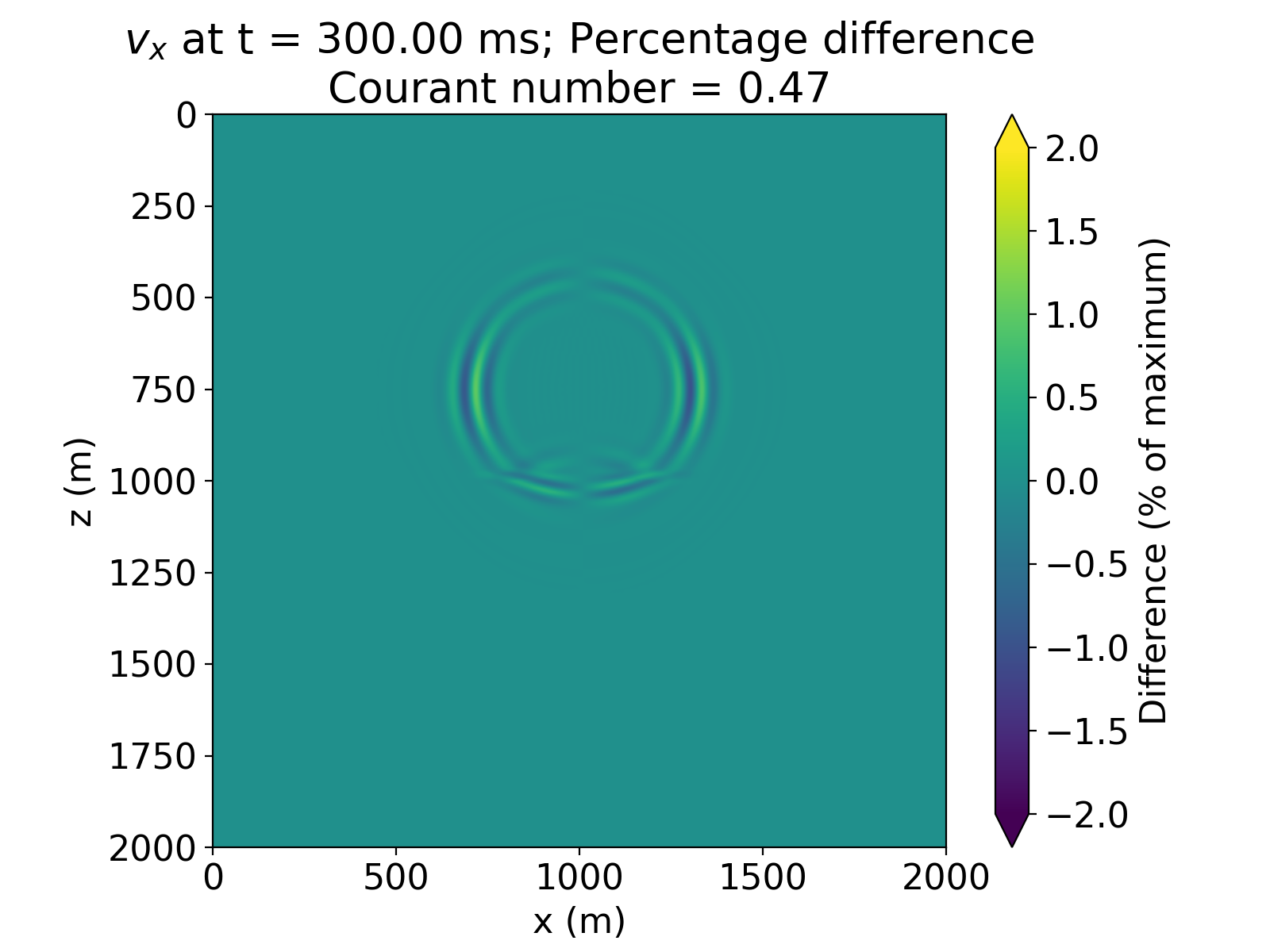}
				\caption{}
			\end{subfigure}
			\begin{subfigure}{0.45\textwidth}
				\centering\includegraphics[width=\textwidth]{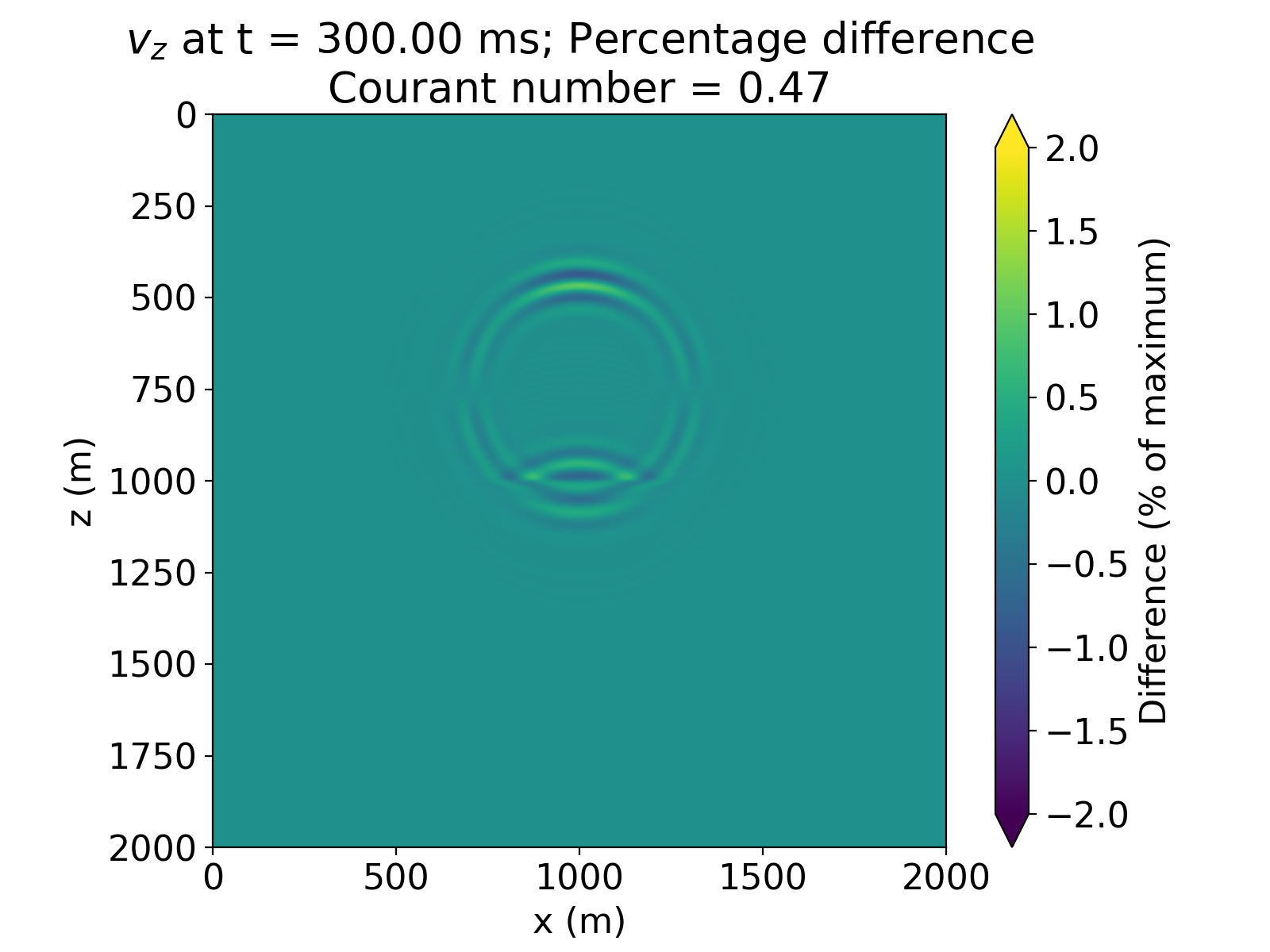}
				\caption{}
			\end{subfigure}
			\caption{Wavefields for velocity components $v_{x}$ and $v_{z}$ calculated using conventional and optimized schemes at 300ms. Difference is normalized against the largest amplitude present in the wavefield.}
			\label{2D_elastic_300}
		\end{center}
	\end{figure*}
	\begin{figure*}[h!]
		\begin{center}
			\begin{subfigure}{0.45\textwidth}
				\centering\includegraphics[width=\textwidth]{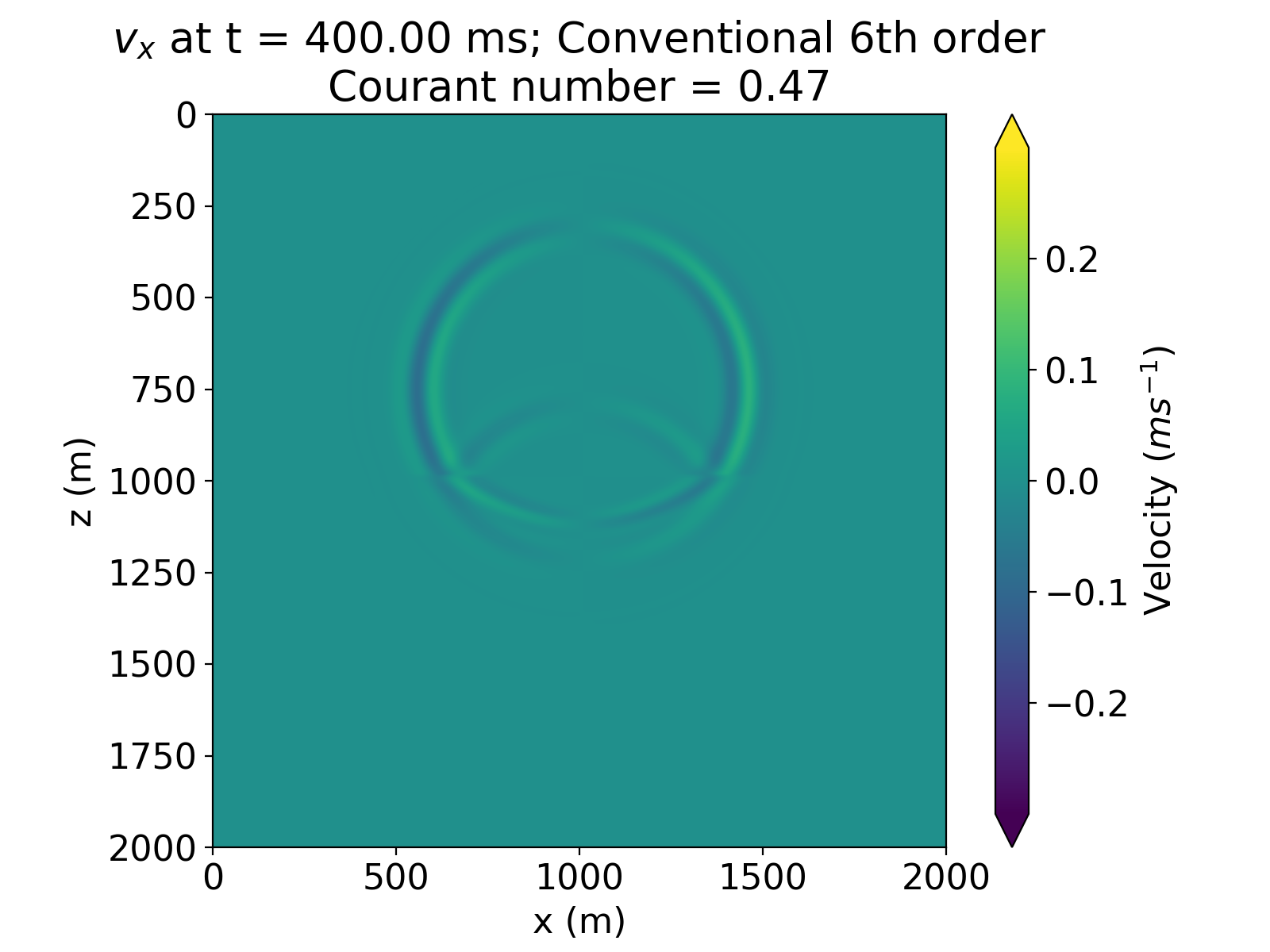}
				\caption{}
			\end{subfigure}
			\begin{subfigure}{0.45\textwidth}
				\centering\includegraphics[width=\textwidth]{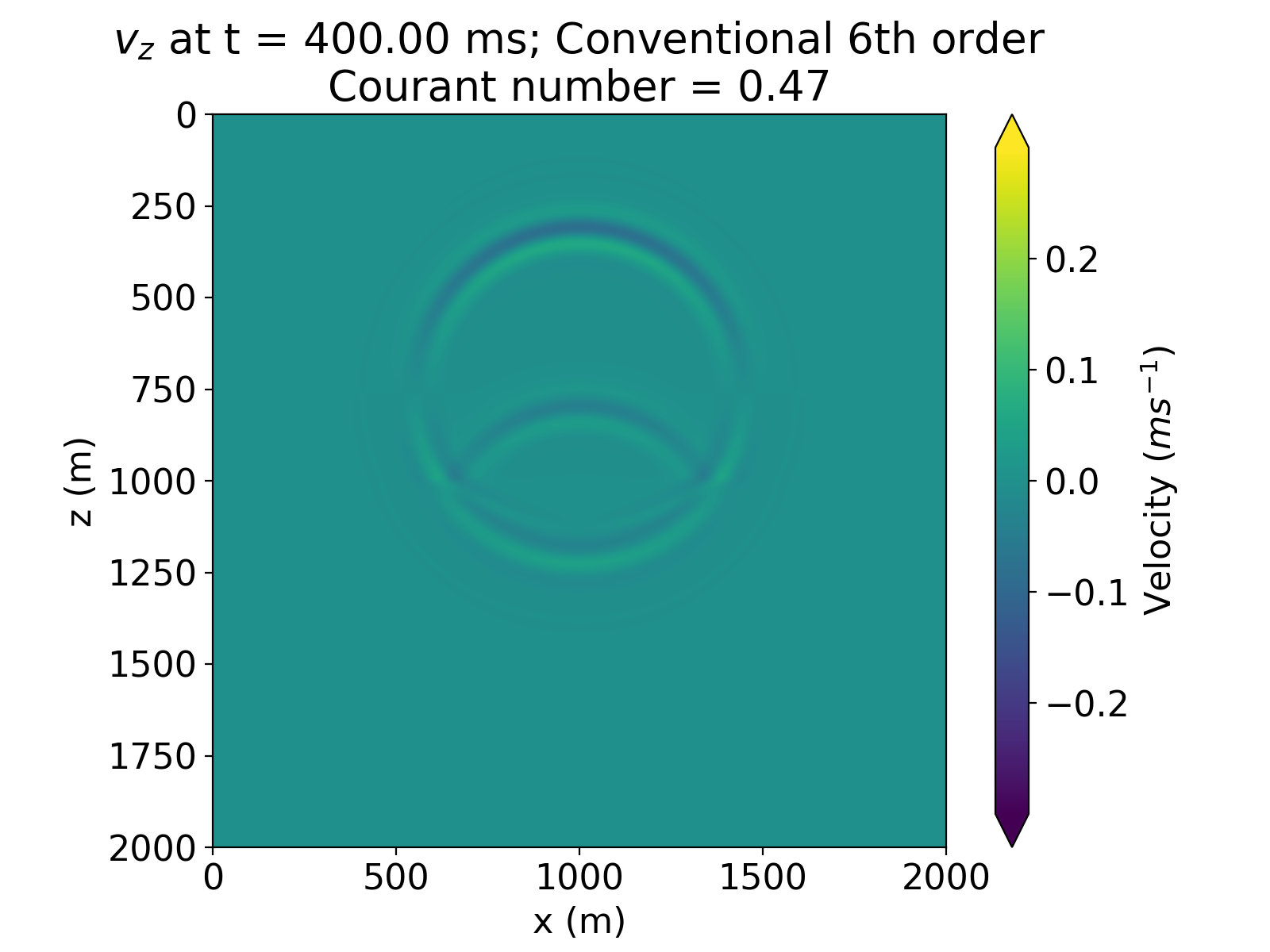}
				\caption{}
			\end{subfigure}\\
			\begin{subfigure}{0.45\textwidth}
				\centering\includegraphics[width=\textwidth]{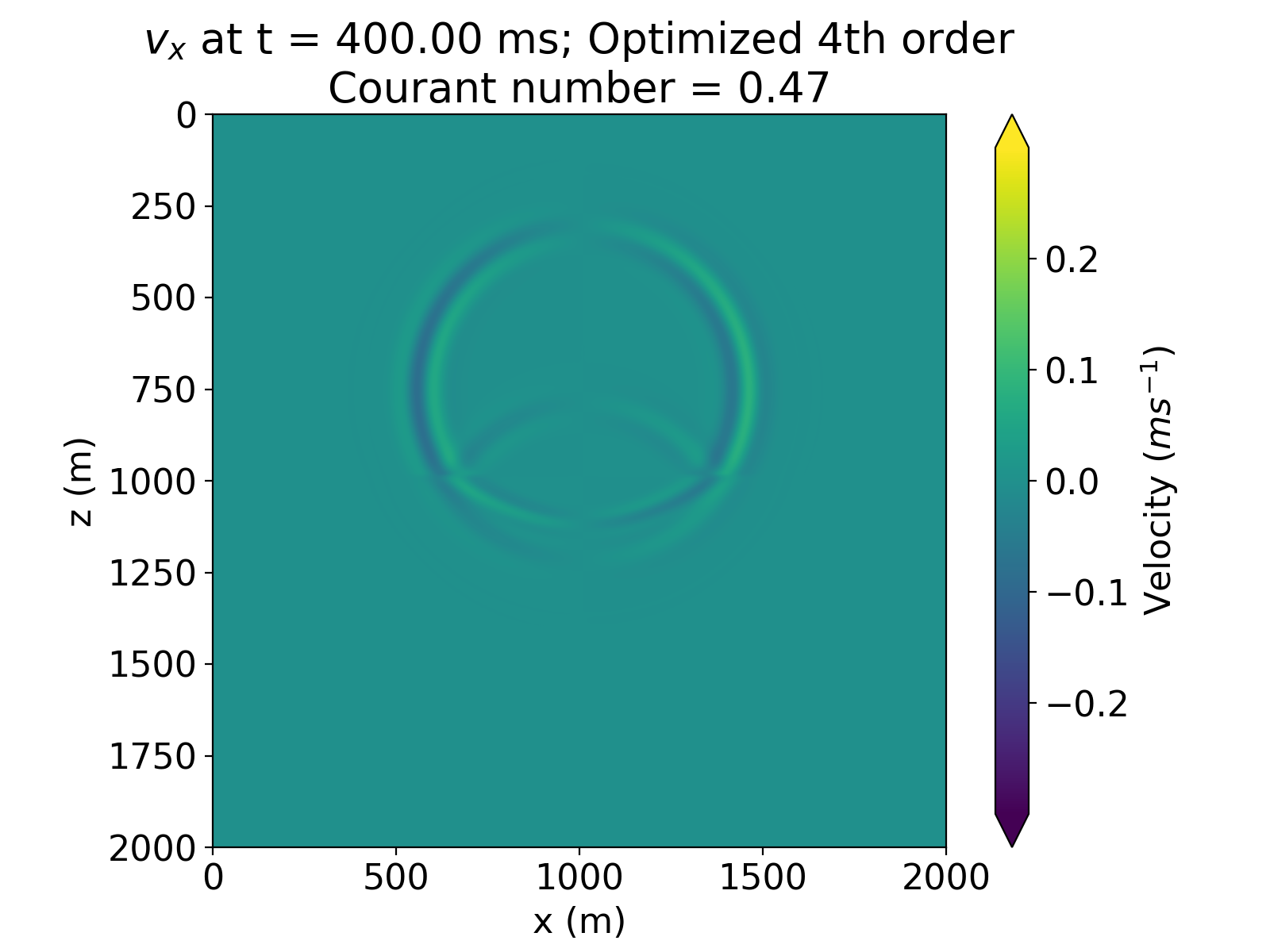}
				\caption{}
			\end{subfigure}
			\begin{subfigure}{0.45\textwidth}
				\centering\includegraphics[width=\textwidth]{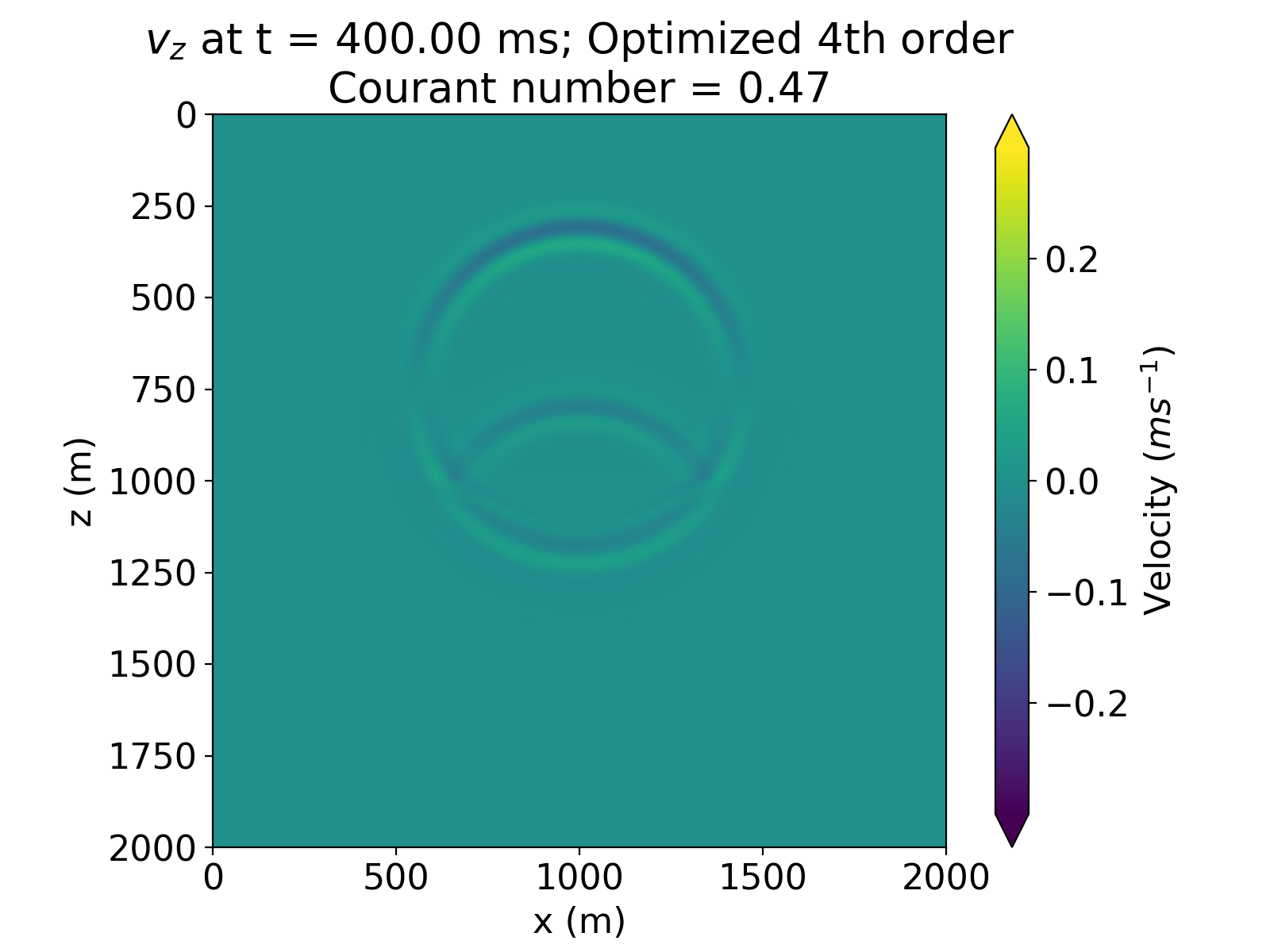}
				\caption{}
			\end{subfigure}
			\begin{subfigure}{0.45\textwidth}
				\centering\includegraphics[width=\textwidth]{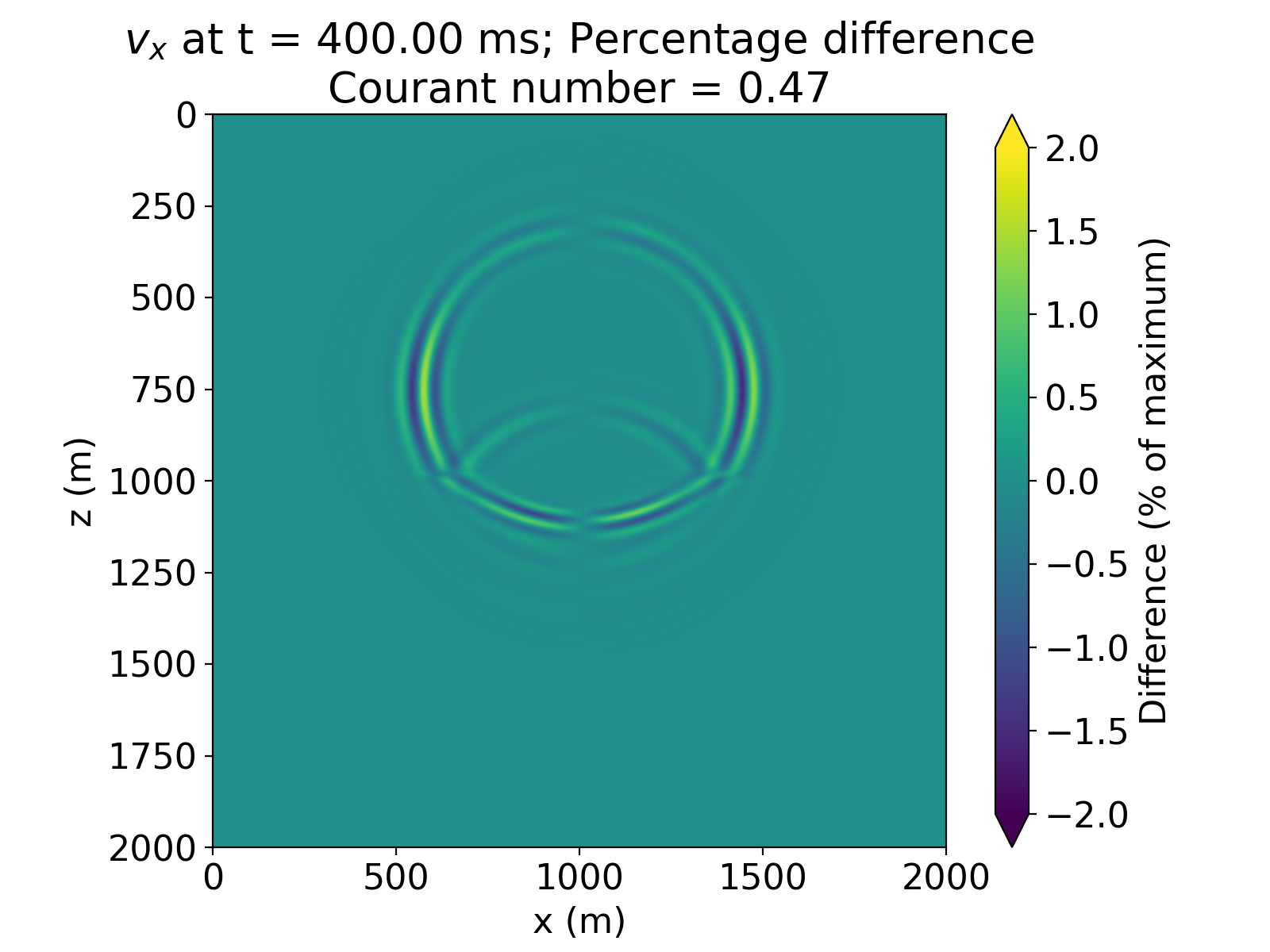}
				\caption{}
			\end{subfigure}
			\begin{subfigure}{0.45\textwidth}
				\centering\includegraphics[width=\textwidth]{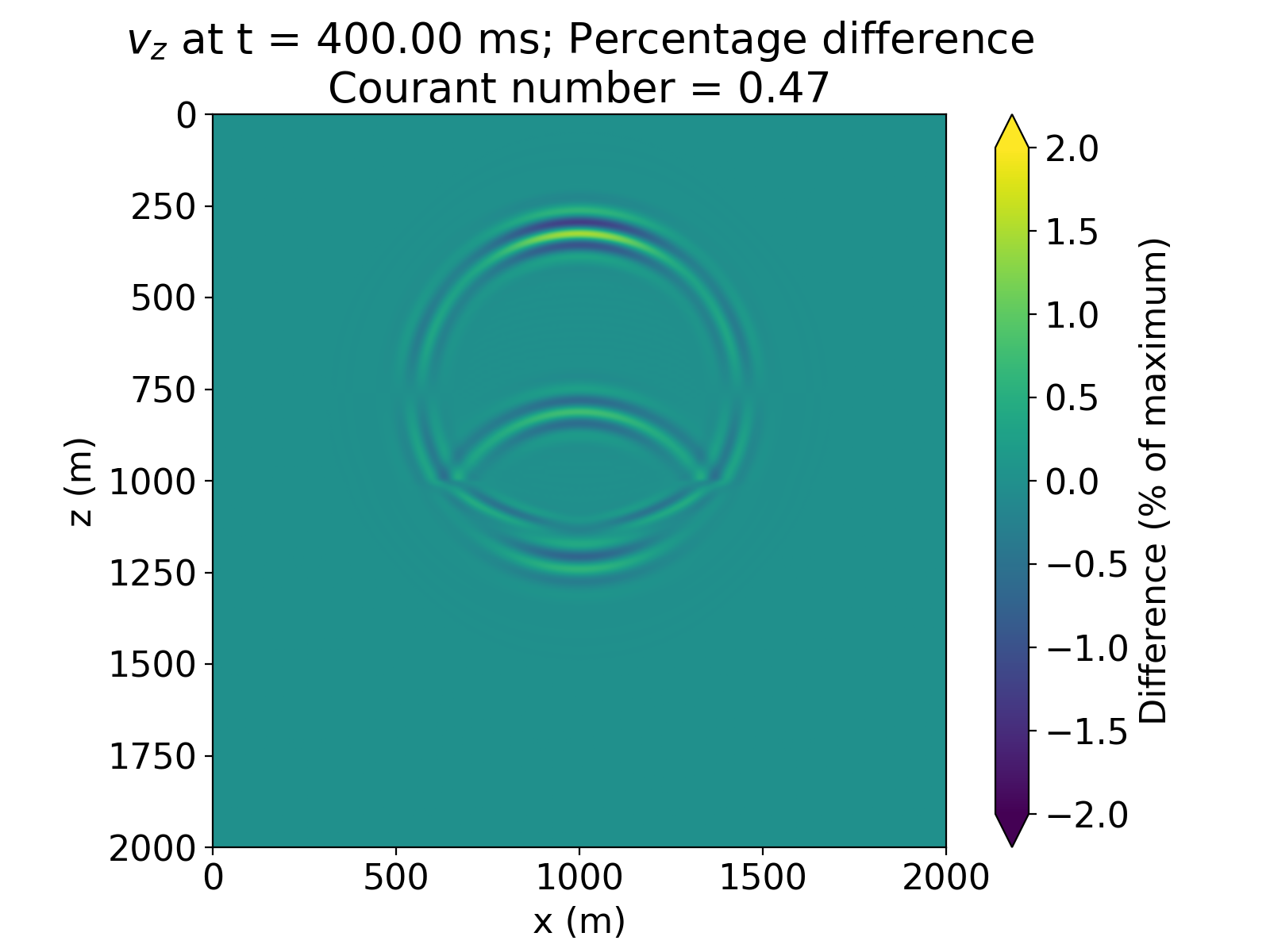}
				\caption{}
			\end{subfigure}
			\caption{Wavefields for velocity components $v_{x}$ and $v_{z}$ calculated using conventional and optimized schemes at 400ms. Difference is normalized against the largest amplitude present in the wavefield.}
			\label{2D_elastic_400}
		\end{center}
	\end{figure*}
	\begin{figure*}[h!]
		\begin{center}
			\begin{subfigure}{0.45\textwidth}
				\centering\includegraphics[width=\textwidth]{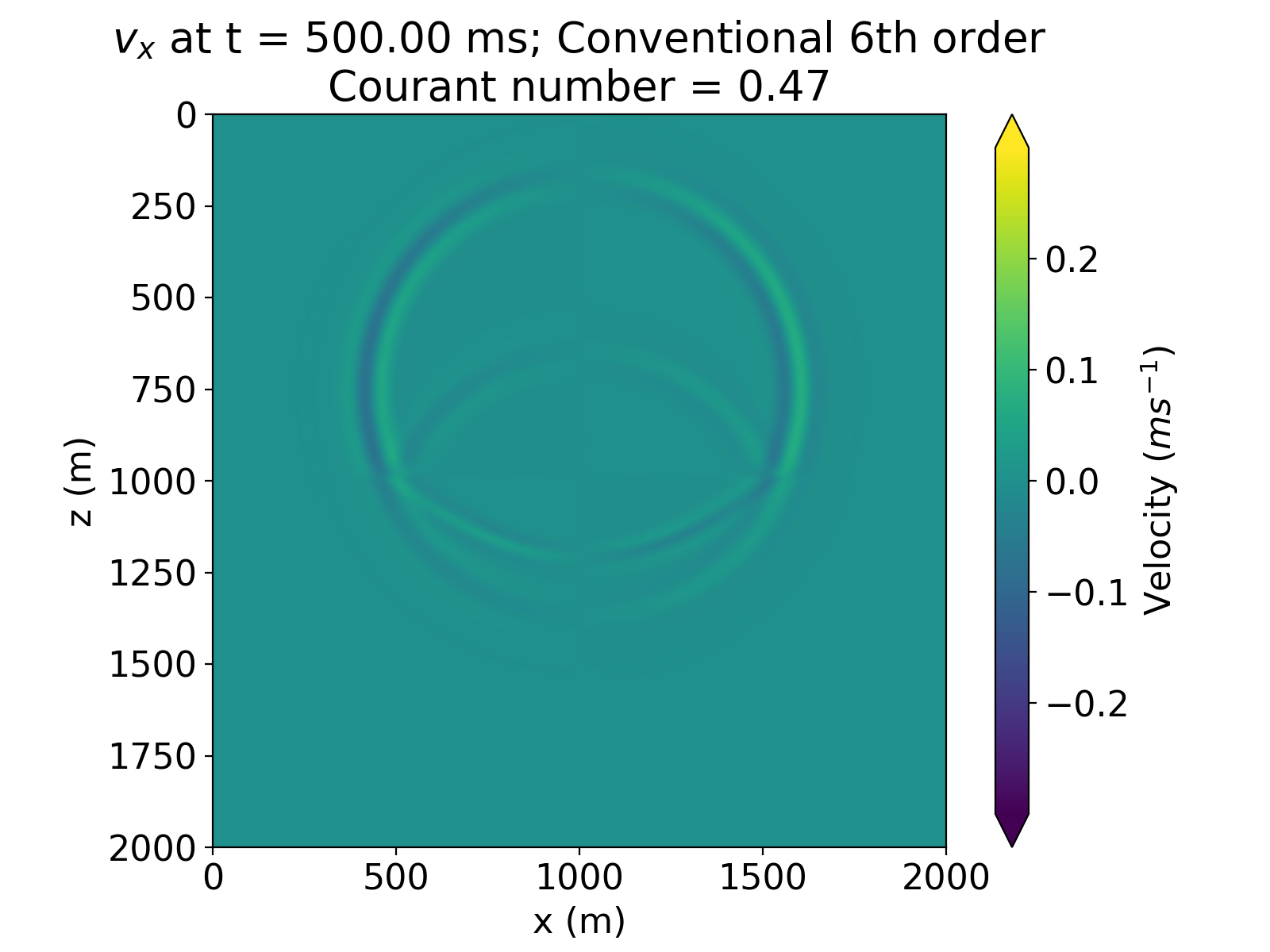}
				\caption{}
			\end{subfigure}
			\begin{subfigure}{0.45\textwidth}
				\centering\includegraphics[width=\textwidth]{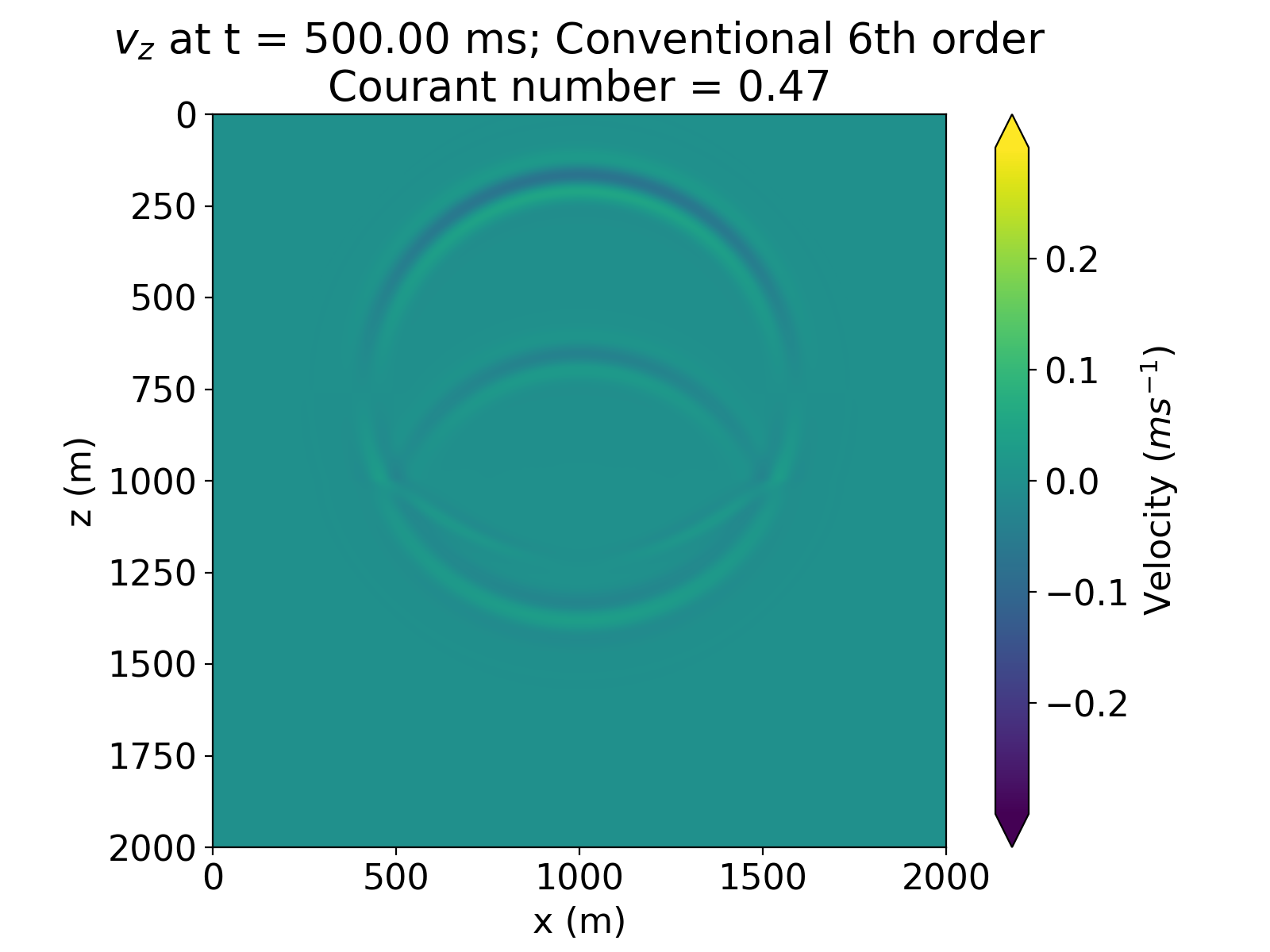}
				\caption{}
			\end{subfigure}\\
			\begin{subfigure}{0.45\textwidth}
				\centering\includegraphics[width=\textwidth]{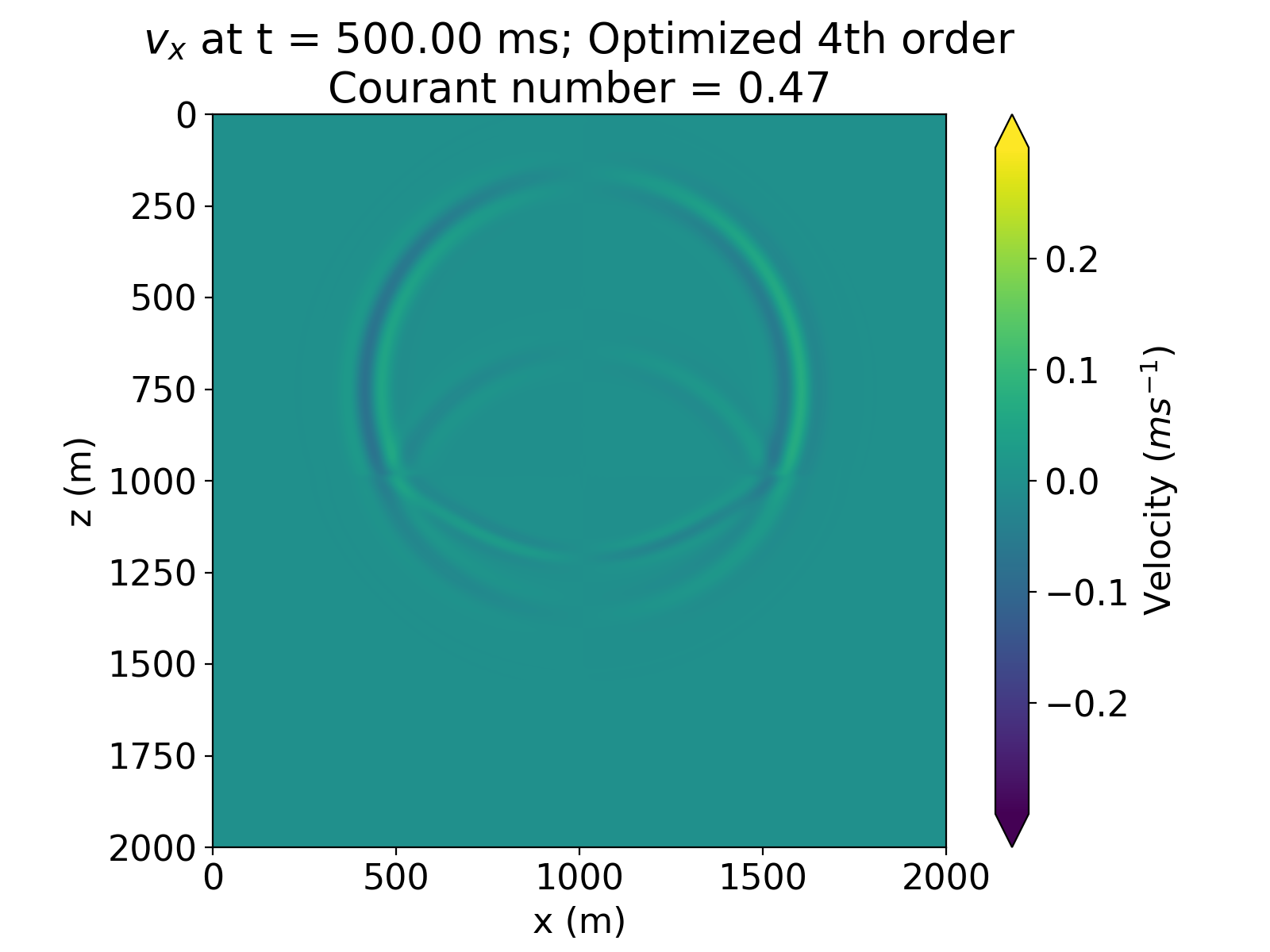}
				\caption{}
			\end{subfigure}
			\begin{subfigure}{0.45\textwidth}
				\centering\includegraphics[width=\textwidth]{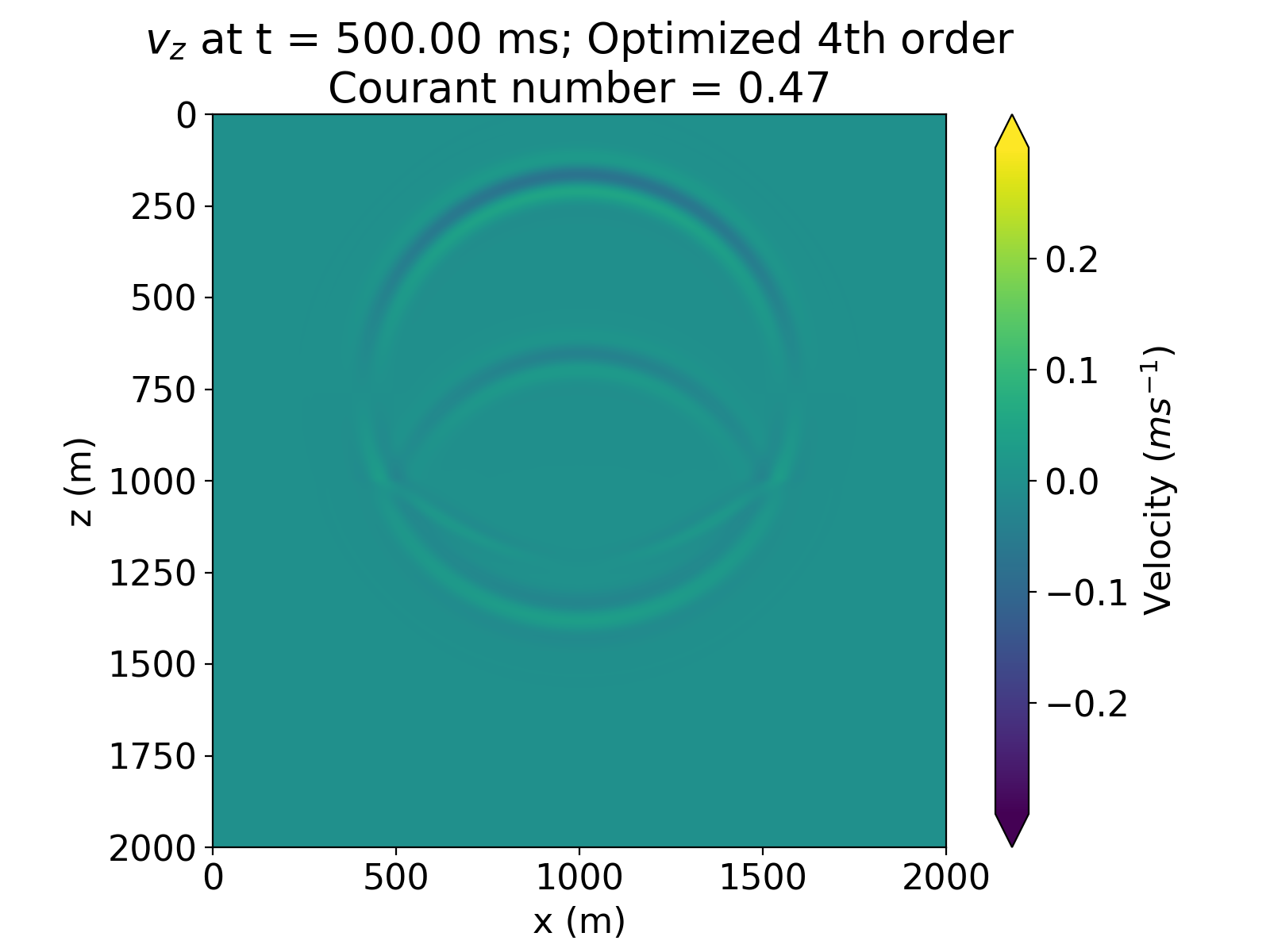}
				\caption{}
			\end{subfigure}
			\begin{subfigure}{0.45\textwidth}
				\centering\includegraphics[width=\textwidth]{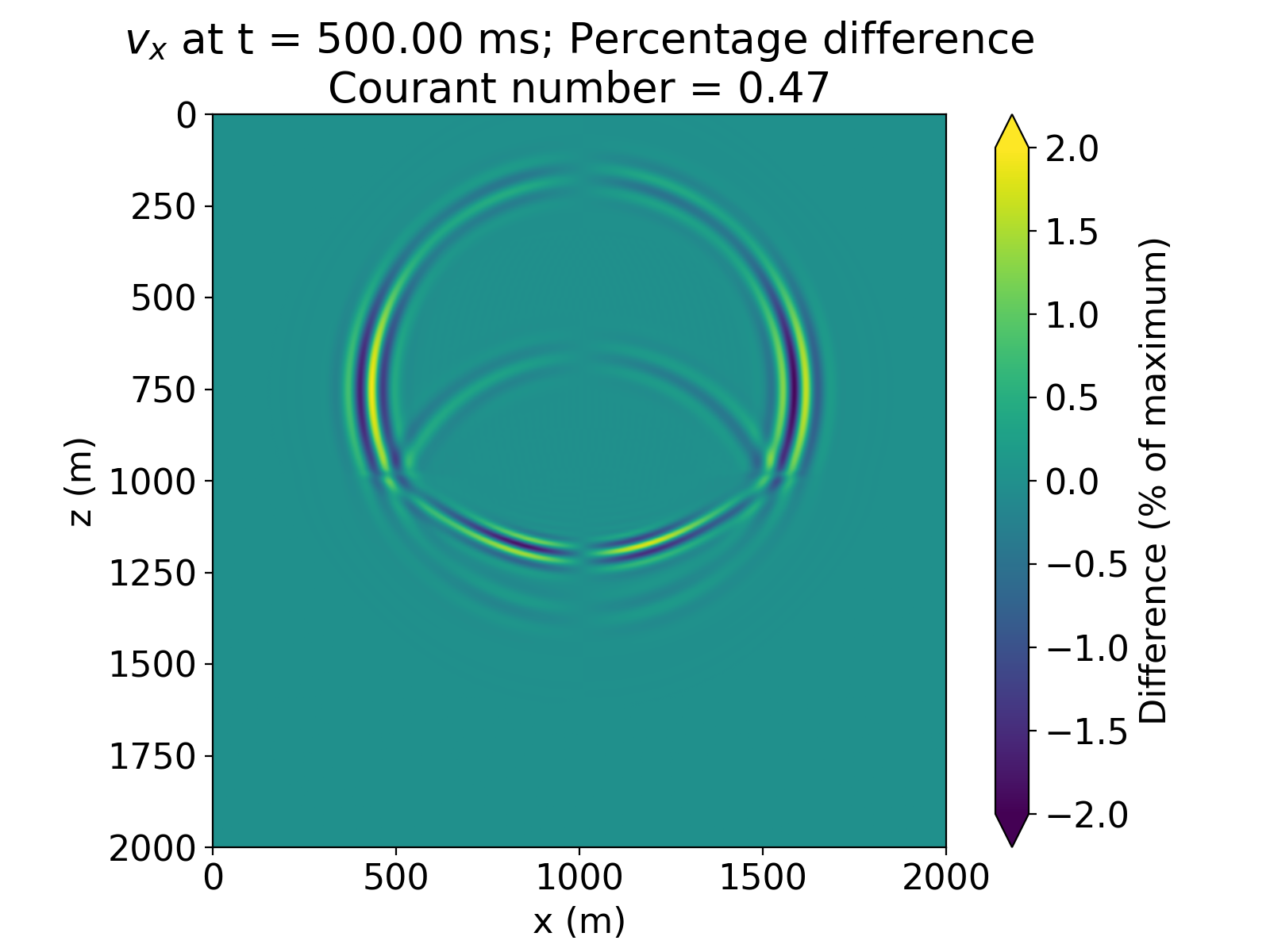}
				\caption{}
			\end{subfigure}
			\begin{subfigure}{0.45\textwidth}
				\centering\includegraphics[width=\textwidth]{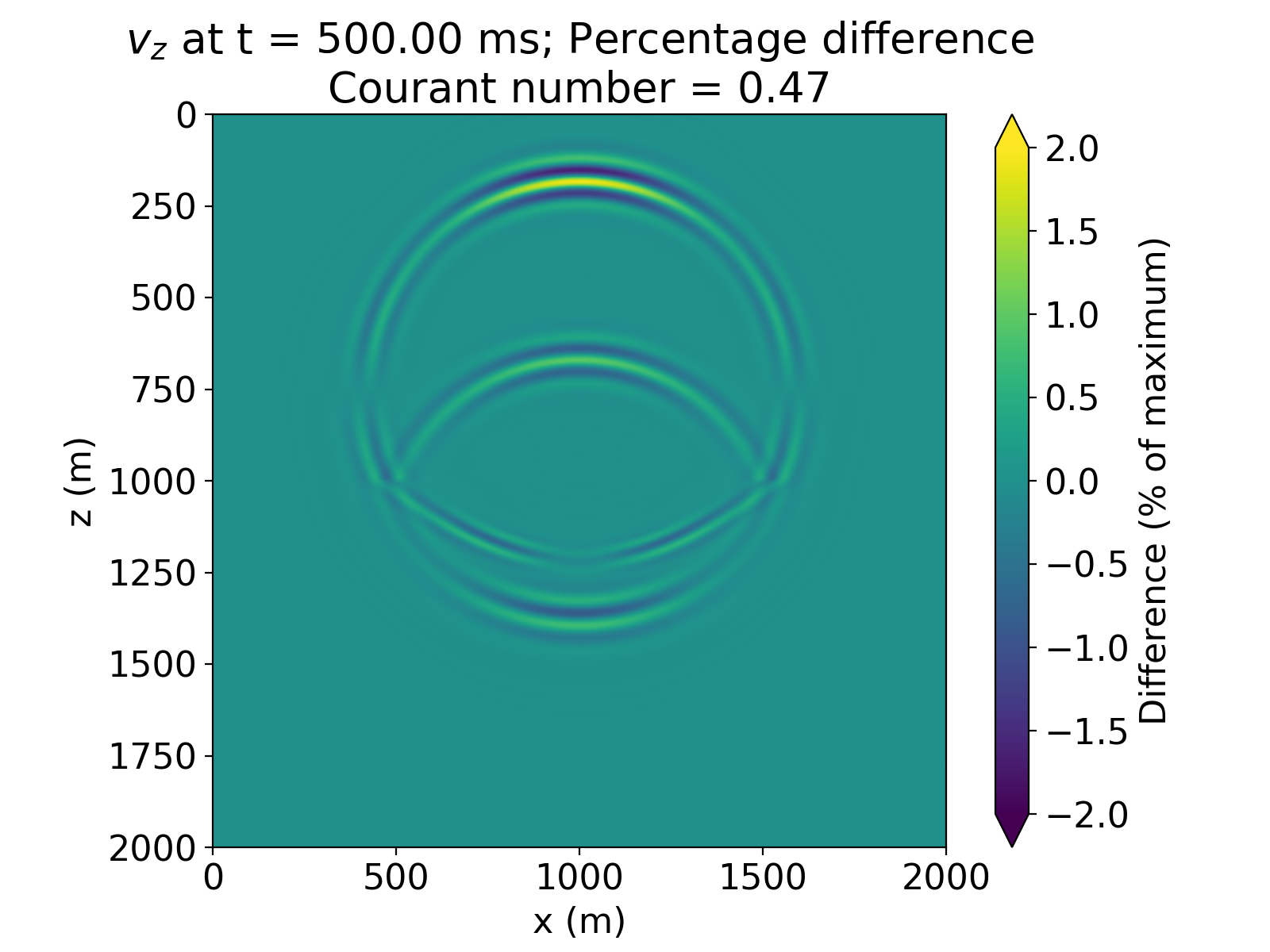}
				\caption{}
			\end{subfigure}
			\caption{Wavefields for velocity components $v_{x}$ and $v_{z}$ calculated using conventional and optimized schemes at 500ms. Difference is normalized against the largest amplitude present in the wavefield.}
			\label{2D_elastic_500}
		\end{center}
	\end{figure*}
	\begin{figure*}[h!]
		\begin{center}
			\begin{subfigure}{0.45\textwidth}
				\centering\includegraphics[width=\textwidth]{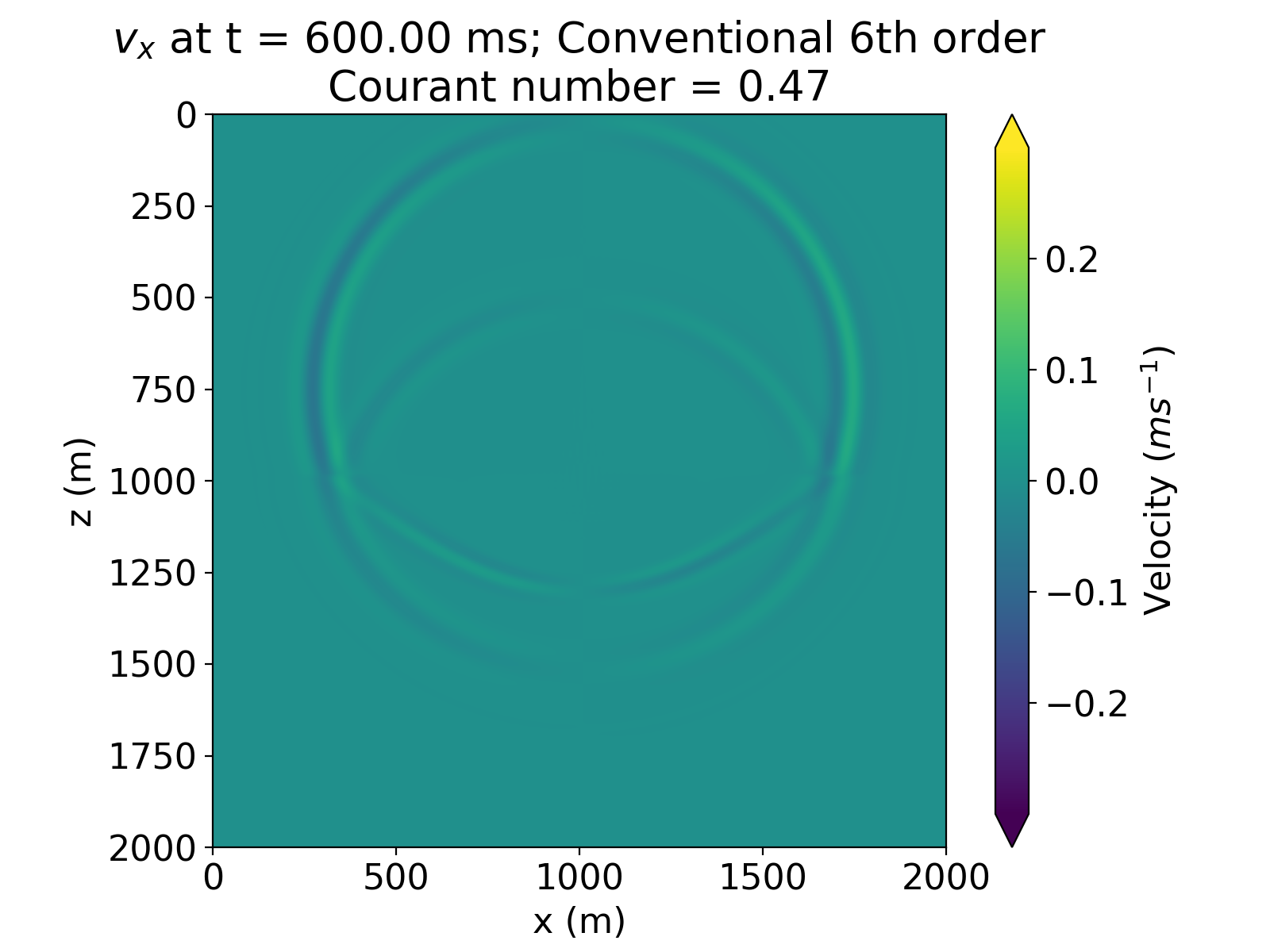}
				\caption{}
			\end{subfigure}
			\begin{subfigure}{0.45\textwidth}
				\centering\includegraphics[width=\textwidth]{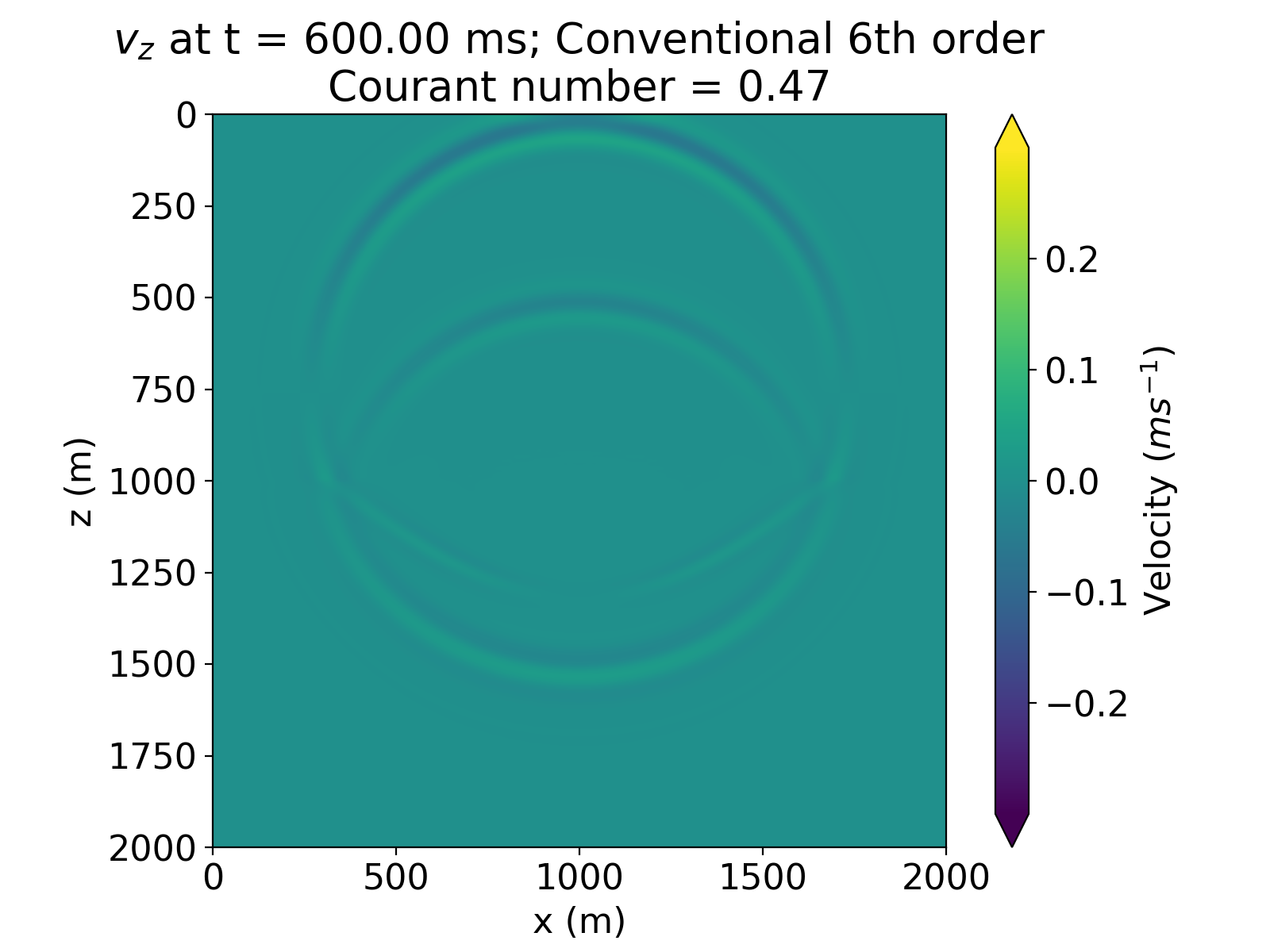}
				\caption{}
			\end{subfigure}\\
			\begin{subfigure}{0.45\textwidth}
				\centering\includegraphics[width=\textwidth]{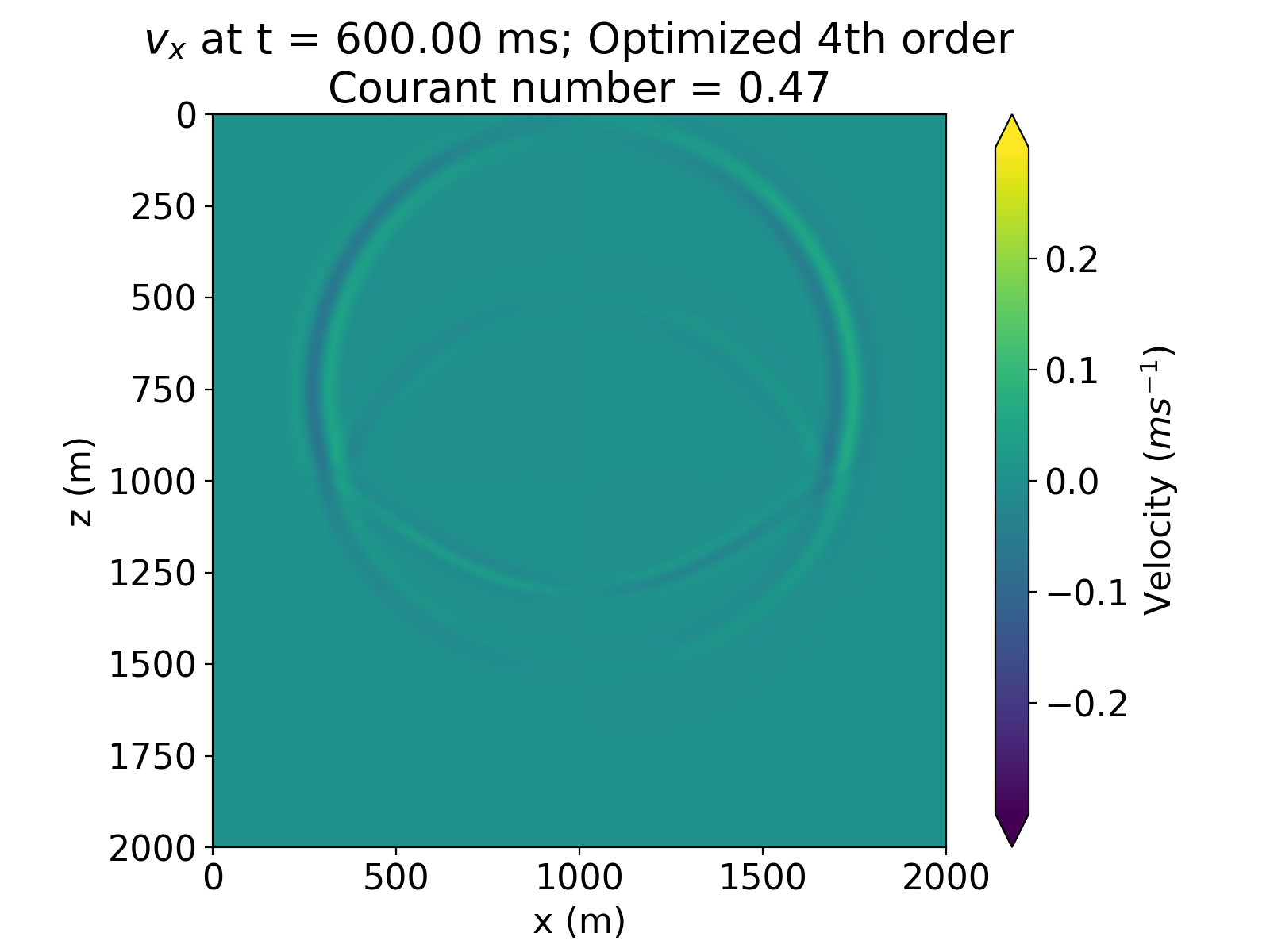}
				\caption{}
			\end{subfigure}
			\begin{subfigure}{0.45\textwidth}
				\centering\includegraphics[width=\textwidth]{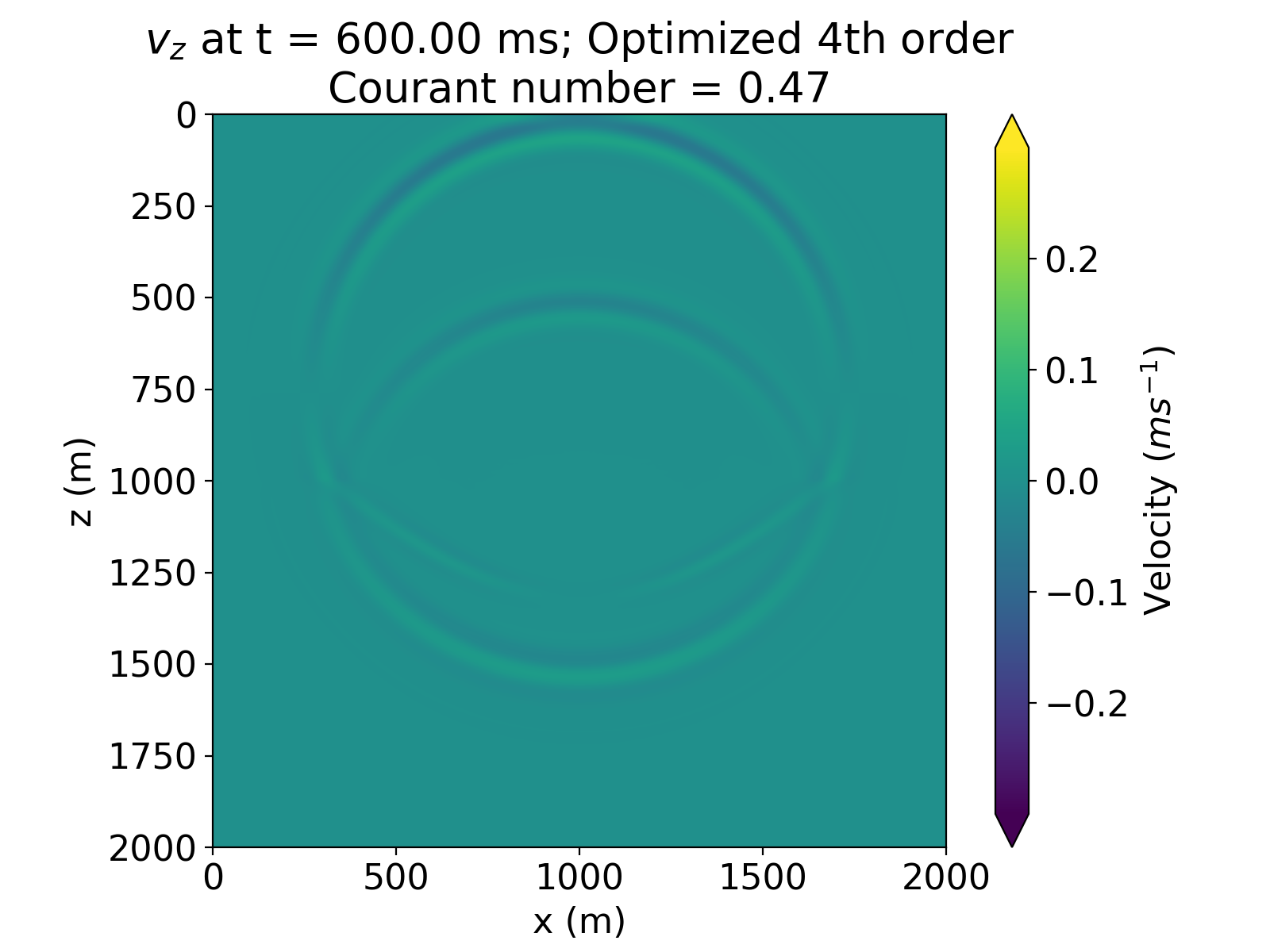}
				\caption{}
			\end{subfigure}
			\begin{subfigure}{0.45\textwidth}
				\centering\includegraphics[width=\textwidth]{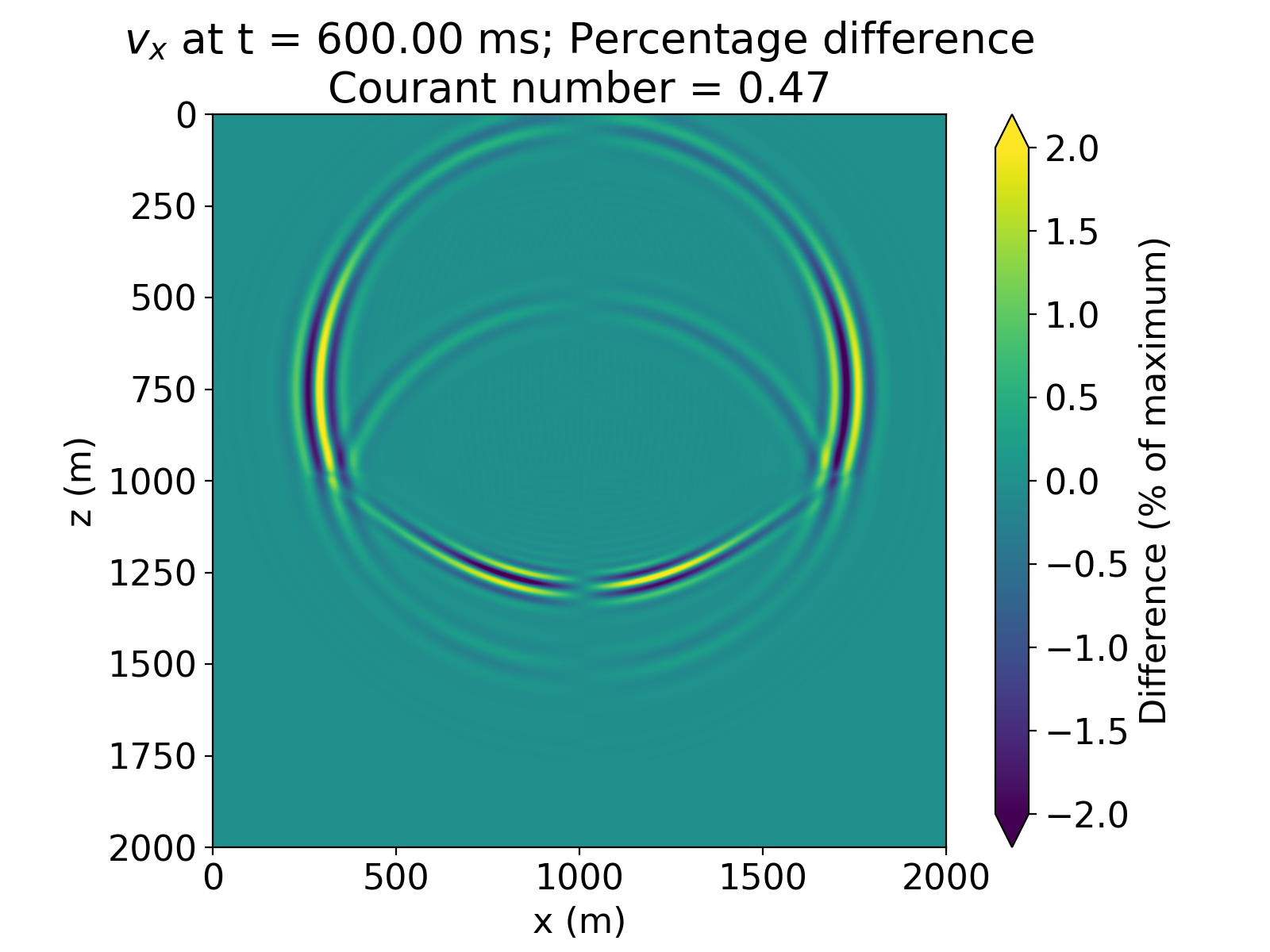}
				\caption{}
			\end{subfigure}
			\begin{subfigure}{0.45\textwidth}
				\centering\includegraphics[width=\textwidth]{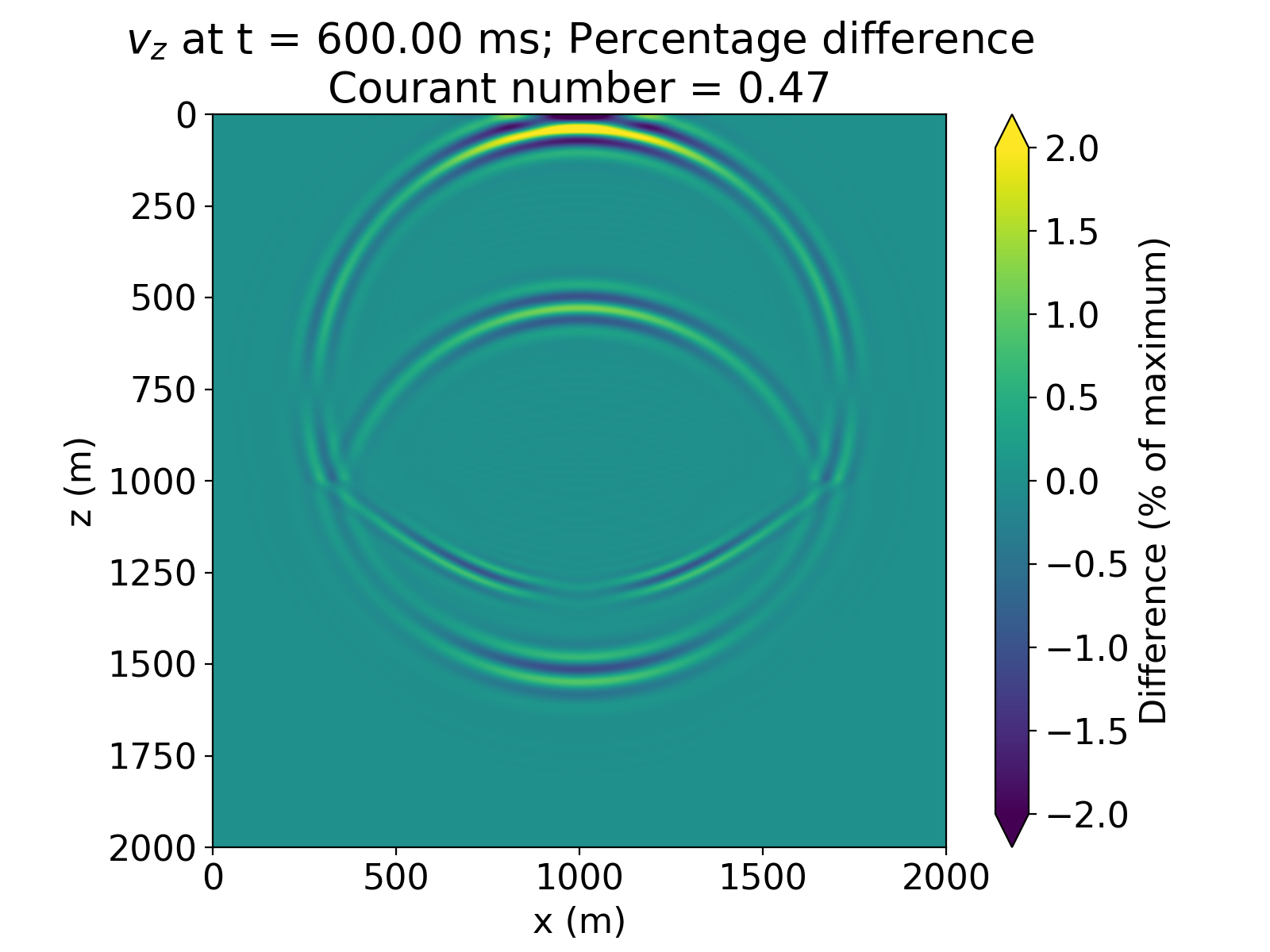}
				\caption{}
			\end{subfigure}
			\caption{Wavefields for velocity components $v_{x}$ and $v_{z}$ calculated using conventional and optimized schemes at 600ms. Difference is normalized against the largest amplitude present in the wavefield.}
			\label{2D_elastic_600}
		\end{center}
	\end{figure*}
	\begin{figure*}[h!]
		\begin{center}
			\begin{subfigure}{0.45\textwidth}
				\centering\includegraphics[width=\textwidth]{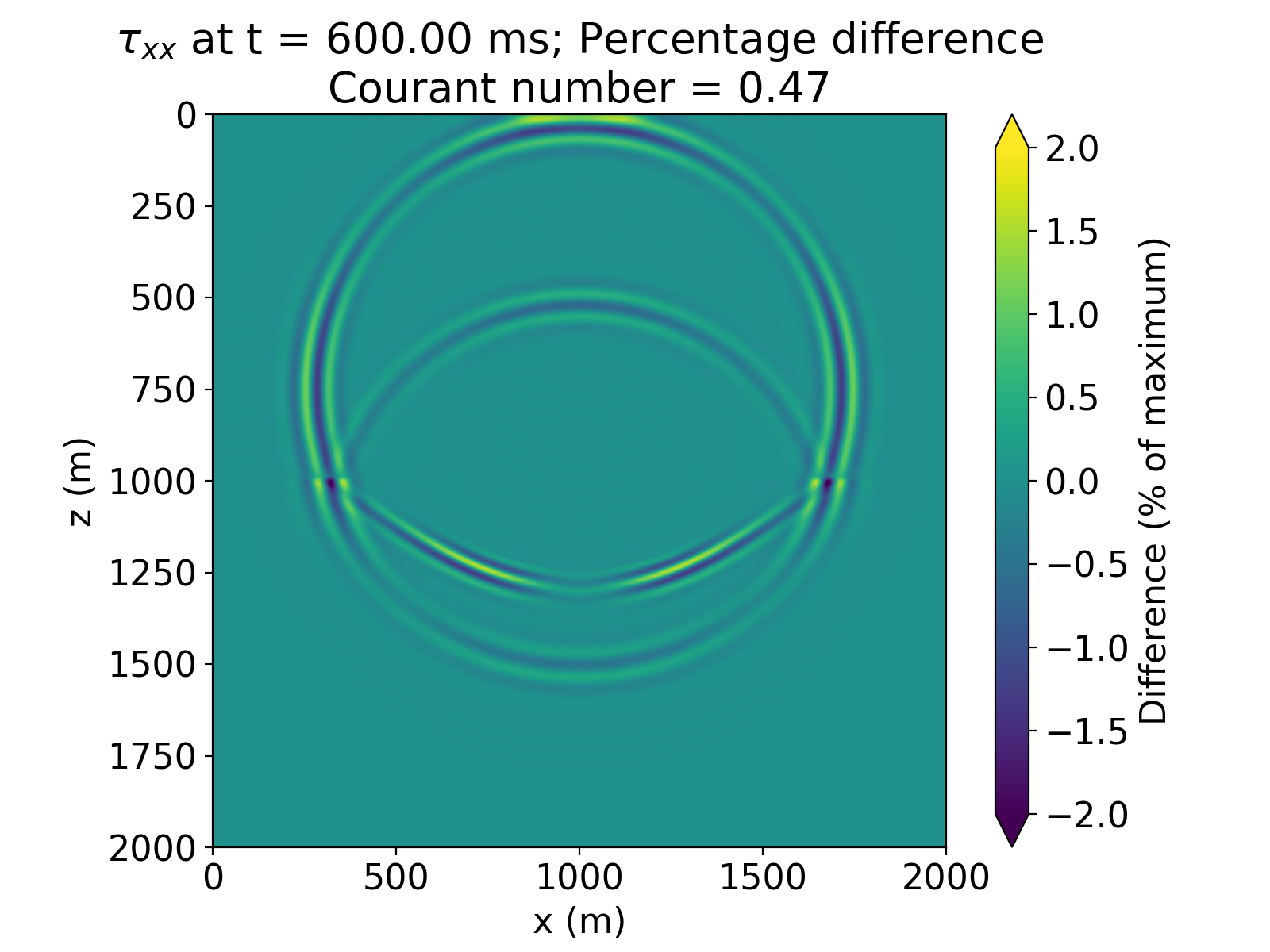}
				\caption{}
			\end{subfigure}
			\begin{subfigure}{0.45\textwidth}
				\centering\includegraphics[width=\textwidth]{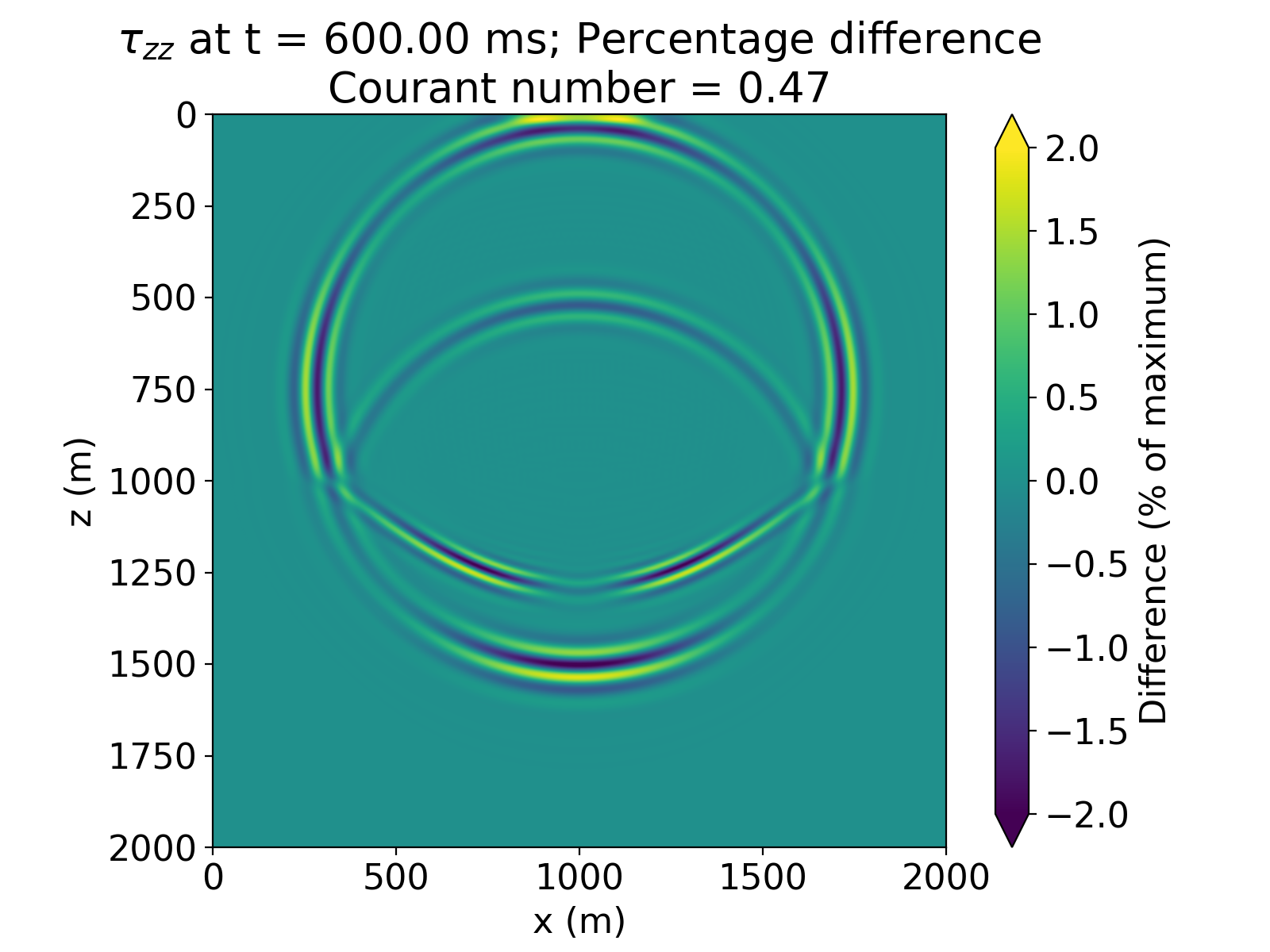}
				\caption{}
			\end{subfigure}\\
			\begin{subfigure}{0.9\textwidth}
				\centering\includegraphics[width=\textwidth]{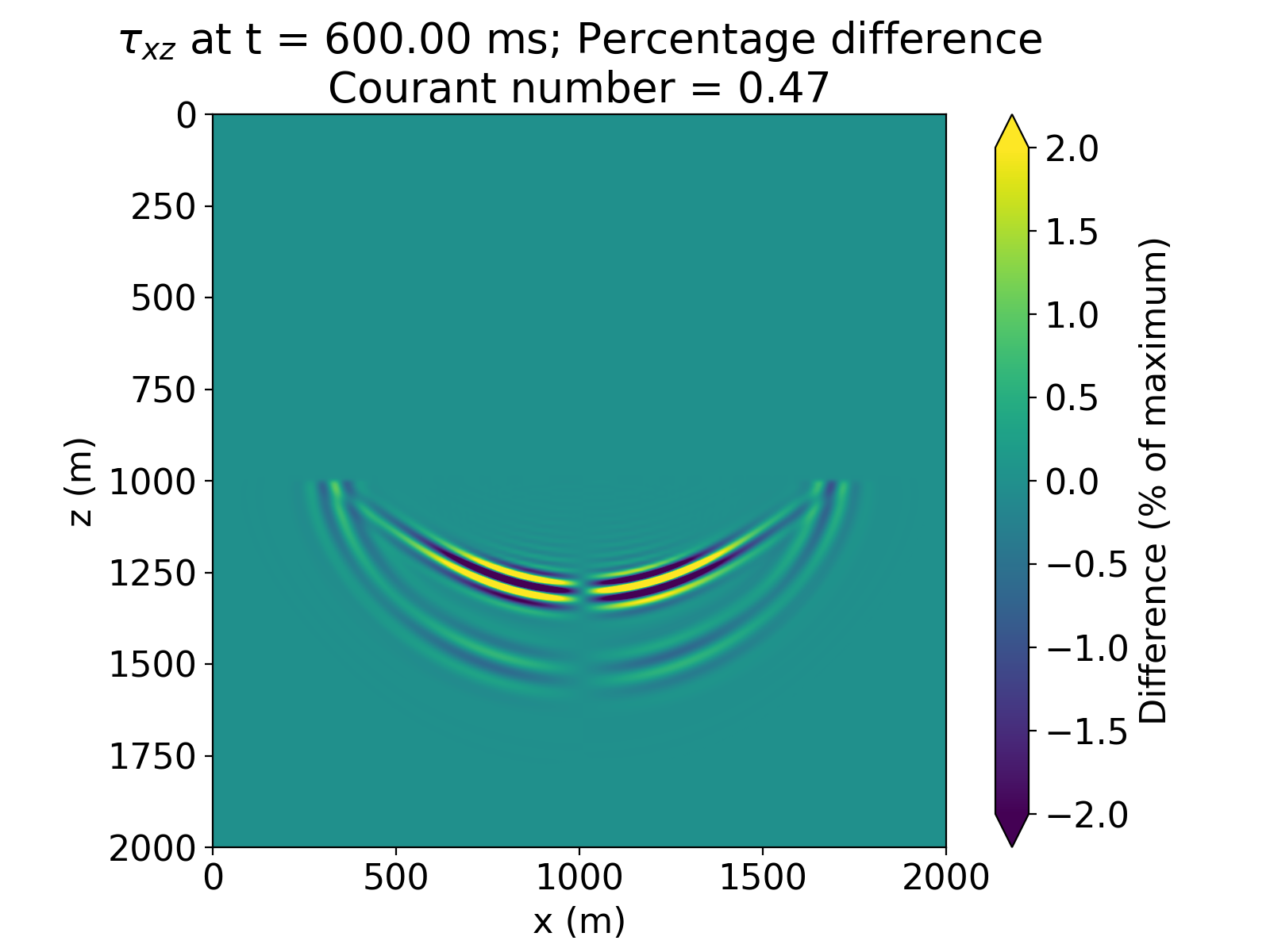}
				\caption{}
			\end{subfigure}
			\caption{Difference between stress component wavefields for conventional and optimized schemes at 600ms (post-incidence on the shale layer), normalized against the largest amplitude present in the wavefield. S-waves are typically subject to increased numerical dispersion, as lower propagation speeds correspond with shorter wavelengths, thus coarsening wave sampling (\citealp{Moczo2000}). Figure \ref{2D_elastic_stress} shows the highest-amplitude artifacts removed by the optimized scheme present in the $\uptau_{xz}$ field, consistent with expectation. Artifact extent is also greatest in the $\uptau_{xz}$ field, with trains of high-frequency ripples trailing behind the shear wavefront. Increased dispersion in the shear stress field is liable to induce inaccuracies in other fields if left to propagate, due to the coupled nature of stress and velocity for P-SV formulations.}
			\label{2D_elastic_stress}
		\end{center}
	\end{figure*}
	\clearpage
	Particle velocity components $v_{x}$ and $v_{z}$ over time are shown in figures \ref{2D_elastic_100} through \ref{2D_elastic_600}. Little difference is discernible between the schemes at any time: not unexpected as both have similar stencil coefficients, producing near-identical dispersion characteristics. However, normalized difference between the wavefields for each scheme reveals a ring shaped artifact with wavelength half that of the main wavefront, and amplitude between 1\% and 2\% of the maximum amplitude. This artifact appears in the same position as the main body of the wave at all times. Upon inspection, this artifact appears squared, diagnostic of numerical dispersion (\citealp{Juntunen2000}), and extends inwards from the wavefront, characteristic of lagging high frequencies. Thus difference between the conventional and optimized FD solutions presumably consists primarily of numerical dispersion artifacts. The magnitude of the difference between fields generated by the conventional and numerical schemes is generally in agreement with differences in the expected degree of numerical dispersion suppression provided by each.
	\\
	\\Both stencils generate satisfactory results for this model, thereby validating use of staggered spatially-optimized stencils for elastic wave propagation.

	\section{Discussion}
	\subsection{Summary of test cases}
	For both test cases, the use of spatially-optimized 4\textsuperscript{th} order stencils in the place of conventional 6\textsuperscript{th} order yields notable improvements to simulation quality for $\Delta x$ near or below the critical value for dispersion relation preservation, including for heavily undersampled wavefields. However, reduced truncation error in the higher-order conventional scheme results in better approximation of the governing PDE for finer grids, as the highest frequency wavefield components are sufficiently sampled to retain their dispersion characteristics. Note that whilst the optimized scheme yielded the smaller error at very small $\Delta x$ for the first test case, this likely arose due to the parameters or nature of the test case in absence of a robust explanation.
	\\
	\\Upon return to inital or reversed wave positions for the acoustic test case, the optimized solution tended to contain the smaller error, with smaller mean absolute error at simulation end for all $\Delta x$ tested within the window $0.0200\mathrm{m}\leq\Delta x\leq 0.0250\mathrm{m}$. However, during transition between these states, misfit between the analytical and optimized FD solutions grows considerably. Whilst apparent in both schemes, the contrast appears more pronounced for the optimized scheme. This periodicity implies more accurate propagation of frequency components within the wave by the optimized scheme, minimizing distortion when the wave returns to its initial position. During the intermediate phase, the optimized scheme displays a tendency to generate high amplitude spikes adjacent to steep gradients which propagate with the position of these gradients. This should be noted when optimized schemes are used with models capable of producing steep gradients, as interpretation of these spurious spikes as meaningful data is liable to lead to inaccurate conclusions. This behavior could possibly be attributed to numerical coupling phenomena described by \citealp{Tam1993}, causing erroneous driving of frequency components when only the temporal or spatial stencils of a DRP scheme are used. In this case, use of the corresponding temporal derivative detailed in \citealp{Tam1993} should address these artifacts. 
	\\
	\\Application of spatially-optimized stencils to staggered grids demonstrates viability of non-conventional FD stencils for elastic modeling. The optimized FD solution was found to be of similar quality to the conventional, although the difference between the two revealed limited increase in numerical dispersion in the conventional scheme. Otherwise, evolution of the wave was realistic and values observed in all fields were sensible for both schemes. The squared appearance of the dispersion artifact, with faces perpendicular to the x and z axis is consistent with maximum dispersion parallel with the x and z axis (\citealp{Levander1988}). This spurious directional dependence is apparent in orientations of the highest amplitude regions of the dispersion artifact, at left and right-most extremities for the $v_{x}$ field, and top and bottom for $v_{z}$. The susceptibility of S-waves to numerical dispersion highlighted during testing underlines the benefits of optimized schemes for multi-phase wavefields and dispersion-free handling of high-spatial-frequency waves originating from conversions at layer boundaries with reduced dispersion compared to conventional schemes. This implies more accurate propagation through complex models containing numerous layer boundaries, supported by findings in \citealp{Zhang2013}. The noteworthy reduction in numerical dispersion achieved demonstrates the transferability of the scheme's DRP nature across multiple applications.
	\\
	\\Despite considerable benefits of optimized schemes when applied to undersampled wavefields, the magnitude of error present was significantly greater than in adequately sampled wavefields. In the absence of computational constraints, use of high-order conventional FD schemes seems the obvious choice for numerically solving PDEs; fine grid spacings minimize the comparative benefits of optimized schemes with equivalent stencil extent. However, as memory constraints are commonplace in modeling applications, it is typically beneficial or even necessary to minimize the number of nodes within the computational domain. A smaller grid also ensures that simulation runtimes are not excessive, reducing the number of calculations per timestep. If the Courant number is held constant, then larger grid spacings yield larger timesteps, further reducing required calculation. If expediency is a priority, use of optimized schemes is preferable, not due to any inherent reduction in the number of computational operations, but because coarser spatial and temporal discretizations are possible.
	\\
	\\For seismic applications, spatial optimization should typically be beneficial, despite only minimal improvement apparent in the second test case. Due to large memory requirements of the commonly used P-SV formulation (\citealp{Virieux1986}, \citealp{Bartolo2015}), the ability to run simulations at near critical $\Delta x$ with marginal penalty is a boon. Combined with variable grid and time increments, this is exploitable to achieve an approximately 25\% reduction in computational cost and memory demand (\citealp{Zhang2013}) without negatively impacting simulation quality. Note that as order increases, benefits of spatially optimizing a staggered scheme become increasingly pronounced (\citealp{Liu2014}). With the common requirement of accurate propagation over large distances necessitating good approximation of dispersion characteristics of the governing PDE, it is evident that optimized FD schemes should generally be chosen for seismic applications.
	\\
	\\It is worth considering that seismic surveys are usually designed in accordance with available computational resources, such that high-frequency components liable to result in numerical dispersion during processing are excluded (\citealp{Dablain1986}). Thus in practice, spatial optimization may be unnecessary to ensure a good-quality solution in applications such as migration or forward modeling.
	\\
	\\Note that extension to further spatial dimensions is straightforward for both standard and staggered grid formulations (\citealp{Moczo2000}). Beyond reduction of the formal accuracy order of the stencil, no clear disadvantages of spatial optimization were highlighted during testing.

	\subsection{Comparative performance and grid spacing}
	In the 1D acoustic test case, it was expected that error in both conventional and optimized schemes would be near-identical at the maximum grid spacing for acceptable propagation accuracy using the conventional scheme: $\Delta x= 0.0213$. The actual point of convergence was found to reside at an intermediate value between maximum grid spacings for the conventional and optimized schemes. This discrepancy resulted from interplay between truncation and numerical dispersion error. 
	\\
	\\Whilst studies have presented methods to quantify the truncation error in FD approximations (e.g. \citealp{Lantz1971}, \citealp{Warming1974}), comparisons of dispersion and truncation error are scarce. Quantification of total error proportion contributed by each is thus difficult. Due to the tedious, time-consuming nature of existing error analysis methods (\citealp{Warming1974}), error is typically determined empirically for some test case, evaluated relative to exact solutions or existing FD methods (e.g. \citealp{Finkelstein2007}, \citealp{Zhang2013}). Additionally, whilst quantification of resolving power is commonplace in studies of DRP schemes (e.g. \citealp{Tam1993}, \citealp{Liu2014}), these methods are not directly comparable with truncation error investigations, evaluating accuracy via different metrics.
	\\
	\\Meaningful evaluation of respective contributions of dispersion and truncation error regarding the point of equal performance between the two schemes is consequentially near-impossible. However, test case results imply that conventional FD stencils should be considered when $\Delta x$ is greater than the limit for good approximation of the dispersion relationship (so long as it is smaller than that for the optimized scheme). In practice, testing both schemes for $\Delta x$ within these bounds and selecting the scheme which yields superior results is likely prudent. This approach is supported by the results of the 2D elastic test, demonstrating that for some problems, the optimized scheme will yield reduced error at maximum $\Delta x$ for the conventional. Additionally, optimized schemes are found to achieve more accurate solutions under these circumstances in test cases described in \citealp{Moczo2000} and \citealp{Finkelstein2007}, supporting the notion that ideal scheme choice is somewhat dependent upon problem specification. Further investigation is required to develop a more robust understanding however.

	\subsection{Optimized scheme behavior for oversampled wavefields}
	The superior performance of the optimized scheme at very small values of $\Delta x$ observed in the first test case is inconsistent with expectations based upon resolution characteristics and truncation error for the two schemes. At $\Delta x=0.0100\mathrm{m}$, sampling of the highest frequency component of the standing wave is well in excess of the number of grid spacings per wavelength required for preservation of propagation characteristics, producing minimal numerical dispersion for either scheme. It follows that higher-order truncation of the conventional stencil should result in reduced error under these circumstances. However, errors observed for both schemes between the start of the simulation and 10s were near-identical, with a reduction in error accumulation rate after this point, making apparent another peculiar aspect of the observed behavior of the optimized scheme: error in the numerical solution appears not to increase after this point.
	\begin{figure*}[h!]
		\begin{center}
			\includegraphics[angle=0, width=0.75\textwidth, keepaspectratio=true]{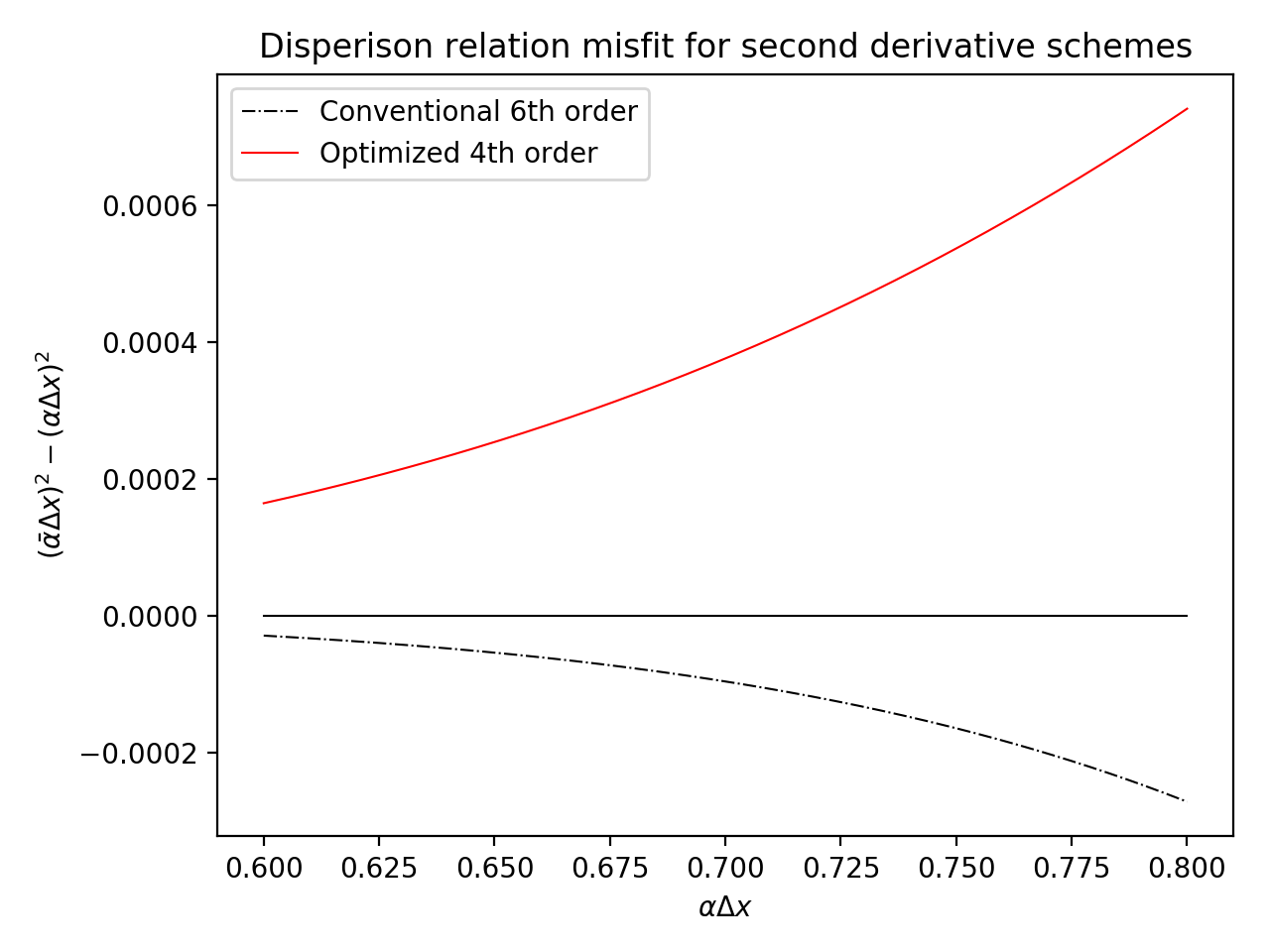} 
			\caption{$(\bar{\alpha}\Delta x)^{2}-(\alpha\Delta x)^{2}$ versus $\alpha\Delta x$ for the 4\textsuperscript{th} order optimized and 6\textsuperscript{th} order conventional central-difference schemes. The solid black line designates the ideal.}
			\label{2nd_dispersion_plot_difference}
		\end{center}
	\end{figure*}
	\\
	\\It is established that FD schemes will still produce some small degree of numerical dispersion, even for particularly fine grid spacings (\citealp{Wang2003}). Additionally, it has been demonstrated that DRP schemes can still marginally reduce numerical dispersion compared to conventional FD schemes for heavily oversampled wavefields (\citealp{Wang2003}). Considering the magnitude of the discrepancy between the two schemes observed in the 1D acoustic wave test case for small $\Delta x$, this presents a possible mechanism by which the optimized scheme could generate the smaller error. However, when the relationship between $\bar{\alpha}\Delta x$ and $\alpha\Delta x$ is examined against the ideal over the region of interest (see figure \ref{2nd_dispersion_plot_difference}), it is clear that not only does the conventional scheme better approximate the second derivative, but the magnitude of misfit present is entirely inconsistent with that of the errors observed. Consequentially, it is unlikely that this is behind behavior exhibited by the first test case at small $\Delta x$.
	\begin{figure*}[h!]
		\begin{center}
			\includegraphics[angle=0, width=0.75\textwidth, keepaspectratio=true]{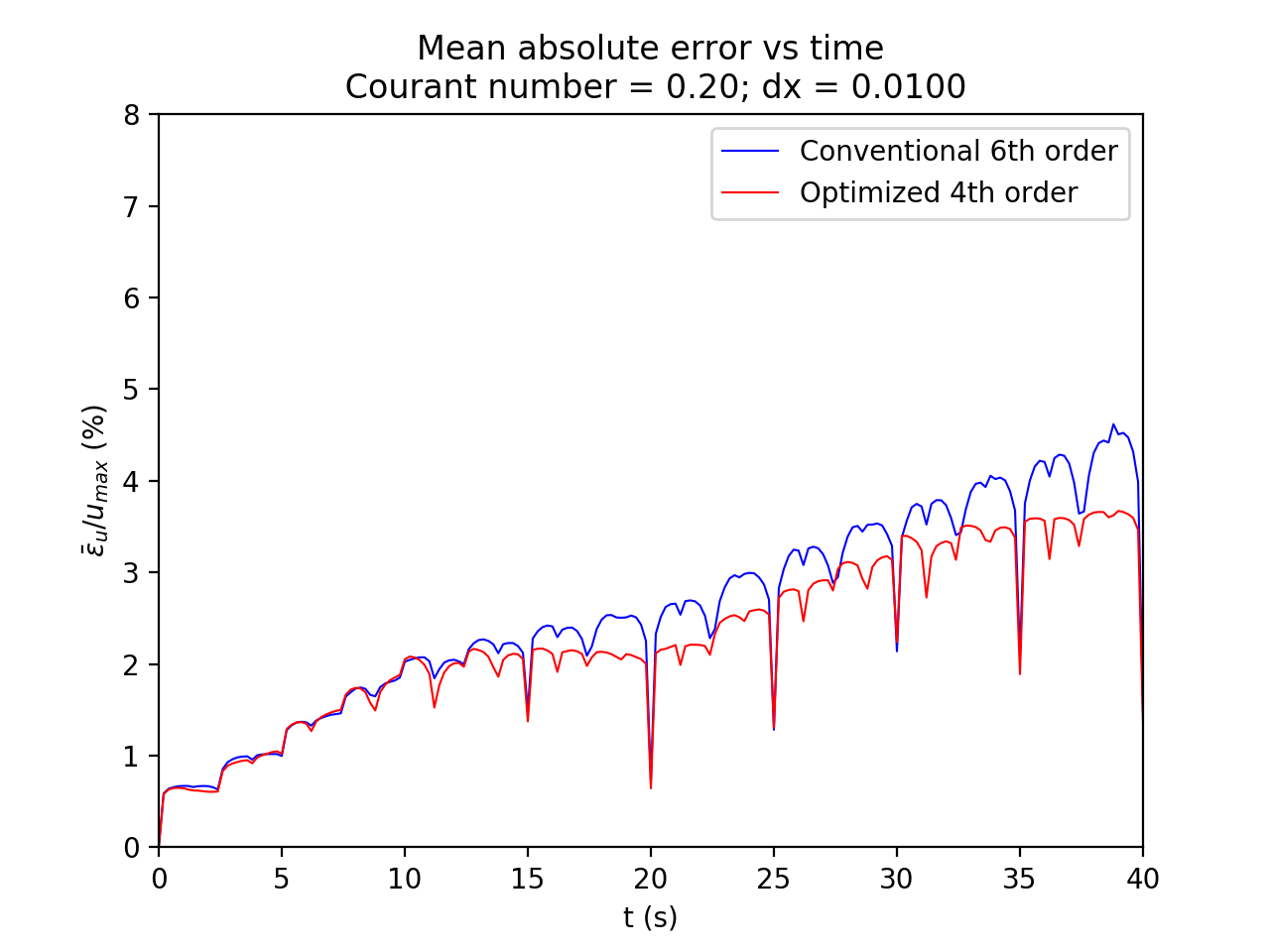} 
			\caption{Evolution of mean absolute error over time for the 1D acoustic wave test case using $\Delta x=0.0100\mathrm{m}$. The simulation has been extended to t=40s. Note that the reduced rate of error introduction apparent in figure \ref{percent_mean_err_over} is only temporary, and normal buildup of error resumes after 20s.}
			\label{percent_mean_err_time_long}
		\end{center}
	\end{figure*}
	\\
	\\An alternative proposition is that the apparent reduction in the rate of change of error with time merely results from the chosen runtime, and expected behavior resumes at some later time outside the window studied. As the relation between mean error and time is fairly rough throughout, this is likely not an unreasonable assumption. To confirm, the first test case for $\Delta x =0.0100$ was run to 40s as opposed to the 20s initially used. From the extended plot of mean absolute error against time in figure \ref{percent_mean_err_time_long}, it is clear that the apparent pause in error accumulation ends at approximately 20s. However, the unexplained superior performance of the optimized scheme continues until simulation end.
	\\
	\\In the absence of a rigorous explanation for the superior performance of the optimized scheme for fine grid spacings, the observed behavior is assumed to result from some aspect of problem specification for the 1D acoustic test case. Similar behavior is reported for a comparable scheme in \citealp{Tam1994}, although its cause is not speculated on. Without further testing and error analysis of the FD schemes (e.g. \citealp{Warming1974}), the cause of these observations remains unclear.
	
	\subsection{Limits to optimization benefits for staggered schemes}
	Dispersion characteristics of the conventional staggered scheme (see figure \ref{1st_dispersion_plot_staggered}) are consistent with values published for comparable schemes in \citealp{Saenger2000} and \citealp{Moczo2000}. Whilst the relationship between $\bar{\alpha}\Delta x$ and $\alpha\Delta x$ deviates from that of the underlying partial differential at a comparatively small value of approximately $\alpha\Delta x=0.5$, it is apparent that the minimizations detailed in equations \ref{eq:25} and \ref{eq:26} yield little further improvement. The 6\textsuperscript{th} order conventional staggered scheme for first derivatives can thus be considered near-optimized, posing the question as to whether this is inherent in staggered schemes, or merely coincidental.
	\\
	\\The least-squares method of coefficient optimization for staggered grids detailed by \citealp{Liu2014} yields identical coefficients to those derived for the 2D elastic test case. In assuming the equivalence of these methods, it follows that similarity between coefficients of the standard and optimized staggered schemes should extend only to stencils of small extent. As approximation order is increased, quality of the optimized FD solution differs increasingly from the conventional of equivalent stencil extent. The apparent limitations of spatial optimization for staggered grids in the derivation of stencils for the 2D elastic test case are therefore a product of stencil extent rather than an ubiquitous feature of staggered grid formulations. As high-order schemes are commonplace in seismic applications (e.g. \citealp{Hong-wei2010}), optimized schemes should typically offer considerable advantages over conventional stencils in practice.

	\section{Conclusions, evaluation, and further work}
	In this work, the spatial optimization method for developing DRP schemes proposed by \citealp{Tam1993} has been extended for use with second derivatives and staggered grids. The optimized second derivative stencil was found to generate reduced error versus the conventional stencil of equivalent extent for wavefields with components of wavelength smaller than $4.4\Delta x$: somewhat larger than the theoretical value of $4.1\Delta x$. This discrepancy is postulated to result from the effects of truncation error. For lower frequency waves or finer grids, both schemes were capable of accurately resolving wave propagation characteristics, resulting in the conventional stencil yielding reduced error due to higher-order Taylor series truncation. However, for heavily oversampled wavefields, the optimized scheme was found to generate less error than the conventional. Whilst reports of this behavior exist in the literature, its causes remain unclear, although it is speculated to result from formulation of the test case or its parameters. It is proposed that notable error reduction at coarser grid spacings yielded by optimization be exploited to minimize computational cost, potentially augmented by variable temporal and spatial increments. This would yield considerable practical benefits in a suite of modeling applications.
	\\
	\\The staggered optimized stencil for first derivatives was found to moderately improve numerical dispersion suppression versus the conventional stencil for wave components with wavelengths at the theoretical minimum for both schemes ($12.6\Delta x$). This was most pronounced in the $\uptau_{xz}$ field, consistent with shorter wavelengths of slower S-wave components. Whilst optimized and conventional staggered stencils of 4\textsuperscript{th} and 6\textsuperscript{th} order respectively have similar coefficients and performance, benefits of spatial optimization are more pronounced for higher-order staggered stencils, implying suitability for seismic modeling and processing applications where high-order schemes are favored. Successful application of spatially-optimized schemes to the P-SV elastic wave formulation demonstrates viability of the method as a means of suppressing numerical dispersion in varied contexts. Consideration of spatially-optimized stencils is consequentially recommended for all applications featuring wavelengths near or below minimum values for accurate propagation.
	\\
	\\Whilst the 1D example presented in the first test case is not dissimilar to several published examples, establishment of the method's validity benefited from a simple initial test case. Additionally, despite usage of spatially-optimized second-derivative stencils detailed in several publications, derivation of coefficients is rarely elaborated on. Whilst prevented by time and space constraints, further testing with dispersive 1D examples may have aided better understanding of the level of numerical dispersion reduction offered by the optimized scheme, and the conditions under which it is most pronounced. The method detailed for spatial optimization of staggered grid stencils represents a considerably more straightforward means than existing methods, and has been demonstrated to provide tangible improvements when applied to simple seismic modeling applications. Without an exact solution, improvements in propagation accuracy for the optimized scheme are subject to a degree of uncertainty. However, as comparison to other FD schemes for performance evaluation is commonplace in published studies, this is unlikely cause for concern.
	\\
	\\The derived spatially-optimized stencils would be complemented by temporal optimization to produce a full DRP scheme for staggered grid applications, further improving propagation accuracy of high-frequency components. It may also be of interest to investigate feasibility of stencil truncation upon accuracy of high-order optimized spatial derivatives, as this could potentially present means to considerably increase the accuracy of an FD scheme without additional computational cost. Viable avenues of inquiry also include integration of spatially-optimized FD methods into forward propagators for full-waveform inversion or migration algorithms, investigating potential benefits of increased propagation accuracy, or reduced runtime enabled by coarser grids. 
	\newpage
	\bibliographystyle{apalike}
	\bibliography{drp_fd_schemes_report_final}
\end{document}